\documentclass [11pt]{article}
\usepackage{amsmath, amsthm, amssymb}
\usepackage{graphicx}
\newfont{\bbb} {msbm10}
\newcommand{\R}{\Bbb{R}}
\newcommand{\bS}{\Bbb{S}}
\newcommand{\bB}{\Bbb{B}}
\newcommand{\T}{\Bbb{T}}
\newcommand{\C}{\Bbb{C}}
\newcommand{\qH}{\Bbb{H}}
\newcommand{\Oc}{\Bbb{O}}

\newcommand{\D}{\Bbb{D}}
\newcommand{\cT}{{\cal{T}}}

\newcommand{\sbs}{\subset}
\newcommand{\ra}{\rightarrow}

\newcommand{\bg}{{\bar{g}}}

\newcommand{\dsigma}{{\dot{\sigma}}}

\newcommand{\hg}{{\hat{g}}}
\newcommand{\hh}{{\hat{h}}}

\newcommand{\cF}{{\cal{F}}}

\newcommand{\cZ}{{\cal{Z}}}
\newcommand{\cY}{{\cal{Y}}}
\newcommand{\cX}{{\cal{X}}}
\newcommand{\ccP}{{\mbox{{\Large $\wp$}}}}

\newcommand{\p}{\partial}

\newcommand{\cL}{{\cal{L}}}

\newcommand{\cS}{{\cal{S}}}
\newcommand{\cH}{{\cal{H}}}

\newcommand{\cG}{{\cal{G}}}
\newcommand{\cA}{{\cal{A}}}
\newcommand{\cE}{{\cal{E}}}

\newcommand{\cW}{{\cal{W}}}

\newcommand{\B}{\Bbb{B}}

\newcommand{\HH}{\Bbb{H}}

\newcommand{\ssl}{{_\lambda}}
\newcommand{\sd}[1]{\Delta_{_{\bS^{#1}}}}

\newcommand{\s}[1]{{\sf{#1}}}

\newcommand{\0}[1]{_{_{#1}}}
\newcommand{\rC}{{\rm{C}}\,}

\newcommand{\sN}{{\sf{N}}}
\newcommand{\sL}{{\sf{Link}}}

\newcommand{\sB}{{\sf{B}}}
\newcommand{\sA}{{\sf{A}}}

\newcommand{\dX}{\dot{X}}
\newcommand{\dDelta}{\dot{\Delta}}
\newcommand{\dsquare}{\dot{\square}}
\usepackage[all]{xy}

\begin{document}

\title{Riemannian Hyperbolization}
\author{Pedro Ontaneda\thanks{The author was
partially supported by a NSF grant.}}
\date{}

\maketitle

Classical flat geometry has formed part of basic human knowledge since ancient times. It is characterized by the almost universally known condition that the sum of the internal angles of a triangle $\triangle$ is equal to $\pi$. We write $\Sigma(\triangle)=\pi$.
Other fundamental geometries are defined by replacing the equality
$\Sigma(\triangle)=\pi$ by inequalities; thus positively curved geometries and
negatively curved geometries are determined by the inequalities
 $\Sigma(\triangle)>\pi$ and $\Sigma(\triangle)<\pi$, respectively, where
$\triangle$ runs over all small non-degenerate triangles in a space.
It is natural then to try to find spaces that admit such geometries, and this
task has been a driving force in Riemannian Geometry for many decades.
But surprisingly there are not too many examples of smooth closed manifolds
that support either a positively curved or a negatively curved metric.
For instance, besides spheres, in dimensions $\geq 17$ (and $\neq 24$) the only positively curved simply connected known
examples are complex and quaternionic projective spaces. In negative
curvature the situation is arguably more striking because negative curvature
has been studied extensively in many different areas in mathematics.
Indeed,  from the ergodicity of their geodesic flow in
Dynamical Systems to their topological rigidity in Geometric Topology;
from the existence of harmonic maps in Geometric Analysis to the
well-studied and greatly generalized algebraic properties of their fundamental groups,
negatively curved Riemannian manifolds are the main object in many 
important and well-known results in mathematics. Yet the fact remains that very
few examples of closed negatively curved Riemannian manifolds are known.
Besides the hyperbolic ones ($\R$, $\C$, $\qH$, $\Oc$), the other known examples are  the Mostow-Siu examples (complex dimension 2) which are local branched covers of complex hyperbolic space (1980, \cite{MS}), 
the Gromov-Thurston examples (1987, \cite{GT}) which are branched covers of  real hyperbolic
ones, the exotic Farrell-Jones examples (1989, \cite{FJ1}) which are homeomorphic but not diffeomorphic to real hyperbolic manifolds
(and there are other examples of exotic type), and the three
examples of Deraux (2005, \cite{Deraux}) which are of the Mostow-Siu
type in complex dimension 3. Hence, excluding the Mostow-Siu and
Deraux examples (in dimensions 4 and 6, respectively), all known examples
of closed negatively curved Riemannian manifolds are homeomorphic to
either a hyperbolic one or a branched cover of a hyperbolic one.\vspace{.1in}

%\noindent {\bf Remark.} The Mostow-Siu and Deraux examples have
%covers that have  maps
%to complex hyperbolic space that look locally like branched covers.
%It is not known whether all these examples are global branched covers.\vspace{.1in}

This lack of examples in negative curvature changes dramatically if we allow singularities,
and a very rich and abundant class of negatively curved spaces (in the geodesic sense)
exists due to the strict hyperbolization process of Charney and Davis \cite{ChD}. The
hyperbolization process was originally introduced by M. Gromov \cite{G}, and later
studied by Davis and  Januszkiewicz \cite{DJ}, and Charney-Davis strict hyperbolization is built on these previous versions. The hyperbolization process is conceptually (but not technically) quite
simple since it has a lego type flavor: in the same way as simplicial complexes
and cubical complexes are built from a basic set of pieces,
basic ``hyperbolization pieces" are chosen, and anything that can be built or assembled 
with these pieces will be negatively curved. This conceptual simplicity could be in some sense a bit deceptive because hyperbolization produces an enormous class of examples
with a very fertile set of properties. But the richness and complexity of the hyperbolized
objects are matched by the richness and complexity of the singularities obtained,
and hyperbolized smooth manifolds are very far from being Riemannian.
Interestingly one can relax and lose even
more regularity and consider negative curvature from the algebraic point of
view, that is consider Gromov's hyperbolic groups, and it can be argued
\cite{Ol} that ``almost every group" is hyperbolic. So, negative curvature
is in some weak sense generic, but Riemannian negative curvature seems very scarce.
It is natural then to inquire about the difference between the class of manifolds
with negatively curved metrics with singularities and its subclass of more
regular Riemannian counterparts. More specifically we can ask whether
the strict hyperbolization process can be brought into the Riemannian universe. 
In this paper we give a positive answer to this question, and we do this
by proving that all singularities of the Charney-Davis strict hyperbolization of
a closed smooth manifold can be smoothed, provided the ``hyperbolization
piece" is large enough (which can always be done). Moreover we prove that we can do this process in a $\epsilon$-pinched way. Here is the statement of our Main Theorem.\vspace{.1in}

\noindent {\bf Main Theorem.} {\it Let $M^n$ be a closed smooth manifold
and let $\epsilon>0$. Then there is a closed Riemannian manifold $N^n$ and
a smooth map $f:N\ra M$ such that}
\begin{enumerate}
\item[{\it (i)}] {\it The Riemannian manifold $N$ has sectional curvatures in the interval
$[-1-\epsilon,-1]$.}
\item[{\it (ii)}] {\it The induced map  $f_*:H_*(N,R)\ra H_*(M,R)$
is surjective, for every (untwisted) $R$.}
\item[{\it (iii)}] {\it If $M$ is $R$-orientable then $N$ is $R$-orientable. In this case
(ii) implies that $f$ has degree one and
$f^*:H^*(M,R)\ra H^*(N,R)$ is injective.}
\item[{\it (iv)}] {\it The map $f^*$ sends the rational Pontryagin classes of $M$ to the
rational Pontryagin classes of $N$. }
%\item[{\it (i)}]
\end{enumerate}
%\vspace{.1in}

\noindent {\bf Addendum to Main Theorem.} {\it The manifold $N$ is
the Charney-Davis strict hyperbolization of $M$ but with a different smooth structure. The hyperbolization is done with a sufficiently
``large" hyperbolization piece $X$.}\vspace{.1in}

By ``large" above we mean that the width of the normal neighborhoods of the faces
of $X$ are very large. These large pieces always exist (see 9.1).
Corollaries 1, 2 and 3 below are the $\epsilon$-pinched Riemannian versions of classical applications
of hyperbolization.\vspace{.1in}

\noindent {\bf Corollary 1.} {\it Every closed smooth manifold is smoothly cobordant
to a closed Riemannian manifold with sectional curvatures in the interval
$[-1-\epsilon,-1]$, for every $\epsilon>0$.}\vspace{.1in}

\noindent {\bf Corollary 2.} {\it The cohomology ring of any finite $CW$-complex 
embeds in the cohomology ring of closed  Riemannian manifold
with sectional curvatures in the interval
$[-1-\epsilon,-1]$, for every $\epsilon>0$.}\vspace{.1in}

\noindent {\bf Proof.} Let $X$ be a finite $CW$-complex. Embed $X$ in some
$\R^n$ and let $P$ be a compact neighborhood of $X$ that retracts to $X$. Let $M$ be the double of $P$. Then there is a retraction $M\ra X$, and Corollary 2 follows from (iii) in the Main Theorem.\vspace{.1in}

Since degree one maps between closed orientable manifolds are $\pi_1$-
surjective we obtain the following result.\vspace{.1in} 

\noindent {\bf Corollary 3.} {\it For every finite CW-complex $X$
there is a closed Riemannian manifold $N$ and a map $f:N\ra X$ such
that: (i) $N$ has sectional curvatures in the interval
$[-1-\epsilon,-1]$, (ii) $f$ is $\pi_1$-surjecive, (iii) $f$ is homology
surjective.}\vspace{.1in}

All known examples of closed negatively curved Riemannian manifolds
with less than $\frac{1}{4}$-pinched curvature have zero rational Pontryagin classes
(for the Gromov-Thurston branched cover examples this was proved by
S. Ardanza \cite{A}). 
The next corollary gives examples of such manifolds with
nonzero rational Pontryagin classes. \vspace{.1in}

\noindent {\bf Corollary 4.} {\it For every $\epsilon>0$ and $n\geq 4$ there
is a closed Riemannian $n$-manifold with sectional curvatures in the interval
$[-1-\epsilon,-1]$ and nonzero rational Pontryagin classes.}\vspace{.1in}

\noindent {\bf Proof.} Take $M$ in the Main Theorem orientable with nonzero
Pontryagin classes.\vspace{.1in}

All manifolds given in Corollary 4 are new examples of closed negatively
curved manifolds.
This follows from Novikov's topological invariance of the rational Pontryagin classes \cite{N}, and the $\frac{1}{4}$-pinched rigidity results given in
(or implied by)
the work
of Hern\'andez \cite{Her}, Yau and Zheng \cite{YZ}, Corlette \cite{Cor},
Gromov \cite{G3} and Mok-Siu-Yeung \cite{MSY}. We state this
in the next corollary.\vspace{.1in}

\noindent {\bf Corollary 5.} {\it For any $\epsilon >0$ and $n\geq 4$ there are
closed Riemannian $n$-manifolds with sectional curvatures in the interval
$[-1-\epsilon,-1]$ that are neither homeomorphic to a hyperbolic manifold
($\R$, $\C$, $\qH$, $\Oc$) nor homeomorphic to the Gromov-Thurston branched cover of a real hyperbolic one, nor homeomorphic to one of the Mostow-Siu or Deraux examples.}\vspace{.1in}

The next application was suggested to us by Stratos Prassidis some time ago and deals with cusps of negatively curved manifolds.
Recall that if $M$ is a complete finite volume noncompact real hyperbolic manifold
then $M$ has finitely many cusps isometric to manifolds of the form $Q\times [b,\infty)$
with metric $e^{-2t}h+dt^2$, where $(Q, h)$ is a closed flat manifold
and $b\in \R$. If $M$ has exactly one cusp diffeomorphic to
$Q\times [b,\infty)$
 we say that {\it the manifold
$Q$ bounds geometrically a hyperbolic manifold.}\vspace{.1in}

More generally, in 1978 M. Gromov defined almost flat manifolds in \cite{G2} and
similar facts hold for them replacing hyperbolic manifolds by pinched
negatively curved manifolds. That is,
let $M$ be a complete finite volume noncompact manifold with pinched negative 
curvature (i.e all sectional curvatures lie in a fixed interval $[-a,-b]$, $0<b\leq a<\infty$).
Then $M$ has finitely many cusps diffeomorphic to manifolds of the form $Q\times [b,\infty)$, 
where $Q$ is an almost flat manifold. 
If $M$ has exactly one cusp diffeomorphic to 
 $Q\times [b,\infty)$
we say that {\it the manifold
$Q$ bounds geometrically a negatively curved manifold.}
Of course a necessary condition for $Q$ to bound geometrically as above is to
smoothly bound a compact manifold. \vspace{.1in}

\noindent {\bf Remark.} Here we do not assume $Q$ to be connected.
Hence if $Q$ bounds geometrically then the number of connected components of $Q$ is the same as the number of connected cusps.\vspace{.1in}

It was proved by Hamrick and Royster \cite{HR} that every closed flat manifold bounds smoothly. This together with the work of Gromov in \cite{G1}, \cite{G2} motivated
Tom Farrell and Smilka Zdravkovska to make the following well-known 
conjectures in \cite{FZ}  30 years ago.\vspace{.1in}

\noindent {\bf Conjecture 1.} Every closed almost flat manifold bounds smoothly. This conjecture was also proposed, independently, by S.-T. Yau in \cite{Y}.\vspace{.1in}

\noindent {\bf Conjecture 2.} Every closed
flat manifold bounds geometrically a hyperbolic manifold.\vspace{.1in}

\noindent {\bf Conjecture 3.} Every closed almost flat manifold bounds geometrically a negatively curved manifold.\vspace{.1in}

It was showed by Long and Reid \cite{LR} that Conjecture 2 is false
by giving examples of three dimensional flat manifolds that do not bound.
The following result says Conjecture 1 implies Conjecture 3.
\vspace{.0001in}

\noindent{\bf Theorem A.} {\it Let $Q$ be a closed almost flat manifold. Assume that
$Q$ bounds smoothly. Then $Q$ bounds geometrically a negatively curved manifold $M$.}\vspace{.1in}

Conjecture 1 has generated a lot of research in the last 30 years and it is
known to be true for an almost flat manifold in the following cases.
Let $Q$ be almost flat. Then $Q$ is covered by a nilmanifold, that is,
the quotient of a simply connected nilpotent Lie group $L$ by a uniform lattice. Denote by $G$ the holonomy of $Q$.\vspace{.1in}

\noindent\,\,\,\,(a)\, The manifold $Q$ is a nilmanifold.\vspace{.05in}

\noindent\,\,\,\,(b)\, The holonomy $G$ has order $k$ or $2k$, where $k$ is odd,
due to Farrell-Zdravkovska \cite{FZ}.\vspace{.05in}

\noindent\,\,\,\,(c)\, The holonomy $G$ of $Q$ acts effectively on the center of $L$, also due to
Farrell-Zdravkovska \cite{FZ}.\vspace{.05in}

\noindent\,\,\,\,(d)\, The holonomy $G$ is cyclic, due to J. Davis and F. Fang \cite{DF}.
Also  Upadhyay \cite{U} had proved that

\hspace{.1in} Conjecture 1 it true when
the following conditions hold: $G$ is cyclic, $G$ acts trivially on the center 

\hspace{.1in} of $L$, and $L$ is 2-step nilpotent. \vspace{.1in}

Hence in all of the above cases $Q$ bounds geometrically a pinched
negatively curved manifold. Note that for any closed $Q$ we have
$\p ( Q\times I)=Q\coprod Q$. Thus we get the following corollary of Theorem A.\vspace{.1in}

\noindent {\bf Corollary 6.} {\it Let $Q$ be a closed almost flat manifold.
Then $Q\coprod Q$ bounds geometrically a pinched negatively curved
manifold.}\vspace{.1in} 

In other words, for every closed connected almost flat manifold there is
a complete finite volume pinched negatively curved manifold with
exactly {\sf two} connected cusps, each diffeomorphic to $Q\times [b,\infty)$.\vspace{.1in}

A complete pinched negatively curved metric $g$ on $Q\times \R$ is called a {\it (pinched negatively curved) cusp metric} if the $g$-volume of $Q\times [0,\infty)$
is finite. And we say that a cusp
metric $g$ on $Q\times \R$ {\it is an eventually warped cusp metric} if $g=e^{-2t}h+dt^2$, for $t<c$, for some $c\in \R$ and a metric $h$ on $Q$.
I. Belegradek and V. Kapovitch \cite{BK} (see also \cite{B}) show, based on earlier work by Z.M. Shen \cite{Shen}, that if $Q$ is almost flat then $Q\times\R$ admits an eventually warped cusp metric.\vspace{.1in}

\noindent {\bf Addendum to Theorem A.} {\it 
Let $g$ be an eventually warped cusp metric on $Q\times\R$. If the sectional curvatures
of $g$ lie in $(a,b)$, with $a<-1<b$, then we can take $M$ in Theorem A with sectional curvatures also in $(a,b)$. Moreover the sectional curvatures of $M$ away from a cusp can be taken in $[-\epsilon-1,-1]$, for any $\epsilon>0$.}\vspace{.1in}

Even though a flat manifold  may not necessarily bound geometrically a hyperbolic
manifold the next corollary says it does bound geometrically an $\epsilon$-pinched to -1
manifold, for any $\epsilon>0$. It follows from the Hamrick and Royster result \cite{HR},
Theorem A and its addendum.\vspace{.1in}

\noindent {\bf Corollary 7.} {\it Every closed flat manifold
bounds geometrically a manifold with sectional curvatures in
$[-\epsilon-1,-1]$, for any $\epsilon>0$.}\vspace{.1in}

Here is a brief description of the paper. In Section 1 
we introduce some notation and basic concepts, including
the definition of $\epsilon$-close to hyperbolic metrics.
This is a slightly technical but important concept. The idea is to
try to measure how close a metric is to being hyperbolic;
we do this in a chart by chart fashion. In Section 2 we
define and study the ``hyperbolic extension" of a metric (or space), which is
a key geometric construction. In this section there are no proofs and
we essentially collect the main results of \cite{O5}.
In Section 3 we describe another key geometric construction,
hyperbolic forcing; it is the composition of two
deformations: warp forcing and the two-variable deformation,
which are studied with more detail in \cite{O4} and \cite{O3},
respectively. Section 4 is a family version of Section 3.
Again, in sections 3 and 4 there are essentially no proofs and
we mostly collect the main results of \cite{O4}, \cite{O3}, and
\cite{O6}. In Section 5 we study neighborhoods of 
simplices of all-right spherical complexes. In this section we introduce a technical device that we called {\it sets of widths}. These are sets of positive real numbers
 that are used as widths for normal neighborhoods of simplices of all-right spherical complexes.
We prove that there are sets of widths, independent of the complex, that
satisfy very useful properties. These are fundamental objects that make all matching processes work. Section 6 is a sort of a ``cone version"
of Section 5; in it we study (all-right) piecewise hyperbolic cone complexes, which are just cones over all-right spherical
complexes with metric warped by $\sinh$. In Section 7 we deal with the smoothing issue for cubical and
all-right spherical complexes; here we collect the
main concepts and results of \cite{O1}. We put everything together
in Section 8 to smooth hyperbolic cones. Section 9 is dedicated to the Charney-Davis
strict hyperbolization process; in this section we collect the results in
\cite{O2}, in particular we mention that strictly hyperbolized smooth manifolds have ``normal differentiable structures". 
Finally we prove the Main Theorem in
Section 10 and Theorem A in Section 11. Subsections have been added
at the end of sections 7, 8, and 9, that deal with generalizations to the case
of manifolds with codimension zero singularities. These subsections
are used in Section 11. \vspace{.1in}

We are grateful to C.S. Aravinda,  I. Belegradek, M. Deraux, 
F.T. Farrell, R. Geoghegan, and
J.-F. Lafont  for their comments
and suggestions. \vspace{.1in}

\begin{center} {\bf \large  Section 1. Some Notation,
Definitions, and Metrics $\epsilon$-Close to Hyperbolic.}\end{center} 

In this paper $\rho$ will denote a fixed smooth function $\rho:\R\ra[0,1]$ such that: {\bf (i)}\, $\rho|_{(-\infty, 0+\delta]}\equiv 0$, and\, \,{\bf (ii)}\, $\rho|_{[1-\delta\, ,\infty)}\equiv 1$, where $\delta>0$ is small.\vspace{.1in}

The standard flat metric on $\R^l$ will be denoted by
$\sigma\0{\R^{l}}$. Similarly, $\sigma\0{\HH^{l}}$ and
$\sigma\0{\bS^{l-1}}$ will denote the standard hyperbolic
and round metrics on $\HH^l$ and $\bS^{l-1}$,
respectively.\vspace{.1in}

Let $\bB=\bB^{l-1}\sbs\R^{l-1} $ be the unit ball, 
with metric $\sigma\0{\R^{l-1}}$.  Write  $I_\xi=(-(1+\xi),1+\xi)\sbs\R$, $\xi>0$.
Our basic models are $\T^{l}_\xi=\T_\xi=\bB\times I_\xi\sbs \R^l$,  with hyperbolic metric $\sigma=e^{2t}\sigma\0{\R^{l-1}}+dt^2$. 
In what follows we may sometimes suppress the subindex $\xi$, if the context is clear.
The number $\xi$ is the {\it excess} of $\T_\xi$. \vspace{.05in}

\noindent {\bf Remarks.}

\noindent {\bf 1.}  One of the reasons to introduce the excess  is that the process of 
hyperbolic extension (see Section 2) decreases the excess of the charts, as shown in the statement of Theorem 2.6.

\noindent {\bf 2.} In the applications we may actually need warp metrics with warping functions that are multiples of hyperbolic
functions. All these functions are close to the exponential $e^t$ (for $t$ large), so instead of introducing one model for each hyperbolic function
we introduced only the exponential model.\vspace{.1in}

Let $|.|$ denote the uniform $C^2$-norm of $\R^s$-valued functions on $\T_\xi=\bB\times I_\xi\sbs\R^{l}$.
Given a metric $g$ on $\T$, $|g|$ is computed considering $g$ as the $\R^{l^2}$-valued function $(x,t)\mapsto (g_{ij}(x,t))$,
$(x,t)\in\T$, where, as usual,
$g_{ij}=g(e_i,e_j)$, and the $e_i$'s are the canonical vectors in $\R^{l}$. Let $\epsilon>0$. We will say that a metric $g$ on $\T$ is $\epsilon$-{\it close hyperbolic}
if $|g-\sigma|<\epsilon$.\vspace{.1in}

Let $\epsilon >0$. A Riemannian manifold $(M^l,g)$ is $\epsilon$-{\it close hyperbolic} if there is $\xi>0$ such that for every $p\in M$ there is an $\epsilon$-close to hyperbolic 
chart with center $p$, that is, there is a chart
$\phi :\T_\xi\ra M$, $\phi(0,0)=p$,  such that $\phi^*g$ is $\epsilon$-close to hyperbolic. Note that all charts are defined on the same model space
$\T_\xi$. The number $\xi$ is called the {\it excess } of the charts
(which  is fixed).
More generally, a subset $S\sbs M$ is $\epsilon$-close to hyperbolic if every $p\in S$ is the center of an $\epsilon$-close to hyperbolic chart in $M$ with fixed
excess $\xi$. \vspace{.1in}

\noindent {\bf Remark 1.1.} Note that if a metric is
$\epsilon$-close to hyperbolic with charts of excess $\xi$ then
it is $\epsilon$-close to hyperbolic with charts of excess $\xi'$,
with $0<\xi'\leq \xi$. \vspace{.1in}

Let $M^n$ be a complete Riemannian manifold. We say that
a point $o\in M$ is a {\it center of} $M$ if the exponential map $exp_o:T_oM\ra M$ is a diffeomorphism.
In particular $M$ is diffeomorphic to $\R^n$. 
For instance if $M$ is Hadamard manifold every point is a center point. In this paper we will always use the same symbol ``$o$" to denote a center of a Riemannian manifold, unless it is necessary to
specify the manifold $M$, in which case we will write $o_{_M}$.
Using the diffeomorphism $exp_o$ onto $M$ and an identification
of $T_oM$ with
$\R^n$ via some fixed choice of an orthonormal basis in $T_oM$,
we can identify $M$ with $\R^n$ and $M-\{o\}$ with $\bS^{n-1}\times\R^+$.
Therefore
 the metric of $M$, restricted to $M-\{ o\}$, can be written 
as  $h_t+dt^2$ on $\bS^{n-1}\times \R^+$, where $\{h_t\}\0{t>0}$
is a one-parameter family of metrics  on $\bS^{n-1}$.
The set of rays
$t\mapsto (x,t)\in \bS^{n-1}\times \R^+$ are geodesics on
$M-\{o\}=\bS^{n-1}\times \R^+$, and we call this set  {\it the ray structure of $M$ with respect to o}.\vspace{.1in}

If $N^l$ has center $o$ we say that $S\sbs N$ is {\it radially $\epsilon$-close to hyperbolic (with respect to $o$)} if it $\epsilon$-close to hyperbolic, and the  
$\epsilon$-close to hyperbolic charts $\phi$ respect the product structure of $\T$
and $N-\{o\}=\bS^{l-1}\times\R^+$, that is $\phi(. , t)=(\phi_1(.), t+a)$, for some $a$ depending on the $\phi$. The term ``radially" in the definition above refers to the
 decomposition of the manifold $N-\{o\}$ as a product $\bS^{l-1}\times\R^+$.
\vspace{.1in}

\noindent {\bf Remarks 1.2.} 

\noindent {\bf 1.} The definition of radially $\epsilon$-close to hyperbolic metrics 
is well suited to studying metrics of the form $g\0{t}+dt^2$ for $t$ large, but for
small $t$ this definition is not useful because: (1) we need some space to fit the charts, and (2) the form of our specific fixed model $\T$ (note that the warping function used in the metric of
$\T$ is the exponential).
An undesired consequence is that even punctured hyperbolic space $\HH^n-\{ o\}=\bS^{n-1}\times\R^+$
(with warp metric $sinh^2(t)\sigma\0{\bS^{n-1}}+dt^2$) is not radially $\epsilon$-close to hyperbolic
for $t$ small. In fact there is ${\sf a}={\sf a}(n, \epsilon)$ such that hyperbolic
$n$-space is $\epsilon$-close to hyperbolic for $t>{\sf a}$ (and not for all $t\leq {\sf a}$), see 3.9 \cite{O3}.
This is not essential for what follows.

\noindent {\bf 2.}
For every $n$ there is a function $\epsilon'=\epsilon'(\epsilon,\xi,n)$ such that:  if a Riemannian metric $g$ on a 
manifold $M^{n}$ is $\epsilon'$-close to hyperbolic, with charts of excess $\xi$,  then the sectional curvatures of $g$ all lie $\epsilon$-close to -1.
This choice is possible, and depends only on $n$ and $\xi$, because the curvature depends only of the derivatives up to order 2 of $\phi^*g$ on $\T_\xi$,
where $\phi$ is an $\epsilon$-close to hyperbolic chart with excess $\xi$.\vspace{.1in}

\noindent {\bf Lemma 1.3.} {\it Let $\phi:\T_\xi\ra M$ be a radially 
$\epsilon$-close to hyperbolic chart centered at $p\in M$. Then, for every $q\in \T_\xi$ we have $d\0{M} (\phi(q),p)\,\leq\,  (2 +\xi)+n^2\,\epsilon$.}
\vspace{.05in}

\noindent {\bf Proof.} Write $q=(x_0,t_0)\in\B\times I_\xi$. Consider the path
$\alpha(t)=(t\,x_0,0)$, $t\in [0,1]$, $\beta(t)=(x_0,t\,t_0)$, $t\in [0,1]$, and
$\gamma=\alpha\ast\beta$.
Write $g'=\phi^*g$ and we have $g'=\sigma+h$, with $|h|<\epsilon$. Then the $g$-length $\ell_g(\phi\circ\gamma)$ of $\phi\circ\gamma$ is
$\ell_{g'}(\gamma)=\ell_{g'}(\alpha)+\ell_{g'}(\beta)\leq
\ell_\sigma(\alpha)+\ell_h(\alpha)+(1+\xi)\leq
1+\epsilon n^2+ (1+\xi)$.
Hence $d\0{M} (\phi(q),p)\,\leq\,\ell_{g}(\phi\circ\gamma)\leq\,
(2+ \xi)+n^2\epsilon$. This proves the lemma. \vspace{.1in}

Next we deal with a natural and useful class of metrics.
These are metrics on $\R^n$ (or on a manifold with center)
that are already hyperbolic on the ball
$B_{a}=B_a(0)$ of radius $a$ centered at $0$, and are radially $\epsilon$-close to hyperbolic outside
$B_{a'}$
(here $a'$ is slightly less than $a$). Here is the detailed
definition. Let $M^n$ have center $o$ and
let $B_{a}=B_a(o)$ be the ball in $M$ of radius $a$ centered at $o$. We say
that a metric $h$ on $M$ is $(B_a,\epsilon)$-{\it close to hyperbolic, with charts of excess $\xi$}, \, if 
\begin{enumerate}
\item[  (1)]  On $B_{a}-\{o\}=\bS^{n-1}\times (0,a)$
we have $h=sinh^2(t)\sigma\0{\bS^{n-1}}+dt^2$. Hence $h$
is hyperbolic on $B_a$.
\item[(2)]  the metric $h$ is 
radially $\epsilon$-close to hyperbolic outside $B_{a-1-\xi}$,
with charts of excess $\xi$.
\end{enumerate}

\noindent {\bf Remarks 1.4.}

\noindent {\bf 1.} We have dropped the word ``radially" 
to simplify the notation. But it does appear in condition (2),
where ``radially" refers to the center of $B_a$.

\noindent{\bf 2.} We will always assume $a>{\sf a}+1$,
where {\sf a} is as in 1.2 (1). Therefore conditions (1), (2) 
and Remark 1.2 (1) imply
a stronger version of (2) above:\vspace{.1in}

\noindent \,  (2') the metric $h$ is 
radially $\epsilon$-close to hyperbolic outside $B_{\sf a}$,
with charts of excess $\xi$.\vspace{.1in}

\noindent This is the reason why we demanded radius $a-1-\xi$
in (2), instead of just $a$: we need a ``common annulus" (i.e the closure of $B_a-B_{a-\xi-1}$) where
the metric is both hyperbolic and $\epsilon$-close to hyperbolic. This is important in the proof of Theorem 2.6.
 \vspace{.1in}

\noindent {\bf 3.} Let $\epsilon'$ be as in Remark 1.2.(2).
Then the following is also true:  if a Riemannian metric $g$ on a 
manifold $M^{n}$ is $(B_a,\epsilon')$-{\it close to hyperbolic, with charts of excess $\xi$},  then the sectional curvatures of $g$ all lie $\epsilon$-close to -1. Here there is no restriction on $a$, provided
it is not too small (see Remark 1.4 (2)).\vspace{.1in}

\noindent {\bf 4.} It follows from the definition above
and from 1.1 that if a metric is
$(B_a,\epsilon)$-close to hyperbolic with charts of excess $\xi$ then
it is $(B_a,\epsilon)$-close to hyperbolic with charts of excess $\xi'$,
with $0<\xi'\leq \xi$. \vspace{.1in}

Let {\s c} $>1$.
A metric $g$ on a compact manifold $M$ is {\s c}{\it-bounded} if $|g|<$ {\s c} and $|\, det \,g\,|_{C^0} > 1/${\s c}.
A set of metrics $\{ g\0{\lambda}\}$ on the compact manifold $M$ is {\s c}{\it-bounded} if every $g\0{\lambda}$
is {\s c}-bounded.  \vspace{.1in}

\noindent {\bf Remarks.}

\noindent {\bf 1.} Here the uniform $C^2$-norm $| .|$ is taken with respect to a fixed finite atlas $\cA$. %Hence the definition of
%a $c$-bounded family depends on the choice of the atlas $\cA$. 

\noindent {\bf 2.} We will assume that the finite atlas $\cA$ is
``nice", that is, it has``extendable" charts, i.e.
charts that can be extended to the (compact) closure of their domains.\vspace{.1in}

\begin{center} {\bf \large  Section 2. Hyperbolic Extensions}\end{center}

 Recall that hyperbolic $n$-space $\HH^n$ is isometric to $\HH^{k}\times \HH^{n-k}$ with
warp metric $(cosh^2\, r)\,\sigma\0{\HH^{k}}+\sigma\0{\HH^{n-k}}$, where $\sigma\0{\HH^{l}}$ denotes the hyperbolic metric of
$\HH^{l}$, and  $r:\HH^{n-k}\ra[0,\infty)$ is the distance to a fixed point in $\HH^{n-k}$.
For instance, in the case $n=2$, since $\HH^1=\R^1$ we have that $\HH^2$ is isometric to $\R^2=\{(u,v)\}$ with warp metric $cosh^2v\, du^2+ dv^2$. In the following paragraph we give a generalization of this
construction.\vspace{.1in}

Let $(M^n,h)$ be a complete Riemannian manifold with {\it center} $o=o_{_M}\in M$. The warp metric 
\begin{center}$g=(cosh^2 r)\, \sigma\0{\HH^k}+h$
\end{center} on $\HH^k\times M$
is the {\it hyperbolic extension (of dimension $k$)} 
of the metric $h$. Here $r$
is the distance-to-$o$ function on $M$.
We write $\cE_k(M,h)=(\HH^k\times M,g)$, and $g=\cE_k(h)$.
We also say that $\cE_k(M)=\cE_k(M,h)$ is the {\it hyperbolic extension
(of dimension $k$) of $(M,h)$} (or just of $M$).
Hence, for instance, we have $\cE_k(\HH^l)=\HH^{k+l}$.
For $S\sbs M$ and $A\sbs\HH^k$ we write $\cE_A(S)=A\times S\sbs\cE_k(M)$.
Also write $\HH^k=\HH^k\times \{o\0{M}\}\sbs\cE_k(M)$ and
we have that any $p\in\HH^k$ is a center of $\cE_k(M)$
(see \cite{O5} or 2.3 below).\vspace{.1in}

Note that $\HH^k$ and every $\{y\}\times M$ are convex in $\cE_k(M)$ (see \cite{BisOn}, p.23).
Let $\eta$ be a complete geodesic line in $M$ passing though $o$
and let $\eta^+$ be one of its two geodesic rays (beginning at $o$) . Then $\eta$ is 
a totally geodesic subspace of $M$ and $\eta^+$ is convex (see \cite{O5}). Also, let $\gamma$ be a complete geodesic line in $\HH^k$. The following two results are proved in \cite{O5}.
\vspace{.1in}

\noindent {\bf Lemma 2.1.} {\it 
We have that\,\, $\gamma\times \eta^+$ is a convex subspace of $\cE_k(M)$
and $\gamma\times \eta$ is totally geodesic in $\cE_k(M)$.}
\vspace{.1in}

\noindent {\bf Corollary 2.2.} {\it We have that \,$\HH^k\times\eta^+$  and $\gamma\times M$ are convex in $\cE_k(M)$. Also  \,$\HH^k\times\eta$
is totally geodesic in $\cE_k(M)$.}\vspace{.1in}

\noindent {\bf Remarks 2.3.}  

\noindent {\bf 1.}
Note that $\HH^k\times\eta$ (with metric
induced by $\cE_k(M)$)
is isometric to $\HH^k\times \R$ with warp metric $cosh^2 v\, \sigma\0{\HH^k}+dv^2$, which is just hyperbolic $(k+1)$-space $\HH^{k+1}$. Also $\gamma\times\eta$
is isometric to $\R\times \R$ with warp metric $cosh^2 v\, du^2+dv^2$, which is just hyperbolic 2-space $\HH^2$.
In particular every point in $\HH^k=\HH^k\times\{ o\} \sbs\cE_k(M)$ is a center point.

\noindent {\bf 2.} It follows from Lemma 2.1 and Remark 2.3(1) that
the ray structure of $\cE_k(h)$ with respect to any center $o\0{\HH^k}\in\HH^k
\sbs\cE_k(M)$ only depends on the ray structure of $M$ and the center $o\0{\HH^k}$.

\noindent {\bf 3.} Denote by $\B_r(M)$ the ball of radius $r$ of $M$.
Note that if $h$ and $h'$ on $M$ have the same ray structures
then the balls $\B_r(M)$ coincide.

\noindent {\bf 4.} Recall that $\HH^k$ is convex in $\cE_k(M)$.
Moreover, for $l\leq k$, we also have $\HH^l\sbs\HH^k\sbs\cE_k(M)$ is
convex. If $h$ and $h'$ on $M$ have the
same ray structures then the $r$-neighborhoods (with respect to
$h$ and $h'$) of the convex subset
$\HH^l$ coincide.\vspace{.1in}

As before (see Section 1) we use $h$ to identify $M-\{ o\}$ with $\bS^{n-1}\times \R^+$. Sometimes we will denote a point
$v=(u,r)\in \bS^{n-1}\times\R^+=M-\{ o\}$ by $v=ru$.
Fix a center $o\in \HH^k\in \cE_k(M)$.  
Since $\HH^k$ is convex in $\cE_k(M)$ we can write $\HH^k-\{ o\}=\bS^{k-1}\times \R^+\sbs\bS^{k+n-1}\times \R^+$
and $\bS^{k-1}\sbs \bS^{k+n-1}$.
Then, for $y\in\HH^k-\{ o\}$ we can also write $y=t\,w$, $(w,t)\in \bS^{k-1}\times\R^+$. 
Similarly, using the exponential map we can identify $\cE_k(M)-\{ o\}$
with $\bS^{k+n-1}\times \R^+$, and for $p\in\cE_k(M)-\{ o\}$
we can write $p=s\,x$, $(x,s)\in\bS^{k+n-1}\times\R^+$.
We denote the metric on $\cE_k(M)$ by $f$ and we can write $f=f_s+ds^2$.
\vspace{.1in}

A point $p\in\cE_k(M)\, -\, \HH^k$ has two sets of coordinates: the {\it polar coordinates}
$(x,s)=(x(p),s(p))\in \bS^{k+n-1}\times \R^+$ and the {\it hyperbolic extension coordinates} $(y,v)=(y(p), v(p))\in \HH^k\times M$. Write $M_o=\{o\}\times M$.
Therefore we have the following functions:\vspace{.1in}

\hspace{.5in}{\small 
\begin{tabular}{|l|l|l|}\hline
the distance to {\it o} function:  & $s:\cE_k(M)\ra [0,\infty)$, & $s(p)=d\0{\cE_k(M)}(p,o)$\\\hline
the direction of {\it p} function:  & $x:\cE_k(M)-\{o\}\ra \bS^{n+k-1},$ & $p=s(p)x(p)$\\\hline
the distance to {\it $\HH^k$} function: & $r:\cE_k(M)\ra [0,\infty),$ &$ r(p)=d\0{\cE_k(M)}(p,\HH^k)$\\\hline
the projection on $\HH^k$ function:  & $y:\cE_k(M)\ra \HH^k, $&\\\hline
the projection on $M$ function:  & $v:\cE_k(M)\ra M, $& \\\hline
the projection on $\bS^{n-1}$ function: & $u:\cE_k(M)-\HH^k\ra \bS^{n-1},$ & $v(p)=r(p)u(p)$\\\hline
the length of $y$ function: & $ t:\cE_k(M)\ra [0,\infty), $&$ t(w)=d_{\HH^k}(y,o)$\\\hline
the direction of $y$ function:  & $w:\cE_k(M)-M_o\ra \bS^{k-1}, $&$ y(p)=t(p) w(p)$\\\hline
\end{tabular}}\vspace{.1in}

Note that $r=d_M(v, o)$. Note also that, by 2.1, the functions $w$ and $u$ are constant on geodesics emanating from $o\in\cE_k(M)$, that is
$w(sx)=w(x)$ and $u(sx)=u(x)$.\vspace{.1in}

Let $\p_r$ and $\p_s$ be the gradient vector fields of $r$ and $s$, respectively. Since the $M$-fibers $M_y=\{ y\}\times M$ are convex
the vectors $\p_r$ are the velocity vectors of the speed one geodesics of the form $a\mapsto (y, a\, u)$, $u\in\bS^{n-1}\sbs M$. These geodesics 
emanate from (and  orthogonally to) $\HH^k\sbs \cE_k(M)$.
Also the vectors  $\p_s$ are the velocity vectors of the speed one geodesics 
emanating from $o\in\cE_k(M)$. For $p\in\cE_k(M)$, denote by $\bigtriangleup =\bigtriangleup (p)$ the right triangle with vertices $o$, $y=y(p)$, $p$
and sides the geodesic segments $[o,p]\in\cE_k(M)$, $[o,y]\in\HH^k$, $[p,y]\in\{ y\}\times M\sbs\cE_k(M)$.
(These geodesic segments are unique and well-defined because:\, (1) $\HH^k$ is
convex in $\cE_k(M)$,\, (2) $(y,o)=o_{_{\{ y\}\times M}}$ and $o$ are centers in $\{ y\}\times M$ and $\HH^k\sbs\cE_k(M)$, respectively.)\vspace{.1in}

Let $\alpha:\cE_k(M)-\HH^k\ra \R$ be the angle between 
$\p_s$ and $\p_r$ (in that order), thus  $cos\, \alpha=f(\p_r,\p_s)$, $\alpha\in [0,\pi]$. 
Then $\alpha=\alpha(p)$ is the interior angle, at $p=(y,v)$, of the right triangle $\bigtriangleup =\bigtriangleup (p)$.
We call  $\beta(p)$ the interior angle of this triangle at $o$, that is $\beta(p)=\beta(x)$ is the spherical distance 
between $x\in \bS^{k+n-1}$ and the totally geodesic sub-sphere $\bS^{k-1}$. Alternatively, $\beta$ is the angle between the geodesic segment
$[o,p]\sbs\cE_k(M)$ and the convex submanifold $\HH^k$.
Therefore $\beta$ is constant on geodesics emanating from $o\in\cE_k(M)$, that is
$\beta(sx)=\beta(x)$. The following corollary follows from 2.1 (see 2.1 in \cite{O5}).\vspace{.1in}

\noindent {\bf Corollary 2.4.} {\it Let $\eta^+$  (or $\eta$) be a geodesic ray (line) in $M$ through $o$ containing
$v=v(p)$ and $\gamma$ a geodesic line in $\HH^k$ through $o$ containing $y=y(p)$. Then $\bigtriangleup (p)\sbs \gamma\times \eta^+\sbs \gamma\times \eta$.}\vspace{.1in}

Note that the right geodesic triangle $\bigtriangleup (p)$ has sides of length $r=r(p)$, $t=t(p)$ and $s=s(p)$. By Lemma 2.1 and Remark 2.3
we can consider $\bigtriangleup$ as contained in hyperbolic 2-space.
Hence using hyperbolic trigonometric identities
we can find relations between $r$, $t$, $s$, $\alpha$ and $\beta$. For instance, using the hyperbolic law of cosines we get:
$cosh\, (s)\, =\, cosh\, (r)\,\, cosh\, (t)$\vspace{.1in}

\noindent Note that this implies $t\leq s$.
Here is an application of this equation. \vspace{.1in}

\noindent {\bf Proposition 2.5 (Iterated hyperbolic extensions)} {\it We have that $$\cE_l\big( \cE_k(M)   \big)=\cE_{l+k}(M)$$
where we are identifying $\HH^{l+k}$ with $\HH^l\times\HH^k$ with warp metric $(cosh^2t)\,\sigma\0{\HH^l}+\sigma\0{\HH^k}$.}\vspace{.1in}

This proposition is proved in \cite{O5}.\vspace{.1in}

\noindent {\bf Remarks.} 

\noindent {\bf 1.} Note that the identification of
$\HH^{l+k}$ with $\HH^l\times\HH^k$ (with warp metric) depends on the order of $l$ and $k$, that is, on the order in which
the hyperbolic extensions are taken.

\noindent {\bf 2.} As before, here the function $t:\HH^k\ra [0,\infty) $ is the distance in $\HH^k$ to the  point $o\in\HH^k$. \vspace{.1in}

We next explore the relationship between
hyperbolic extensions and metrics $\epsilon$-close to hyperbolic. Since $\cE_k(\HH^l)=\HH^{k+l}$ one would
expect that if $M$ is ``close" to $\HH^l$, then
$\cE_k(M)$ would be close to $\HH^{k+l}$.
This motivates the following question.

\vspace{.1in}

\noindent {\bf Question.} {\it What can we say about the hyperbolic extension of a $(B_a,\epsilon)$-close to
hyperbolic metric?}\vspace{.1in}

(Recall that metrics $(B_a,\epsilon)$-close to
hyperbolic are metrics that are already hyperbolic on the ball
$B_{a}=B_a(0)$ of radius $a$, and are radially $\epsilon$-close to hyperbolic outside $B_{a'}$, see Section 1.)\vspace{.1in}

The next result answers this question; it is Theorem B in 
\cite{O5}.\vspace{.1in}

\noindent {\bf Theorem 2.6.}  {\it  Let $M^n$ have center $o$. Assume $M$ is $(B_a,\epsilon)$-close to hyperbolic, with charts of excess $\xi>0$. Then $\cE_k(M)$ is $(B_a,C\epsilon)$-close to hyperbolic, with charts of excess $\xi'$,
provided $a$ is sufficiently large. Explicitly we want
\newline\hspace*{2.5in}$a\,\geq \,R\,=\,R(\epsilon,k,\xi)$
\vspace{.05in}

\noindent Here $C=C(n,k,\xi)$, and  \,$\xi'=\xi-e^{-a/2}>0$.}\vspace{.1in}

This theorem is proved in \cite{O5}. Explicit formulas for
$C$ and $R$ are given in \cite{O5}
(the constant $C$ here is called $C_2$ in \cite{O5}). Note that the excess of the charts decreases.
This is one of the main reasons to introduce the excess. In Section
3 (see also \cite{O4}) we describe another
geometric process, warp forcing, which also reduces the excess of the charts.\vspace{.2in}

\begin{center} {\bf \large  Section 3. Deformations of Metrics}\end{center}

The goal of this section is to describe the ``hyperbolic forcing" method. It has as input a metric on $\R^n$ of
the form $g=g_r+dr^2$ (or, more generally a metric on
a manifold with center) and as output a metric still of the form
$g'_r+dr^2$, but which is hyperbolic on a ball centered
at the origin.\vspace{.1in}

Hyperbolic forcing is defined as the composition of
two other metric deformations:
the two-variable deformation and warp forcing. We present these first.\vspace{.2in}

\noindent {\bf \large 3.1.  The Two Variable Warping Deformation.} 
\vspace{.1in}

Let $g'$ be a metric on the $(n-1)$-sphere $\bS^{n-1}$ and consider the warp metric
$g=sinh^2t\, g' +dt^2$ on $\bS^{n-1}\times \R^+$. Recall that
$\rho:\R\ra[0,1]$ is a fixed smooth function with $\rho (t)=0$ for $t\leq 0$ and $\rho(t)=1$ for $t\geq 1$.  
Given positive numbers $a$ and $d$ define $\rho\0{a,d}(t)=\rho(2\,\frac{t-a}{d})$. 
Also fix an atlas $\cA\0{\bS^n}$ on $\bS^{n-1}$ as before (see remarks at the end 
of Section 1).
All norms and boundedness constants will be taken with respect to
this atlas. Recall that $\sigma\0{\bS^{n-1}}$ is the round metric on $\bS^{n-1}$. Write \vspace{.1in}

\hspace{2in}$g\0{t}=\big(\,1-\rho\0{a,d}(t)\, \big) \sigma\0{\bS^{n-1}} +\rho\0{a,d}(t) \,g'$\vspace{.1in}

\noindent  and define the metric\vspace{.1in}

\hspace{2.25in}$
\cT_{_{a,d}}\, g\, =\, sinh^2\, t\,\,g\0{t}+dt^2
$\vspace{.1in}

We call the correspondence $g\mapsto \cT\0{a,d}g$ the {\it two variable warping deformation}. By construction we have that 
$\cT\0{a,d}g$ satisfies the following property:\vspace{.1in}

$\hspace{1.4in}\cT\0{a,d}g\,=\,\left\{ \begin{array}{lllll}
sinh^2\, (r)\sigma\0{\bS^{n-1}}\, +\, dr^2&& {\mbox{on}}& & B_{a}\\
g&&{\mbox{outside}}& & B_{a+\frac{d}{2}}
\end{array}\right.$\vspace{.1in}

\noindent 
Hence, the two variable warping deformation
changes a warp metric $h$ inside the ball $B_{a+\frac{d}{2}}$ 
making it (radially) hyperbolic on the smaller ball $B_a$.
The warp metric $h$ does not change outside 
$B_{a+\frac{d}{2}}$.\vspace{.1in}

\noindent {\bf Remarks. 3.1.2.} 

\noindent{\bf 1.} Note that if we choose $g$ to be
the warped-by-sinh hyperbolic metric, that is,
$g=sinh^2t\, \sigma\0{\bS^{n-1}} +dt^2$, then
$\cT\0{a,d}g=g$.

\noindent{\bf 2.} To be able to define $\cT\0{a,d}g$ 
the metric $g$ does not need to be a warp metric
everywhere. It only needs to be a warp metric in the ball $B_{a+\frac{d}{2}}$.
%Also to define $\cW\0{r\0{0}} g$, the metric $g$ needs to be defined only
%outside the interior of $B_{r_0}$. 
%\vspace{.1in}

%The following  is Corollary 4.5 in \cite{O3}.}
%Recall that the concept of a $c$-bounded metric was 
%introduced in Section 1.\vspace{.1in}

%\noindent {\bf Theorem 3.1.5.} {\it Let the metric $h$ on $\bS^{n-1}$ be $c$-bounded. Write $g=sinh^2(t)h
%+dt^2$.  Then the metric $\cT_{_{a,d}}\, g\,$ is  $(B_a,\epsilon)$-close to hyperbolic, with charts of excess $\xi$, provided} 
%\begin{center}$\,C_2\,\Big(e^{-a}+\frac{1}{d}\,\Big)\leq\epsilon$
%\end{center}
%\noindent {\it Here $C_2$ is a constant depending on
%$c$, $n$ and $\xi$.} \vspace{.1in}

%\noindent For an explicit formula of $C_2$ see
%\cite{O3}. %The  Corollary in the Introduction \cite{O3} is obtained from 4.5 of \cite{O3} (Theorem 2.4 here)
%by taking $\xi=0$. 
\vspace{.2in}

\noindent {\bf \large 3.2.  Warp Forcing.} 
\vspace{.1in}

Let $(M^n, g)$ be a complete Riemannian manifold with center $o\in M$. Recall that we
can write the metric on $M-\{ o\}=\bS^{n-1}\times\R^+$ as \,\,$g=g_r+dr^2$.  We denote by $\bS_r=\bS_r(M)=\bS^{n-1}\times\{r\}$ the sphere of radius $r$. For a fixed $r\0{0}>0$ we can think of the metric $g\0{r\0{0}}$
as being obtained from $g=g_r+dr^2$ by ``cutting" $g$ along
the sphere of radius $r\0{0}$, so we call
$g\0{r\0{0}}$ the {\it warped spherical cut of
$g$ at $r\0{0}$.}  In the same vein, we call the metric
\begin{center}$\hat{g}\0{r\0{0}}\, =\, \bigg(\frac{1}{sinh^2(r_0)}\bigg)\, g \0{r\0{0}} $\end{center}
\noindent the {\it (unwarped by sinh) spherical cut of $g$ at $r\0{0}$ }.
Note that in the particular case where $g=g_r+dr^2$ is 
already a warped-by-$sinh$
metric (that is, $g_r=sinh^2(r)g'$ for some fixed
$g'$ independent of $r$) we have that the warped spherical cut
of $g=sinh^2(r)g'+dt^2$ at $r\0{0}$ is $sinh^2(r\0{0})g'$,
and the spherical cut at $r\0{0}$ is
$\hat{g}\0{r\0{0}}=g'$. Hence the terms ``warped" and
``unwarped" (usually we will omit the term ``unwarped").
\vspace{.1in}

Fix $r\0{0}>0$.
We define the warped-by-$sinh$ metric $\bg\0{r\0{0}}$ by:\vspace{.1in}

\hspace{1.4in}$\bg \0{r\0{0}}\,=\,sinh^2(t)\hg\0{r\0{0}}\,+\, dr^2     \,=\, sinh^2\, (t)\big(  \frac{1}{sinh^2(r\0{0})}\big) g\0{r\0{0}}\, +\, dr^2
$\vspace{.1in}

\noindent We now force the metric $g$ to be equal to $\bg\0{r\0{0}}$ on  $B_{r\0{0}}=\B_{r\0{0}}(M)$ and stay equal to $g$ outside $B_{r\0{0} +\frac{1}{2}}$.
For this we define the {\it warped forced } (on $B_{r\0{0}}$) metric as:\vspace{.1in}

\hspace{2in}
$\cW\0{r\0{0}}\, g\, =\, (1-\rho\0{r\0{0}}) \, \bg\0{r\0{0}}\, +\, \rho\0{r\0{0}} \, g
$\vspace{.1in}

\noindent where $\rho\0{r\0{0}}(t)=\rho(2t-2r\0{0})$, and $\rho :\R\ra [0,1]$ is as before (see Section 1). Hence we have\vspace{.1in}

\hspace{1.8in}$\cW\0{r\0{0}} g\,=\,\left\{ \begin{array}{lllll}
\bg\0{r\0{0}}&& {\mbox{on}}& & B_{r\0{0}}\\
g&&{\mbox{outside}}& & B_{r\0{0}+\frac{1}{2}}
\end{array}\right.$\vspace{.1in}

Hence warp forcing changes
the metric only on $B_{r\0{0}+\frac{1}{2}}$, making it a warp metric inside $B_{r\0{0}}$. The metric
$g$ does not change outside
$B_{r\0{0}+\frac{1}{2}}$.
We call the process $g\mapsto\cW g$ {\it warp forcing}.
%The next result is the Main Theorem of \cite{O4}. It states that if $g$ is {\it $\epsilon$-close to a hyperbolic metric} then
%the warp forced metric  $\cW\0{r\0{0}} g$ is also close to hyperbolic. 
\vspace{.1in}

%\noindent {\bf Theorem 3.2.1.} {\it  Let $(M,g)$ have center $o$, and $S\sbs M$. If $g$ is radially $\epsilon$-close to hyperbolic on $S$ 
%with charts of excess $\xi$,
%then $\cW\0{r\0{0}} g$ is radially $\eta$-hyperbolic on $S-B\0{r\0{0}-(1+\xi)}$ with charts of excess $\xi-1$, provided $\eta\geq e^{16+6\xi}\big(e^{-2r\0{0}}+\epsilon\big)$.}
%\vspace{.1in}

%The prove of this result is given in \cite{O4}.\vspace{.1in}

\noindent{\bf Remarks 3.2.1.} 

\noindent {\bf 1.} Notice that to define $\cW\0{r\0{0}} g$ we only need
$g\0{r}$ to be defined for $r\geq r\0{0}$. 
%But to apply 3.2.4 we need 
%$g\0{r}$ to be defined for $r\geq r\0{0}-1-\xi$.

%\noindent {3.} Notice that $\cW\0{r\0{0}} g$ is a warp metric no just on
%$B\0{r\0{0}}$ but on $B\0{r\0{0}+\delta}$, for small $\delta>0$.

\noindent{\bf 2.} Note that if we choose $g$ to be
the warped-by-sinh hyperbolic metric, that is,
$g=sinh^2t\, \sigma\0{\bS^{n-1}} +dt^2$, then
$\cW\0{r\0{0}}g=g$

\vspace{.2in}

\noindent {\bf \large 3.3.  Hyperbolic Forcing.} 
\vspace{.1in}

Let $(M^{n},g)$ have center $o$. As before we write $g=g\0{r}+dr^2$.
Let  $r\0{0}>d>0$. We define
the  metric $\cH\0{r\0{0},d}\,\big( g)$ in the following way.
First warp-force the metric $g$, i.e take $\cW\0{r\0{0}} g$. Recall $\cW\0{r\0{0}}
g$ is a warp metric on $B\0{r\0{0}}$ (see also remarks 3.1.2 and 3.2.1). Hence we can use two variable warping 
deformation (see 3.1)  and  define\vspace{.1in}

\noindent {\bf (3.3.1)}\hspace{1.9in}
$\cH\0{r\0{0},d}\,g\,=\, \cT_{_{(r\0{0}-d),d}}\big( \, \cW\0{r\0{0}} g\,\big)\,\, $
%=\, \cT_{_{(r\0{0}-d),d}}\big( \,\bg\0{r\0{0}}\big) $

\vspace{.1in}

\noindent 
The process $g \mapsto \cH\0{r\0{0},d}$ is called {\it hyperbolic forcing}.
Write $h=\cH\0{r\0{0},d}\,g\,$. 
Note that $h$ also has the form  $h=h_r+dr^2$.
In the next results we explicitly describe $ h_r$ and give
some properties of the metric $h=\cH\0{r\0{0},d}\,g$. These results are proved in
\cite{O6}.\vspace{.1in}

\noindent {\bf Proposition 3.3.2.} {\it We have}
\vspace{.1in}

\noindent \,\,\,\,\,\,\,\,\,\,\,\,\,{\small $
h_r\,\, =\,\,\left\{
\begin{array}{lll}
g_r&&   r\0{0} +\frac{1}{2}\leq r\\
\big(1-\rho\0{r\0{0}} (r)\big) \,sinh^2 (r)\, \hg\0{r\0{0}} \, +\,\rho\0{r\0{0}} (r) \,g\0{r}& 
&r\0{0} \leq r\leq r\0{0} +\frac{1}{2}\\
sinh\,^2 (r) \bigg( \,\big( 1-  \rho_{_{(r\0{0} -d), d}}(r)\, \big)\, \sigma\0{\bS^{n-1}}\,+\, \rho_{_{(r\0{0} -d), d}}(r)  \hg\0{r\0{0}}\,  \bigg)&&  r\0{0} -d \leq r\leq r\0{0}  \\
sinh\,^2(r)\, \sigma\0{\bS^{n-1}}&&r\leq r\0{0}-d
\end{array}
\right.
$}\vspace{.1in}

\noindent {\it where the gluing functions $\rho\0{r\0{0}}$ and $ \rho_{_{(r\0{0} -d), d}} $ are defined in 3.2 and 3.1, respectively.}\vspace{.1in}

\noindent {\bf Proposition 3.3.3.}  {\it The metric $h=\cH\0{r\0{0},d}\,g\,$ has the following properties.}

{\it 
\noindent (i) The metric $h$ is canonically hyperbolic on $B_{_{r\0{0}-d}}$, i.e equal to $sinh^2(r)\sigma\0{\bS^n}+dr^2$ on
$B\0{r\0{0}-d}$.

\noindent (ii) We have that $g=h$ outside $B_{_{r\0{0}+\frac{1}{2}}}.$

\noindent (iii) The metric $h$ coincides with $\cW\0{r\0{0}} \big(g\0{r\0{0}}\big)$ outside $B\0{r\0{0}-\frac{d}{2}}$.

\noindent (iv) The metric $h$ coincides with 
$\cT_{_{(r\0{0}-d),d}} \bar{g}\0{r\0{0}}$
 on $B\0{r\0{0}}$.

\noindent (v) All the $g$-geodesic rays $r\mapsto ru$, $u\in \bS^n$, emanating from the center are geodesics of
$(M , h)$. Hence, the space  $(M,h)$ has center $o$. Moreover
the function $r$ (distance to the center $o$) is the same on the spaces $(M,g)$ and $(M,h)$. In other words, the spaces $(M,g)$ and $(M,h)$ have the same ray structures.}\vspace{.13in}

Next we discuss the following question:\vspace{.1in}

\noindent  {\it Is the hyperbolically forced metric 
$h=\cH\0{r\0{0},d}\, g$ close to hyperbolic, when
$g$ is close to hyperbolic?}\vspace{.1in}

Notice that from 3.1.2 and 3.2.1 it follows that if we choose $g$ to be
the warped-by-sinh hyperbolic metric, that is,
$g=sinh^2t\, \sigma\0{\bS^{n-1}} +dt^2$, then
$\cH\0{r\0{0},d}\, g=g$. Therefore one would expect
that the answer to the previous question is ``yes".
So, it is better ask a more quantitative question:\vspace{.1in}

\noindent {\it To what extend is the hyperbolically forced metric 
$h=\cH\0{r\0{0},d}\, g$ close to hyperbolic, when
$g$ is close to hyperbolic?}\vspace{.1in}

The next theorem deals with this question.
This theorem is proved in \cite{O6}.\vspace{.1in}

\noindent {\bf Theorem 3.3.4.} {\it Let $M^{n}$ have center $o$ and metric $g=g\0{r}+dr^2$. Assume the spherical cut $\hg\0{r\0{0}}$ is {\s c}-bounded.
If the metric $g$ is radially $\epsilon$-close to hyperbolic outside $B\0{r\0{0}-(1+\xi)}$
with charts of excess $\xi>1$, then the metric $\cH\0{r\0{0},d}\, g$ is $(B_{r\0{0}-d},\eta$)-close to hyperbolic
with charts of excess $\xi-1$, provided} \begin{center}$\eta\geq
C_1\,\Big(\,\frac{1}{d}\,+\,e^{-(r\0{0}-d)}\Big)\,+\, C_2
\,\epsilon$\end{center}
\noindent {\it  Here 
$C_1$ is a constant depending only on $n$, $\xi$,  {\s c}, and $C_2$ depends only on $\xi$.}
\vspace{.1in}

\noindent {\bf  Remarks 3.3.5.}

\noindent {\bf 1.} An important point here is that by taking $r\0{0}$ and $d$ large the metric $\cH\0{r\0{0},d}\, g$ can be made \,$2C_2\epsilon$-close to hyperbolic.
How large we have to take $d$ and $r\0{0}$ depends on {\s c}, which is a $C^2$ bound
for the the $\hat{g}\0{r\0{0}}$, the spherical
cut of $g$ at $r\0{0}$ (see 3.2).

\noindent {\bf 2.}  Note that the excess of the charts decreases by 1.
This is because of warp forcing (see \cite{O4}).
\vspace{.2in}

\begin{center} {\bf \large  Section 4. Deformations of Families of Metrics}\end{center}

In this section we give a one-parameter version of the 
concepts and results presented in Section 3. Let $(M^{n}, g)$ be a complete Riemannian manifold with center $o\in M$. Recall that we
can write the metric on $M-\{ o\}=\R^{n}-\{0\}=\bS^{n-1}\times\R^+$ as \,\,$g=g\0{r}+dr^2$,
where $r$ is the distance to $o$.\vspace{.1in}

Fix $\xi>0$, and let $\lambda_0>1+\xi$.
We say that the collection $\{ g_\ssl\}_{\lambda\geq   \lambda_0}$ is a $\odot${\it-family of metrics on $M$} if each $g_\ssl$ is 
a  metric  of the form $g_\ssl=\big(g_\ssl\big)_r+dr^2$ defined (at least) for $r>\lambda-(1+\xi)$.\vspace{.1in}

\noindent {\bf Remark.} We will always assume that
the family
of metrics $\{ g_\ssl\}$ is smooth, that is, the map
$(x,\lambda)\mapsto g_\ssl(x)$ is smooth, $x\in M$,
$\lambda\geq\lambda_0$.\vspace{.1in}

We say that the $\{g_\ssl\}$ has {\it (spherical) cut limit at $b$} if there is a 
$C^2$ metric $\hg_{_{\infty+b}}$ on $\bS^{n-1}$ such that \vspace{.1in}

\noindent {\bf (4.1)}\hspace{1.4in}$
\big|\,\,{\widehat{\big(g_\ssl\big)}}_{_{\lambda+b}} \,\,-\,\, \hg_{_{\infty+b}}\,\,\big|\,\,\longrightarrow\,\, 0
$  \,\,\,\,\,\,\,as\,\,\,\,\,\, $\lambda\ra \infty$\vspace{.1in}

\noindent {\bf Remarks 4.2.}

\noindent {\bf 1.} Recall that the metric ${\widehat{\big(g_\ssl\big)}}_{_{\lambda+b}}$ is the spherical cut of $g_\ssl$ at
$\lambda+b$. See Section 3.2.

\noindent {\bf 2.} The arrow above means convergence in the $C^2$-norm on the space of $C^2$ metrics on  $\bS^{n-1}$.
See remarks at the end of Section 1.

\noindent {\bf 3.} The definition above implies that $\big(g_\ssl\big)_{_{r+b}}$ is defined for large $\lambda$,
even if $b<-(1+\xi)$.

\noindent {\bf 4.} Note that the concept of cut limit at $b$ depends strongly on the indexation of the family.

\noindent {\bf 5.}  If
an  family $\{ g_\ssl\}$ has cut limits at $b$,
then the family
$\{\,\widehat{(g_\ssl)}\0{\lambda+b}\,\}_\ssl$
is {\s c}-bounded, for some {\s c}; see \cite{O7}.
% It is an ``interval version" of 4.5.
% {\it If
%an  $\odot$-family $\{ g_\ssl\}$ has cut limits on
%the compact interval $I$,
%then the family
%$\{\,\widehat{(g_\ssl)}\0{\lambda+b}\,\}\0{\lambda,b\in I}$
%is $c$-bounded, for some $c$.}
\vspace{.1in}

Consider the  $\odot$-family $\{g_\ssl\}$ and let $d>0$. Apply hyperbolic forcing to
get \begin{center} $h_\ssl=\cH\0{\lambda, d} g_\ssl $\end{center} We say that the family $\{h_\ssl\}$ is the {\it hyperbolically forced family} 
corresponding to the $\odot$-family $\{g_\ssl\}$.
Note that we can write $h_\ssl=( h_\ssl )_r+dr^2$.
Using Proposition 3.3.2 we can explicitly describe $( h_\ssl )_r$:
(see \cite{O6} for more details)\vspace{.1in}

\noindent {\bf (4.3.)}\,\,\,\,\,\,\,\,{\small $
\big( h_\ssl \big)_r\,\, =\,\,\left\{
\begin{array}{lll}
\big( g_\ssl \big)_r&&   \lambda +\frac{1}{2}\leq r\\
\big(1-\rho_\ssl (r)\big) \,sinh^2 (r)\, {\widehat{\big(g_\ssl \big)}}_\ssl \, +\, \rho_\ssl(r)\,\big( g_\ssl \big)_r&
&\lambda \leq r\leq \lambda +\frac{1}{2}\\
sinh\,^2 (r) \bigg( \big( 1-  \rho_{_{(\lambda -d), d}}(r) \big)\, \sigma\0{\bS^{n-1}} \,+\,
\rho_{_{(\lambda -d), d}}(r)  \,{\widehat{\big(g_\ssl \big)}}_\ssl\,
 \bigg)&&  \lambda -d \leq r\leq \lambda  \\
sinh\,^2(r)\, \sigma\0{\bS^{n-1}}&&r\leq \lambda-d
\end{array}
\right.
$}\vspace{.1in}

The next proposition is a one-parameter version of
3.3.3. It is proved in \cite{O6}.\vspace{.1in}

\noindent {\bf Proposition 4.4.}  {\it The metrics $h_\ssl$ have the following properties.

\noindent (i)  The metrics $h_\ssl$ are canonically hyperbolic on $B_{_{\lambda-d}}$, i.e equal to $sinh^2(r)\sigma\0{\bS^{n-1}}+dr^2$ on
$B\0{\lambda-d}$, provided $\lambda >d$.

\noindent (ii) We have that $g_\ssl=h_\ssl$ outside $B_{_{\lambda+\frac{1}{2}}}.$

\noindent (iii)  The metric $h$ coincides with $\cW\0{\lambda} \big(g\0{\lambda}\big)$ outside $B\0{\lambda-\frac{d}{2}}$.

\noindent (iv)  The metric $h$ coincides with 
$\cT_{_{(\lambda-d),d}}\Big( {\overline{(g_\ssl )}}_\ssl\,\Big)$
 on $B\0{\lambda}$.

\noindent (v) 
If the $\odot$-family $\big\{ g_\ssl\big\}$ has  cut limits for $b=0$ then
$\big\{ h_\ssl\big\}$ has cut limits on $(-\infty,0]$. In fact we have
\begin{center}$
\hh\0{_{\infty+b}} \,\, =\,\,\left\{
\begin{array}{lll}
\hg_{_{\infty}}&&  b=0\\
\big( 1-  \rho(2+\frac{2b}{d})  \big)\, \sigma\0{\bS^n} 
\,+\, \rho(2+\frac{2b}{d})  \,\hg_{_{\infty}} &&   -d \leq b\leq 0  \\
 \sigma\0{\bS^n}&&b\leq -d
\end{array}
\right.
$\end{center}
\noindent where $\rho$ is as in Section 2.

\noindent (vi)   If we additionally assume that $\{g_\ssl\}$ has  cut limits on $[0,\frac{1}{2}]$, then
$\big\{ h_\ssl\big\}$ has also  cut limits on $[0,\frac{1}{2}]$. In fact, for $b\in [0,\frac{1}{2}]$ we have
\begin{center}$
\hh\0{_{\infty+b}} \,\, =\,\,\big( 1-\rho(b) \big)\,\hg_{_{\infty}}\,\,+\,\, \rho(b) \,\hg_{_{\infty+b}}
$\end{center}
\noindent where $\rho$ is as in Section 3. Of course if $\{g_\ssl\}$ has a cut limit at $b> \frac{1}{2}$ then $\{h_\ssl\}$ 
has the same cut limit at $b$ (see item (ii)).

\noindent (vii)  All the rays $r\mapsto ru$, $u\in \bS^n$, emanating from the origin are geodesics of
$(M , h_\ssl)$. Hence, all spaces  $(M,h_\ssl)$ have center $o\in M$
and have the same geodesic rays emanating from the common center $o$. Moreover
the function $r$ (distance to  $o\in M$) is the same on all spaces $(M,h_\ssl)$.}\vspace{.1in}

We now state one of our most important results. 
It is used in an essential way in smoothing
Charney-Davis strict hyperbolizations. It is proved in 
\cite{O6} using 3.3.4. Before, we need a definition.
We say that
an $\odot$-family $\{g_\ssl\}$  is {\it radially $\epsilon$-close to hyperbolic,
with charts of excess $\xi$}, if each
$g\0{\lambda}$ is radially $\epsilon$-close to hyperbolic outside $B_{\lambda-(1+\xi)}$,
with charts of excess $\xi$.\vspace{.1in}

\noindent {\bf Theorem 4.5.} {\it Let $M$ have center $o$, 
$\{g\0{\lambda}\}$
an $\odot$-family on $M$, and $\epsilon'>0$. Assume that $\{g\0{\lambda}\}$ has cut limits at $b=0$.
If $\{ g\0{\lambda}\}$ is radially $\epsilon$-close to hyperbolic, with charts of excess $\xi>1$, then  $\cH\0{\lambda,d}\,g\0{\lambda}$ is $(B_{\lambda-d},\epsilon'+C_2\epsilon)$-close to hyperbolic,
with charts of excess $\xi-1$, provided\vspace{.05in}

\noindent \,\,\,\,(i) $\lambda-d>ln(\frac{2C_1}{\epsilon'})$\vspace{.05in}
  
\noindent\,\,\,\,(ii) $d\geq \frac{2C_1}{\epsilon'}$.
\vspace{.05in}

\noindent Here $C_1$ and $C_2$ are as in Theorem 3.3.4.}\vspace{.1in}

\noindent {\bf Remarks 4.6.}

\noindent {\bf 1.} Note that we can take $\epsilon'$ as small as we want hence $\epsilon'+C_2\epsilon$
as close as $C_2\epsilon$ as we desire, provided we choose $d$ and $\lambda$
sufficiently large. How large depending on $\epsilon'$
and $c$.

\noindent {\bf 2.} The constant $C_1(c,n,\xi)$ in Theorem 
3.3.4 depends on $c$, which is a $c$ bound for
the limit metric $\hg_{_{\infty+b}}$ (see 4.1). This $c$
exists (see 4.2 (5)).
\vspace{.2in}

\noindent {\bf \large 4.7 Cuts Limits and Hyperbolic Extensions.} 
\vspace{.1in}

At the beginning of this section we gave the definition of cut limit (see 4.1).
More generally, let $I\sbs\R$ be an interval (compact or noncompact). We say the  $\odot$-family $\{g_\ssl\}$ has {\it cut limits on $I$} if the convergence in (4.1) is uniform
in $b\in I$. Explicitly this means: for every $b\in I$ and $\epsilon>0$,
there are $\lambda_*$ and neighborhood
$U$ of $b$ in $I$ such that
{\scriptsize $\big|\,\,{\widehat{\big(g_\ssl\big)}}_{_{\lambda+b'}} \,\,-\,\, \hg_{_{\infty+b'}}\,\,\big|<\epsilon$}, for $\lambda>\lambda_*$ and $b'\in U$.
In particular $\{g_\ssl\}$ has a cut limit at $b$, for every $b\in I$. \vspace{.1in}

If the  $\odot$-family $\{g_\ssl\}$ has cut limits on $\R$
we will just say that $\{g_\ssl\}$  {\it has  cut limits.}\vspace{.1in}

\noindent {\bf Remark 4.7.1.} Let $a\in \R$.  If 
$\{g_\ssl\}\0{\lambda}$ has cut limits then so does
the reparametrized family $\{g\0{\lambda +a}\}\0{\lambda }$.

\vspace{.05in}
Here is a natural question:\vspace{.05in}

\noindent {\it If
$\{h_\ssl\}\0{\lambda}$ has a cut limits, does 
$\{\cE_k(h_\ssl)\}\0{\lambda}$ have cut limits?} \vspace{.1in}

\noindent {\bf Remark.} 
More generally we can ask whether $\{\cE_k(h\0{\lambda})\}\0{\lambda'}$
has cut limits, where $\lambda=\lambda(\lambda')$. Of course the answer
would depend on the change of variables $\lambda=\lambda(\lambda')$.\vspace{.1in}
 
The next result gives an affirmative answer to this
question provided the family $\{h_\ssl\}$ 
is, in some sense, nice near the origin. Explicitly,
we say that $\{h_\ssl\}$ is {\it hyperbolic
around the origin} if there is a $B\in \R$ such that
%\vspace{.1in}

\hspace{2.3in}
${\widehat{\big(h_\ssl\big)}}_{_{\lambda+b}}=
\sigma\0{\bS^{n-1}}$\vspace{.1in}

\noindent for every $b\leq B$ and every
(sufficiently large) $\lambda$.
Note that this implies that each $h_\ssl$ is
canonically hyperbolic on the ball of radius $\lambda +B$.
Examples of $\odot$-families that are
hyperbolic around the origin are families obtained
using hyperbolic forcing, as above.\vspace{.1in}

As mentioned before the next result answers affirmatively the question above.
Moreover it also says that some reparametrized families $\{\cE_k(h\0{\lambda})\}\0{\lambda'}$,
for certain change of variables $\lambda=\lambda(\lambda')$, have cut limits as well.
Write $\lambda=\lambda(\lambda',\theta)=sinh^{-1}(sinh(\lambda')\,sin\,\theta)$, for fixed $\theta$.
Note that $\lambda=\lambda'$ for $\theta=\pi/2$.
We say that $\{\cE_k(h\0{\lambda})\}\0{\lambda'}$ is the {\it $\theta$-reparametrization
of} $\{\cE_k(h\0{\lambda})\}\0{\lambda}$. 
Note that if we
consider an hyperbolic right triangle
with one angle equal to $\theta$ and 
side (opposite to $\theta$) of length
$\lambda$, then 
$\lambda'$ is the
length of the hypothenuse of the triangle. All $\theta$-reparametrizations, in the limit $\lambda'\ra \infty$, differ
just by translations. The next proposition is proved in \cite{O7}.\vspace{.1in}

\noindent {\bf Proposition 4.7.2.} 
{\it Let $M$ have center $o$. Let  $\{h_\ssl\}\0{\lambda}$ be $\odot$-family of metrics
on $M$.
Assume $\{h_\ssl\}\0{\lambda}$ is hyperbolic around the origin.
If $\{h_\ssl\}$ has cut limits, 
then the $\theta$-reparametrization $\{\cE_k(h\0{\lambda})\}\0{\lambda'}$ has cut limits as well. Here 
$\theta\in(0,\pi/2]$.}\vspace{.2in}

\begin{center} {\bf \large  Section 5. Normal Neighborhoods on All-Right Spherical Complexes}\end{center} 

In this section we define  and and give some properties of neighborhoods of simplices in
all-right spherical complexes. The goal is to
define ``natural normal neighbohoods" of simplices
in all-right spherical complexes, and
give some of its properties.\vspace{.1in}

We use the definition and properties of a spherical complex given in Section 1 of \cite{ChD}. Recall that a spherical
complex is an {\it all-right spherical complex} if all of its edge lengths are equal to 
$\pi/2$. 
Given an all-right spherical complex $P$ we will use the same symbol $P$ for the complex itself (the collection of all simplices), and its realization (the union of all its simplices).
 In this paper we shall assume that all spherical complexes satisfy the
``intersection condition" of simplicial complexes: {\sf every two simplices intersect in at most one common face}.\vspace{.1in}

\noindent {\bf Remark 5.0.1.} Let $P$ be an all-right spherical
complex and $\Delta\in P$. The symbol $\dDelta$ denotes the interior of $\Delta$.
In this paper we will use the three definitions of link
$\sL (\Delta,P)$ of $\Delta$ in $P$. The {\it geometric link} 
 $\sL(\Delta,P)$ is the union of the end points of geodesic segments of small length $\beta>0$
emanating perpendicularly (to $\Delta$) from some point  
$x\in \dDelta$. If we want to specify $\beta$ and $x$ we
say that $\sL _\beta(\Delta,P)$ is the $\beta$-link 
{\it based at $x$}.
The {\it geometric star} $\s{Star}(\sigma,K)$ is the union of 
the corresponding segments. The {\it simplicial link} 
is the subcomplex of $P$ formed by all simplices $\Delta'$ such that (1) 
$\Delta'$ is disjoint from $\Delta$, (2) $\Delta'$ and $\Delta$
span a simplex (this simplex is the join $\Delta\ast\Delta'\in P$,
and $\Delta'$ is the {\it opposite face} of $\Delta$ in $\Delta\ast\Delta'$). Note that if we continue a geodesic $[x,u]$, with $u$ in the geometric $\beta$-link at $x$, we will hit a unique point in $\Delta'$. This radial geodesic projection gives a relationship between 
geometric links and simplicial links. The {\it simplicial star}
is the subcomplex of $P$ formed by all simplices $\Delta'$ that  
contain $\Delta$. 
For $x\in\dDelta^k$ the {\it direction link of $\Delta$ in $P$ at $x$}
is the set of all vectors at $x$ perpendicular to $\Delta^k$. Using geodesics
emanating from $x$ perpendicularly to $\Delta$ we also get a relationship
between geometric links and the direction links.
These different definitions of link all come with natural all-right spherical metrics:
the geometric link with the rescaled induced metric, the simplicial link with the induced metric and the direction link with the
angle metric. The relationships between the different definitions of link mentioned
above all respect the metrics.\vspace{.2in}

\noindent {\bf \large 5.1 Sets of Widths of Normal Neighborhoods on the Sphere $\bS^m$.}\vspace{.1in}

We consider the $m$-sphere $\bS^m\sbs\R^{m+1}=\{x=(x_1,...,x_{m+1})\}$ with the {\it canonical all-right spherical structure} whose
$m$-simplices are $\bS^m\cap \{ (-1)^{s_i}x_i\geq 0\}$, for any
choice $s_i\in\{0,1\}$, $i=1,...,m+1$.  Let $\beta\in (0,\pi/2)$
and $\Delta\in\bS^m$. The closed normal neighborhood of 
$\Delta$ in $\bS^m$ of width $\beta$ is the union of
(images of) geodesics of length  $\beta$ emanating 
perpendicularly from $\Delta$. It will be denoted by 
${\sf{N}}_{_{\beta}}(\Delta,\bS^m)$. 
For the special case $dim\,\Delta=m$ we will take
${\sf{N}}_{_{\beta}}(\Delta^m,\bS^m)=\Delta^m$, for any $\beta$.%We make the convention that
%${\sf{N}}_{_{\beta}}(\Delta,\bS^m)$ is equal to the simplicial
%star of $\Delta$ when $\beta\geq\pi/2$.
\vspace{.1in}

Let ${\sf{B}}=\{\beta_k\}_{k=0,1,2...}$ be an indexed set of real numbers with $\beta_k\in(0,\pi/2)$
and $\beta_{k+1}<\beta_k$. We write $\s{B}(m)=\{\beta_0,...\beta_{m-1}\}$.
The set $\s{B}$ determines {\it the set of spherical} \,$\s{B}$-{\it neighborhoods}\,
${\sf{N}}_{\sf{B}}(\bS^m)={\sf{N}}_{\s{B}(m)}(\bS^m)=\{ {\sf{N}} _{_{\beta_k}}(\Delta^k,\bS^m) \}_{\Delta^k\in\bS^m,\, k<m}$,
for any sphere $\bS^m$ (of any dimension). Note that the normal neighborhoods of all
$k$-simplices $\Delta^k$ have the same width $\beta_k$.
The set $\s{B}$ is called a {\it set of widths of spherical normal neighborhoods}   or simply {\it a set of widths}. The set $\s{B}(m)$ is a {\it finite set of widths
of length $m$}. The definitions above still make sense if we replace $\bS^m$ by  $\bS^m_\mu$, the $m$-sphere of radius $\mu$ (for small $\beta_k$'s).\vspace{.1in}

We are interested in pairs of sets of widths  $\big(\s{B},\s{A}\big)$,
$\s{B}=\{\beta_k\}$ and $\s{A}=\{\alpha_j\}$, having the following 
{\it {\bf D}isjoint {\bf N}eighborhood {\bf P}roperty}:\vspace{.1in}

\noindent {\bf (5.1.2.)\, DNP:}\,\,\,\,\,{\it For every $k$ and $m$,
$k<m$,  the following sets 
are disjoint}\vspace{.05in}

\hspace{1.8in}{\small $\bigg\{ \, \s{N}_{\beta_k}(\Delta^k,\bS^m)\, \,-\, \,\bigcup_{j<k} \s{N}_{\alpha_j}(\Delta^j,\bS^m) \,\bigg\}_{\Delta^k\in\bS^m}$} \vspace{.1in}

\noindent The disjoint neighborhood property obtained by fixing $k$ and $m$
above will be denoted by {\bf DNP($k,m$)}. 
In this case we allow the sets of widths to be finite of length at least $k+1$.
Note that the ordering of the pair $\big(\s{B},\s{A}\big)$
is important. It is straightforward to verify that {\bf DNP($k,m$)}
(and {\bf DNP})
is equivalent to the following property. For fixed (any)
$k$ and $m$ we have: for different $k$-simplices
$\Delta^k_1$ and $\Delta^k_2$ we have\vspace{.1in}

\noindent {\bf (5.1.2)'} \hspace{1in}{\small $\s{N}_{\beta_k}(\Delta_1^k,\bS^m)\,\bigcap\, \s{N}_{\beta_k}(\Delta_2^k,\bS^m)\,\,\,\,
{\mbox{{\LARGE $\sbs$}}}\,\,\, \,\,\bigcup_{j<k} \s{N}_{\alpha_j}(\Delta^j,\bS^m) $}\vspace{.1in}

\noindent That is, the $\s{B}$-neighborhoods of different $k$-simplices intersect only inside the $\s{A}$-neighborhood of the $(k-1)$-skeleton (which is equal to \,$\bigcup_{j<k} \s{N}_{\alpha_j}(\Delta^j,\bS^m)$).\vspace{.1in}

\noindent {\bf Proposition 5.1.3.} {\it The pair of
(finite or infinite) sets of widths $\big(\s{B},\s{A}\big)$
satisfy} \,{\bf DNP($k,m$) } {\it if and only if} \,
$\frac{sin\,\beta\0{k}}{sin\,\alpha\0{k-1}}\leq\frac{\sqrt{2}}{2}$.
\vspace{.1in}

\noindent Note that the inequality condition is independent of $m$.
The proposition follows directly from lemmas 5.1.4 (taking $k=l$ and $\beta=\gamma$) and 5.1.4 given below, and the fact that $\{\alpha\0{k}\}$ is decreasing.\vspace{.1in}

\noindent{\bf Lemma 5.1.4.} {\it  Let
 $\Delta^k$, $\Delta^l\in\bS^m$
and $\Delta^j=\Delta^k\cap\Delta^l$. Let $\alpha, \beta, \gamma\in (0,\pi/2)$
such that \,$\frac{sin\,\beta}{sin\,\alpha}, \,\frac{sin\,\gamma}{sin\,\alpha}\,\leq\,\frac{\sqrt{2}}{2}$. Then }
$$\s{N}_{\beta}(\Delta^k,\bS^m)\,\cap\, \s{N}_{\gamma}(\Delta^l,\bS^m)\,\,\,\,
{\mbox{{\Large $\sbs$}}}\,\,\, \s{N}_{\alpha}(\Delta^j,\bS^m) $$
\noindent {\bf Proof.} 
In this proof $\sL(\Delta,\bS^m)$ shall denote the simplicial link and $\s{Star}(\Delta,\bS^m)$ the simplicial star (see 5.0.1). Note that $\s{N}_\beta(\Delta,\bS^m)\sbs
\s{Star}(\Delta,\bS^m)$,  for every $\Delta\in\bS^m$.
Write $\s{S}=\sL(\Delta,\bS^m)$, $\Delta'_1=\s{S}\cap\Delta^k$ and $\Delta'_2=\s{S}\cap\Delta^l$.
Then $\Delta'_i$ is a simplex in the all-right triangulation of $\s{S}$.
Also $\Delta'_1$ and  $\Delta'_2$ are disjoint.
Hence their distance in $\s{S}$ is at least $\frac{\pi}{2}$.\vspace{.1in}

Suppose there is $q\in\s{N}_{\beta}(\Delta^k,\bS^m)\cap
\s{N}_{\gamma}(\Delta^l,\bS^m)$. Since both of these neighborhoods lie in
$\s{Star}(\Delta^j,\bS^m)$ there is geodesic segment $[p,q]$ in $\s{Star}(\Delta^j,\bS^m)$
with $p\in \Delta^j$ and $[p,q]$ perpendicular to $\Delta^j$ (note that
$p$ may lie in $\p\Delta^j$ and that the geodesic may not be unique
if $q\in\s{S}$). Write $\alpha'=d_{\bS^m}(p,q)$.
We have to prove $\alpha'\leq\alpha$. We assume $\alpha'>\alpha$ and get
a contradiction.
Let $q_1$ be the closest point in $\Delta^k$ to $q$
and $q_2$ be the closest point in $\Delta^l$ to $q$. We have
$a_1=d_{\bS^m}(q_1,q)\leq \beta$ and $a_2=d_{\bS^m}(q_2,q)\leq \gamma$.
We get a right (at $q_i$) spherical triangle with one side equal to $a_i$ and hypotenuse equal to $\alpha'$. Let $\theta_i$ be
the angle at $p$, that is, the angle opposite to the side of length $a_i$. Then by the spherical law of sines we get
{\footnotesize $$
sin\,\theta_1\,=\,\frac{sin\, a_1}{sin\,\alpha'} < \frac{sin\,\beta}{sin\,\alpha}
\leq\frac{\sqrt{2}}{2}
$$}
Consequently $\theta_1<\frac{\pi}{4}$. Similarly we get $\theta_2<\frac{\pi}{4}$.
Let $z_i$ be the intersection of $\s{S}$ with the ray at $p$ with direction $q_i$.
Analogously let $q'$ be the intersection of $\s{S}$ with the ray at $p$ with direction $q$. Note that $z_i\in\Delta'_i$. Also note that
$d_\s{S}(q,z_i)$ is equal to the length of the arc
$qz_i$ in $\s{S}$. This together with the fact that $\s{S}$ is a 
$(m-j-1)$-sphere of radius one imply that $d_\s{S}(q,z_i)=\theta_i$.
Therefore 
{\footnotesize $$\frac{\pi}{2}\,\leq\, d_\s{S}(\Delta'_1,\Delta'_2)\,\leq\,
d_\s{S}(z_1,z_2)\,\leq
d_\s{S}(z_1,q')+d_\s{S}(q',z_2)\,\leq \,\big(\theta_1+\theta_2\big)\,
$$}
Hence $\frac{\pi}{2}\leq \theta_1+\theta_2<\frac{\pi}{4}+\frac{\pi}{4}=\frac{\pi}{2}$ which is a contradiction. This proves the lemma. 
\vspace{.1in}

\noindent {\bf Lemma 5.1.4.} {\it  Let
$\Delta^k_1$, $\Delta^k_2\in\bS^m$ be two different $k$-simplices,
and $\Delta^{k-1}=\Delta^k_1\cap\Delta^k_2$. Moreover the $\Delta^k_i$
span a simplex.
Let $\alpha, \beta\in (0,\pi/2)$. Suppose that }
$\s{N}_{\beta}(\Delta_1^k,\bS^m)\,\cap\, \s{N}_{\beta}(\Delta_2^k,\bS^m)\,\,\,\,
{\mbox{{\Large $\sbs$}}}\,\,\, \s{N}_{\alpha}(\Delta^j,\bS^m).$
\noindent {\it Then $\frac{sin\,\beta}{sin\,\alpha}\leq\frac{\sqrt{2}}{2}$. }
\vspace{.1in}

\noindent {\bf Proof.} The lemma is certainly true for $\bS^1$.
Using the spherical law of sines it is straightforward to verify the lemma for
$\bS^2$. The case $\bS^m$, $m>2$, can be reduced to the case $m=2$ using
the orthogonal sphere $S$ to $\Delta^{k-1}$ in $\bS^m$ at the barycenter of $\Delta^{k-1}$, as in the proof of 5.1.4. This proves the lemma. \vspace{.1in}

The next result says that {\bf DNP} implies a seemingly stronger version of itself (see (5.1.2)').\vspace{.1in}

\noindent {\bf Lemma 5.1.5.} {\it Suppose the pair of sets of widths
$(\sB, \sA)$ satisfies} {\bf DNP}. {\it Let $\Delta^j=\Delta^k\cap\Delta^l$,
$j<k\leq l$. Then}
{\footnotesize $$\s{N}_{\beta_k}(\Delta^k,\bS^m)\,\bigcap\, \s{N}_{\beta_l}(\Delta^l,\bS^m)\,\,\,\,
{\mbox{{\LARGE $\sbs$}}}\,\,\, \,\,\bigcup_{i<k} \s{N}_{\alpha\0{i}}(\Delta^i,\bS^m) $$}
\noindent {\bf Remark.} Note that the condition
$\Delta^j=\Delta^k\cap\Delta^l$, $j<k,l,$\, is equivalent to\,
$\Delta^k\not\sbs\Delta^l$ and $\Delta^l\not\sbs\Delta^k$,
where the empty set is considered a simplex of dimension -1.\vspace{.1in}

\noindent {\bf Proof of Lemma 5.1.5.}
From Proposition 5.1.3 we have $\frac{sin\, \beta_l}{sin\,\alpha_j}\leq\frac{sin\, \beta_k}{sin\,\alpha_j}<\frac{\sqrt{2}}{2}$. The lemma
now follows from Lemma 5.1.4.
% we now get
%$$\s{N}_{\beta_k}(\Delta^k,\bS^m)\,\bigcap\, \s{N}_{\beta_l}(\Delta^l,\bS^m)\,\,\,\,
%{\mbox{{\LARGE $\sbs$}}}\,\,\,  \s{N}_{\alpha_j}(\Delta^j,\bS^m)
%\,\,\,\,{\mbox{{\LARGE $\sbs$}}}\,\,\, \,\,
%\bigcup_{i<k} \s{N}_{\alpha\0{i}}(\Delta^i,\bS^m) $$
%\noindent 
This proves the lemma.
\vspace{.2in}

\noindent {\bf \large 5.2. Natural Neighborhoods on the Sphere $\bS^m$.}
\vspace{.1in}

Let $\Delta\in\bS^m$.
In this section the $\beta$-geometric link at the barycenter of $\Delta$ will be called {\it the linked sphere of $\Delta$ of radius $\beta$}, and will be denoted by
$\s{S}_\Delta$. Rescaling gives an identification between
$\s{S}_{\Delta^k}$ and $\bS^{m-k-1}$, thus we will consider
$\s{S}_\Delta$ as an all-right spherical complex (alternatively
we can consider $\s{S}_{\Delta^k}$ with the angle metric).
In this case the simplices of $\s{S}_\Delta$
are $\s{S}_\Delta\cap \Delta'$, for all $\Delta'\supset \Delta$.\vspace{.1in}

%The sphere of radius $\mu>0$
%will be denoted by $\bS^k_\mu$. The closed neighborhoods
%and linked spheres will be denoted by ${\sf{N}}_{_{\beta}}(\bS^k_\mu,\bS^m_\mu)$ and $\s{S}^\mu_\Delta$, respectively.
%\vspace{.1in}

Let ${\sf{B}}=\{\beta_k\}$ be a set of widths.  
Let $\Delta=\Delta^k\in\bS^m$ and  
$\s{S}_\Delta$ be the linked sphere of $\Delta$ of radius $\beta_k$. By intersecting $\s{S}_\Delta$ with each element of the set 
$\s{N}_\s{B}(\bS^m)$ we get the set  
$\s{N}(\s{S}_\Delta,\s{B})\,=\,\{ {\s{S}_\Delta\,\cap\,\sf{N}} _{_{\beta_j}}(\Delta^j,\bS^m) \}_{\Delta^j\in\bS^m}$.
It is straightforward to verify that for simplices $\Delta^j$, with $\Delta\sbs\Delta^j$, there are decreasing 
$\beta'_{j-k-1}>0$ such that
\vspace{.03in}

\hspace{1.8in}{\small $ \s{S}_\Delta\,\cap\,\s{N} _{_{\beta_j}}(\Delta^j,\bS^m)\,=\, \s{N} _{_{\beta'_{j-k-1}}}(\s{S}_\Delta\cap\Delta^j,\s{S}_\Delta)$}\vspace{.01in}

\noindent where the last term is the  $\beta'_{j-k-1}$-normal neighborhood of
the simplex $\s{S}_\Delta\cap\Delta^j$ in $\s{S}_\Delta$
(recall that we are identifying $\s{S}_\Delta$ with $\bS^{m-k-1}$,
or, using the angle metric).
Note that, since the $\beta_i$'s are decreasing, we have
$\beta'_l<\pi/2$.
Hence we can write $\s{N}(\s{S}_\Delta,\s{B})\,=\,\s{N}_{\s{B'}(m-k-1)}(\s{S}_\Delta)$
where $\s{B'}(m-k-1)=\{\beta'_0,...,\beta'_{m-k-2}\}$ and we also say that $\s{N}_{\s{B}'(m-k-1)}(\s{S}_\Delta)$ is the set of $\s{B'}(m-k-1)$-neighborhoods
of $\s{S}_\Delta$. 
%[The set $\s{B'}$ is not strictly speaking a set of widths because
%some $\beta'_i$ may be $\geq \pi/4$.]
Note that $\s{B}'(m-k-1)$ depends only on $\s{B}$ and the dimension $k$ of $\Delta^k$.
The next lemma gives this relation explicitly.\vspace{.1in}

\noindent {\bf Lemma 5.2.1.} {\it For $l=0,...,m-k-2$
we have \,\,$sin(\beta'_l)=\frac{sin(\beta_{k+l+1})}{sin(\beta_k)}$.}\vspace{.1in}
%$$sin\,\bigg(\,\frac{\beta'_i}{sin\,\beta_k}\,\bigg)\,\,=\,\, \frac{sin(\beta_{k+1+i})}{sin(\beta_k)}$$

\noindent {\bf Proof.} Let $p\in\s{S}_{\Delta}\cap\big(\s{N}_{_{\beta_j}}(\Delta^j,\bS^m)\big)$, where $\Delta=\Delta^k\sbs\Delta^j$. Then there is a $q\in\Delta^j$ such that $d=d_{_{\bS^m}}(p,q)=d_{_{\bS^m}}(p,\Delta^j)\leq\beta_j$. We are interested in the case when $d$ is maximum, so we assume $d=\beta_j$.
Let $o$ be the barycenter of $\Delta$. Since  $d_{_{\bS^m}}(o,p)=\beta_k$
we get a right (at $q$) spherical triangle with one side equal to $\beta_j$ and hypotenuse equal to $\beta_k$. The angle opposite to the side of length $\beta_j$ is $\beta'_{j-k-1}$. Then, by the spherical law of sines we get $sin(\beta'_{j-k-1})=\frac{sin(\beta_j)}{sin(\beta_k)}$. 
This proves the lemma.\vspace{.1in}

Therefore the set of widths $\s{B}$ on $\bS^m$ induces the set  $\s{B'}(m-k-1)$ on $\bS^{m-k-1}$,
by considering $\bS^{m-k-1}$ as the (rescaled) link of the $k$ simplices in the all-right  triangulation of $\bS^m$. Lemma 5.2.1 gives a relationship between $\s{B}(m)$ and $\s{B'}(m-k-1)$.\vspace{.1in}

Let $\s{B}=\{ \beta_{_i}\}_{i=0,1,...}$ be a set of widths. 
We say that $\s{B}$ is a {\it natural set of neighborhood widths for all spheres} if $\s{B}(m-k-1)=\s{B'}(m-k-1)$ for all $m $ and $k$ with $m>k$.\vspace{.1in}

\noindent {\bf Corollary 5.2.2.} {\it The set of widths $\s{B}=\{\beta_{_i}\}$ is natural if and only if
$sin(\beta_i)=sin^{i+1}(\beta_{_0})$ and $\beta\0{0}<\pi/4$.}\vspace{.1in}

\noindent {\bf Proof.} It follows from 5.2.1 with $l=0$ that 
$sin(\beta_{_{k+1}})=sin(\beta_{_{k}})\,sin(\beta_{_0})$. This proves the corollary.\vspace{.1in}

Given $\varsigma\in (0,1)$ we define $\s{B}(\varsigma)=\{ \beta_{_i}\}$ by $\beta_{_i}=sin^{-1}\big( \varsigma^{i+1} \big)$.
%Write $\s{N}_{_\varsigma}(\bS^m)=\s{N}_{\s{B}(\varsigma)}(\bS^m)$.
Hence the corollary says that $\s{B}$ is natural if and only if $\s{B}=\s{B}(\varsigma)$, for some $\varsigma\in (0,1)$. 
In fact, in this case we have $\varsigma=sin(\beta_{_0})$.\vspace{.1in}

Let $\varsigma\in (0,1)$ and $c>1$. We denote by $\sB(\varsigma;c)=\{\gamma_i\}$ the set defined by  $\gamma_i=sin^{-1}\big( c\,\varsigma^{i+1}\big)$.
Note that $\s{B}(\varsigma;c)$ is a set of widths provided $c\varsigma<1$. Proposition 5.1.3 implies the next corollary.\vspace{.1in}

\noindent {\bf Corollary 5.2.3.} {\it The pair of sets of 
widths $\big(\s{B}(\varsigma;c),\s{B}(\varsigma;c')\big)$ satisfy}  {\bf DNP} {\it
provided $\frac{c}{c'}\,\varsigma<\frac{\sqrt{2}}{2}$.}\vspace{.1in}

Note that we can take $c=c'=1$ which implies that a natural set of widths $\s{B}(\varsigma)$
satisfy {\bf DNP} with $\s{A}=\s{B}=\s{B}(\varsigma)$.
\vspace{.2in}

\noindent {\bf 5.3. \large Neighborhoods in Piecewise Spherical complexes.}\vspace{.1in}

This subsection is essentially a version of 5.1 in which we 
replace $\bS^m$ by an all-right spherical complex.
Let $P$ be an all-right spherical complex and
$\Delta^k\in P$. %Then the set of all 
%$\s{Link}(\Delta^k,\Delta^j)$, with $\Delta^k\sbs\Delta^j$ is also an all-right spherical complex which we denote
%by $\s{Link}(\Delta^k, P)$. 
We can write $\sL(\Delta^k,P)=\bigcup_{\Delta^k\sbs\Delta^j\in P}\sL(\Delta^k,\Delta^j)$
as sets and complexes.
Hence the set $\{\Delta^j\}_{\Delta^k\,\sbs\,\Delta^j\in P}$ is in one-to-one correspondence
with the set of spherical simplices of $\sL(\Delta^k,P)$, that is $\Delta^j$ corresponds to $\sL(\Delta^k,\Delta^j)$,
which is an all-right  spherical simplex of dimension $j-k-1$ in $\sL(\Delta^k,P)$. 
The all-right  spherical metric on $\s{Link}(\Delta^k,P)$
will be denoted by $\sigma_{_{\s{Link}(\Delta^k,P)}}$.\vspace{.1in}

%Also, the $\s{Star}(\Delta,P)$ is the subcomplex of $P$ formed by all spherical %simplices in $P$ containing $\Delta$\vspace{.1in}

\noindent {\bf Remark 5.3.1.} In the paragraph above we did not specify
the type of link we were using (see 5.0.1).
Any of these will lead to the corresponding definition of $\sL(\Delta,P)$, but
they are all equivalent as metric complexes (note that we are assuming $P$
has the ``intersection condition").  
We will use any of the definitions depending on the situation.\vspace{.1in}

\noindent {\bf Lemma 5.3.2.} {\it We have that}\vspace{.02in}

\hspace{1.8in}$
\sL\big( \sL (\Delta^k,\Delta^j), \sL(\Delta^k, P) \big)\,\,\,=\,\,\,\sL (\Delta^j, P)
$\vspace{.1in}

\noindent {\it provided $\Delta^k\,\sbs\,\Delta^j$. (This is an equality of 
all-right spherical metric complexes.)}\vspace{.1in}

\noindent {\bf Remark.} If we use the simplicial definition of links this equality
is actually an equality of sets.\vspace{.1in}

\noindent {\bf Proof.} Let $\Delta^k\,\sbs\,\Delta^j$. Let $\Delta^l$ be the opposite face of $\Delta^k$ in $\Delta^j$.
The statement in the lemma written using the simplicial definition of link is:\, $\sL\big( \Delta^l, \sL(\Delta^k, P) \big)\,\,\,=\,\,\,\sL (\Delta^j, P)$.
We prove this. We have that $\Delta^i\in \sL\big( \Delta^l, \sL(\Delta^k, P) \big)$ is equivalent to (recall $\Delta^k\cap\Delta^l=\emptyset$)
(1) $\Delta^i\cap\Delta^k=\emptyset$, (2) $\Delta^i\cup\Delta^k$ is contained in a simplex, 
(3) $\Delta^i\cap\Delta^l=\emptyset$, (4) $\Delta^i\cup\Delta^l\cup\Delta^k$ is contained in a simplex.
On the other hand  $\Delta^i\in \sL (\Delta^j, P)$ is equivalent to
(a) $\Delta^i\cap\Delta^j=\emptyset$, (b) $\Delta^i\cup\Delta^j$ is contained in a simplex. 
Since  $\Delta^l$ is opposite to $\Delta^k$ in $\Delta^j$ we have that (a)+(b) if
and only if (1)-(4). This proves the lemma.\vspace{.1in}

Define
the {\it closed normal neighborhood of $\Delta^k$ in $\Delta^j$ of width $\beta$ } as
\,$\s{N}_{_\beta}(\Delta^k,\Delta^j)=\s{N}_{_\beta}(\Delta^k,\bS^m)\cap\Delta^j$. %As mentioned
%at the end of section 6.1, equation (6.0.1) implies that
%$\s{N}_{_\beta}(\Delta^k,\Delta^j)-\Delta^k$ is isometric to $\Delta^k\times\sd{j-k-1}\times (0,\beta)$
%with the doubly warp metric of (6.0.1).
If $\Delta^k$ is a simplex in the all-right  spherical complex $P$, we define
the {\it closed normal neighborhood of $\Delta^k$ in $P$ of width $\beta$ } as%\vspace{.01in}

\hspace{2in}$\s{N}_{_\beta}(\Delta^k,P)=\bigcup_{\Delta^k\,\sbs\,\Delta^j\in P}\s{N}_{_\beta}(\Delta^k,\Delta^j)$\vspace{.1in}

%\noindent and we can identify $\s{N}_{_\beta}(\Delta^k,P)-\Delta^k$ with
%$\Delta^k\times\s{Link}(\Delta^k,P)\times (0,\beta)$ with the doubly warp metric\vspace{.1in}

%\noindent {\bf (5.3.3.)}\hspace{1in}
%$ cos^2(r)\, \sigma_{_{\bS^k}}\,+\,sin^2(r)\,\sigma_{_{\s{Link}(\Delta^k,P)}}\,+\, dr^2$\vspace{.1in}

%\noindent where $r$ is the distance to $\Delta^k$, and $\sigma_{_{\s{Link}(\Delta^k,P)}}$
%is the spherical all-right  metric on $\s{Link}(\Delta^k,P)$.\vspace{.1in}

Let $\s{B}=\{\beta_k\}$ be a set of widths. Then, for any all-right  spherical complex $P$
the set $\s{B}$ induces the set of
neighborhoods $\s{N}_{\sB}(P)=\big\{ \s{N}_{_{\beta_k}}(\Delta^k,P) \big\}_{\Delta^k\in P}$. The next lemma is the spherical complex version
of Lemma 5.1.4.\vspace{.1in}

\noindent{\bf Corollary 5.3.4.} {\it  Let
 $\Delta^k$, $\Delta^l\in P$
and $\Delta^j=\Delta^k\cap\Delta^l$. Let $\alpha, \beta, \gamma\in (0,\pi/2)$
such that \,$\frac{sin\,\beta}{sin\,\alpha}, \,\frac{sin\,\gamma}{sin\,\alpha}\,\leq\,\frac{\sqrt{2}}{2}$. Then }
\,$\s{N}_{\beta}(\Delta^k,P)\,\cap\, \s{N}_{\gamma}(\Delta^l,P)\,\,\,\,
{\mbox{{\Large $\sbs$}}}\,\,\, \s{N}_{\alpha}(\Delta^j,P) $.\vspace{.1in}

\noindent {\bf Proof.} The proof is the same as the proof of Lemma 5.1.5
with minor obvious changes. Just recall that we are assuming $P$ to have the
intersection condition. Also note that if $p$, $q$, $q_i$ are as in the proof of
5.1.5 then they all three lie in an all-right simplex in $P$. This proves the
corollary.\vspace{.1in}

As in Section 5.1, the next two results follow directly from corollary 5.3.4. 
The first is a version of {\bf DNP} (see 5.1.2) for $P$, obtained by replacing $\bS^m$ by $P$.\vspace{.1in}

\noindent {\bf Corollary 5.3.5.} {\it Let the pair of sets of widths $\big(\s{B},\s{A}\big)$
satisfy } {\bf DNP}. {\it Then for any all-right spherical complex $P$ and
$k$ we have that the following sets are disjoint}
{\footnotesize $$\bigg\{ \, \s{N}_{\beta_k}(\Delta^k,P)\, \,-\, \,\bigcup_{j<k} \s{N}_{\beta_j}(\Delta^j,P) \,\bigg\}_{\Delta^k\in P}$$}
\noindent The next is a version of Lemma 5.1.5 for general $P$.\vspace{.1in}

\noindent {\bf Corollary 5.3.6.} {\it Let the pair of sets of widths $\big(\s{B},\s{A}\big)$ satisfy } {\bf DNP}. {\it Then for any all-right spherical complex $P$ and
$\Delta^j=\Delta^k\cap\Delta^l$, $j<k\leq l$, simplices in $P$, we have that}
{\footnotesize $$\s{N}_{\beta_k}(\Delta^k,P)\,\bigcap\, \s{N}_{\beta_l}(\Delta^l,P)\,\,\,\,
{\mbox{{\LARGE $\sbs$}}}\,\,\, \,\,\bigcup_{i<k} \s{N}_{\alpha\0{i}}(\Delta^i,P) $$}%\vspace{.2in}

\begin{center} {\bf \large  Section 6. Normal Neighborhoods on Hyperbolic Cones}\end{center} 

In this section we define and give some properties of  neighborhoods of faces in hyperbolic cones. Hyperbolic cones are
cones over all-right spherical complexes; they admit a
canonical piecewise hyperbolic metric. To define the neighborhoods on hyperbolic cones we will use
the objects and results of Section 5.\vspace{.1in}

\noindent {\bf \large 6.1. Neighborhoods in Piecewise Hyperbolic cones.}\vspace{.1in}

We write 
$\R^{k+1}_+=(0,\infty)^{k+1}$,
$\bar{\R}^{k+1}_+=[0,\infty)^{k+1}$ and $\bar{\HH}^{k+1}_+=\B_{_{\HH}}^{k+1}\cap\bar{\R}^{k+1}_+$,
where $\B_{_{\HH}}^{k+1}$ is the disc model of $\HH^{k+1}$.
The canonical all-right spherical $k$-simplex is
$\Delta\0{\bS^k}=\bS^k\cap\bar{\R}^{k+1}_+$. We denote
the origin of $\HH^{k+1}$ by $o=o\0{\HH^{k+1}}$.
We can identify 
$\bar{\HH}^{k+1}_+-\{ o\}$ 
with $\sd{k}\times\R^+$ with metric $sinh^2(s)\,\sigma\0{\bS^k}+ds^2$, where $s$ is the distance to the ``vertex" $o$.
We say that $\bar{\HH}^{k+1}_+$ is the {\it infinite hyperbolic cone of} $\sd{k}$ and write $\rC\sd{k}=\bar{\HH}^{k+1}_+$.\vspace{.1in}

Let $P$ be an all-right  spherical complex. 
The piecewise spherical metric on $P$ will be denoted by
$\sigma\0{P}$. Recall that $P$ is constructed by
gluing the all-right spherical simplices $\Delta\in P$, where
each $\Delta=\Delta^k$ is a copy of $\sd{k}$. The {\it infinite piecewise hyperbolic cone of $P$} is the space  $\rC P$ 
obtained by gluing the hyperbolic cones $\rC\Delta$, $\Delta\in P$ using the same rules used for obtaining $P$.
Note that all vertex points of the $\rC \Delta$ get glued to a unique {\it vertex} \, $o=o_{_{\rC P}}$.
The cones $\rC \Delta$, $\Delta\in P$, are the {\it cone simplices of}\, $\rC P$ and the
{\it faces } of the cone simplex $\rC \Delta$ are the $\rC \Delta'$, with $\Delta'\sbs\Delta$. The set of all cone simplices 
will also be denoted by $\rC P$. 
The complex $\rC P$ (i.e. $\rC P$ together with its cone faces) is an {\it all-right hyperbolic cone
complex}. \vspace{.1in}

The piecewise hyperbolic metric on $\rC P$ shall be denoted by 
$\sigma_{_{\rC P}}$ and its corresponding geodesic metric by
$d_{_{\rC P}}$. Note that the metric $\sigma\0{\rC P}$ can be deduced from the hyperbolic cone structure. 
Note also that $\rC P$ is smooth and hyperbolic
outside the cone of the codimension 2 skeleton of $P$.
All (constant speed) rays emanating from $o$ are length minimizing geodesics defined on $[0,\infty)$.
Then we can identify $\rC P-\{ o\}$ with $P\times\R^+$ with warp metric $sinh^2(r)\,\sigma_{_P}+dr^2$, where $r$ is the distance
to the vertex $o$.
Even though $\sigma\0{\rC P}$ is not (generally) smooth, the set 
of rays emanating from the vertex $o\0{\rC P}$ gives  a well defined
ray structure as in Section 1.\vspace{.1in}

%For $\Delta\in P$, $\beta>0$ and
%$x\in \rC\dDelta-o$, the {\it $\beta$-link $\sL_\beta(\rC \Delta, \rC P)$ of $\rC\Delta$
%in $\rC P$, at $x$,} is the set of all end points of geodesics of length 
%$\beta$ emanating
%perpendicularly to $\rC\Delta$ at $x$. As in section 5.3, by rescaling 
%$\sL_\beta(\rC \Delta, \rC P)$ has a natural all-right spherical structure.\vspace{.1in}

For $s\geq 0$ we denote the {\it open ball of radius s of} $\rC P$ {\it centered at o } by $\B_{_s}(\rC P)$. Note that
this ball is the ``finite open cone" $P\times (0, s)\cup\{o\}$, where we are using the identification above.
The {\it closed ball} will be denoted by ${\bar{\B}}_{_s}(\rC P)$ and the {\it sphere of radius s}, $s>0$, will be denoted
by $\bS_s(\rC P)$, which we shall sometimes identify with $P\times\{ s\}$ or simply with $P$.\vspace{.1in}

Let $\Delta\in P$. In this section $\s{Star}(\Delta,P)$ will
denote the \s{simplicial} star of $\Delta$ in $P$.  
Since $\s{Star}(\Delta,P)$ is an all-right spherical complex 
then $\rC(\s{Star}(\Delta,P))$ is a well defined all-right hyperbolic
cone complex, which we could interpret as the
the simplicial star of $\rC\Delta$ in 
$\rC P$. To save parentheses we will write
$\rC\s{Star}(\Delta,P)$ instead of $\rC(\s{Star}(\Delta,P))$.
Note that $\rC\s{Star}(\Delta,P)\sbs \rC P$.\vspace{.1in}

% as the complex $\s{Star}(\rC \Delta,\rC P)=\{\rC\Delta'\}_{\Delta'\in\s{Star}(\Delta,P)}$. We also write $\s{Star}(\rC \Delta,\rC P)$
%for the underlying space, i.e. the union of all
%$\rC\Delta'$, $\Delta'\in\s{Star}(\Delta,P)$.\vspace{.1in}

 We will use the following three identifications, given
 in 6.1.2, 6.1.3, and 6.1.4 below. (Notation warning:
 sometimes we consider $\Delta^k\sbs \Delta^j$ and
 sometimes $\Delta^j\sbs \Delta^k$.)
\vspace{.1in}

\noindent {\bf 6.1.2.}
Let $\Delta^j\sbs\Delta^k$ and let $\Delta^l$ be the opposite face of $\Delta^j$  in $\Delta^k$. We have $l=k-j-1$. Since $\rC\Delta^j=\bar{\HH}_+^{j+1}\sbs\HH^{j+1}$, $\rC\Delta^l=\bar{\HH}_+^{l+1}\sbs\HH^{l+1}$
and $\rC\Delta^k=\bar{\HH}_+^{k+1}=\bar{\HH}_+^{j+1}\times \bar{\HH}_+^{l+1}
\sbs\HH^{k+1}$ we can write\vspace{.1in}

\hspace{1.6in}$\rC\Delta^k=\rC\Delta^j\times\rC\Delta^l
\sbs\HH^{j+1}\times\HH^{l+1}=\cE_{k+1}\HH^{l+1}$\vspace{.05in}

\noindent  with warp metric $cosh^2(r)\,\sigma_{_{\HH^{j+1}}}+\sigma_{_{\HH^{l+1}}}$,
where $r$ is the distance in $\HH^{l+1}$ to $o$.
Thus we can write $\rC\Delta^k=\cE_{\rC\Delta^j}(\rC\Delta^l)$.
Note that the order of the decomposition here is important (see 2.5). The identification above can be done explicitly in the 
following way. Let $p\in\rC\Delta^k\sbs\HH^{k+1}=\cE_{\HH^{j+1}}(\HH^{l+1})$. 
We use the functions (or coordinates) given in Section 2: $s$, $r$, $t$, $y$, $v$, $x$,
$u$, $w$. Then $p=sx\in\rC\Delta^k$, $(s,x)\in\R^+\times\Delta^k$, corresponds to
$(y,v)=(tw,ru)\in\rC\Delta^j\times\rC\Delta^l$, $(t,w)\in\R^+\times\Delta^j$,
$(r,u)\in\R^+\times\Delta^l$. Note that $x=[w,u](\beta)$, where
$\beta$ is as in Section 2, i.e. it is the angle between $w$ and $x$, and $[w,u]$ is the spherical segment in $\Delta^j$
from $w$ to $u$.\vspace{.1in}

%\noindent {\bf 4.} If $\rC P$ is an all-right hyperbolic cone complex
%and $\Delta$ is an spherical simplex then $\rC\Delta\times \rC P$ is naturally
%an all-right hyperbolic cone complex. In fact the hyperbolic cone structure
%is given by the identification
%$$
%\rC\Delta\times\rC P\,=\, \rC\big(\,    \Delta\,\ast\, P      \,\big)
%$$
%\noindent where $\Delta\,\ast\, P$ is the all-right spherical complex %whose
%simplices are the joins $\Delta\ast\Delta'$, $\Delta'\in P$. The identification
%above is done using identification 6.1.2.(1) simplexwise. (Note that $\Delta'$ is the opposite face of $\Delta$ in $\Delta\ast\Delta'$.)\\ 

\noindent{\bf 6.1.3.} Using 6.1.2 simplexwise we obtain
the following important identification;
this identification will be used many times. 
Let $\Delta=\Delta^k\in P$.
We have that $\rC\s{Star}(\Delta,P) $ can be identified with
$\rC\Delta\times \rC\big(\sL(\Delta, P)\big)$ with metric $cosh^2(r)\,\sigma_{_{\HH^{k+1}}}\,+\, \sigma_{_{\rC\sL(\Delta, P)}}$,
where $r$ is the distance in $\rC\big(\sL(\Delta, P)\big)$ to the vertex $o\in \rC\big(\sL(\Delta, P)\big)$. 
Note that the vertex of $\rC\s{Star}(\Delta,P)$ is identified with
$(o',o'')$, where $o'$, $o''$ are the vertices of $\rC\Delta$ and $\times \rC\big(\sL(\Delta, P)\big)$, respectively.
The identification here is an identification of
all-right hyperbolic cone complexes. Explicitly using the coordinates $s$, $r$, $t$, $y$, $v$, $x$, $u$, $w$ given in Section 2 we see that an element $p=sx\in\rC\s{Star}(\Delta^k,P)$ can be written as $(tw,ru)\in\rC\Delta^k\times\rC\Delta^l\sbs\rC \Delta^k\times\rC\sL(\Delta^k,P)$, where $\Delta^k\ast \Delta^l$ is a simplex
in the (simplicial) star $\s{Star}(\Delta^k,P)$; that is, $\Delta^l$ is a simplex
in $\sL(\Delta^k,P)$. Since we can write $x=[w,u](\beta)$, $\beta$ is the angle
between $w$ and $x$, the identification is given by
$s\,\big[w,u\big](\beta)\,=\, \big( t\,w\,,\,r\, u \big)$.
\vspace{.1in}

\noindent{\bf 6.1.4.} As mentioned above, even though $\sigma_{_{\rC\sL(\Delta, P)}}$
is not in general smooth it has a well defined ray structure, where
we are taking $o\0{\rC P}=o\0{\rC\sL(\Delta, P)}$ as the {\it center} of $\rC\sL(\Delta, P)$.
Hence it makes sense to consider, as in Section 2, the {\it hyperbolic
extension} $\cE_k(\rC\sL(\Delta, P))=\rC\Delta\times \rC\big(\sL(\Delta, P)\big)$ with metric $cosh^2(r)\,\sigma_{_{\HH^{k+1}}}\,+\, \sigma_{_{\rC\sL(\Delta, P)}}$. Therefore, using 6.1.3,  we can write 
\vspace{.03in}

\hspace{1in}$\rC\s{Star}(\Delta,P)=\cE\0{\rC \Delta}\Big(\rC\big( \sL(\Delta,P)\big)\Big)
\,\,\sbs\,\, \cE\0{k}\Big(\rC\big( \sL(\Delta,P)\big)\Big)$\vspace{.05in}

\noindent where we consider $\rC\s{Star}(\Delta,P)\sbs \rC P$ with metric $\sigma\0{\rC P}$
and $\rC\big(\sL(\Delta,P)\big)$ with metric $\sigma\0{\rC\,\sL(\Delta,P)}$.\vspace{.13in}

For a cone simplex $\rC\Delta\in \rC P$, we define its {\it closed normal neighborhood of width $s$} by \vspace{.1in}

\noindent ${\bf (6.1.4.)} \hspace{1in}
\sN_{_s}(\rC\Delta,\rC P)\,=\, \rC\Delta\times{\bar{\B}}_{_s}(\rC \sL(\Delta,P))\,\,\sbs\,\, 
\rC\s{Star}(\Delta,P)
$\vspace{.1in}

\noindent where we are using the identification given in 6.1.3.
Hence $\sN_{_s}(\rC\Delta,\rC P)$ is the union of (the images of)
all geodesics of length $s$ emanating perpendicularly from $\rC \Delta$. The {\it open normal neighboorhood of width s} will be denoted by
$\stackrel{\circ}{\sN}\0{s}(\rC\Delta,\rC P)$.
\vspace{.1in}

\noindent {\bf Lemma 6.1.5.} {\it \, Let  $\Delta^j\sbs\Delta^k\in P$. Then}\vspace{.07in}

\noindent \hspace{.2in}
$
\rC\s{Star}(\Delta^k, P)\,=\,\rC\Delta^k\times\rC\sL(\Delta^k,P)=
\rC\Delta^j\times\rC\s{Star}(\Delta^l,\sL(\Delta^j,P))\sbs
\rC\Delta^j\times\rC\sL(\Delta^j,P)
$\vspace{.05in}

\noindent  {\it where $\Delta^l=\sL(\Delta^j,\Delta^k)$
is the opposite face of $\Delta^j$ in $\Delta^k=\Delta^j\ast\Delta^l$.}
%A similar statement holds if we replace $\sN$ by $\stackrel{\circ}{\sN}$.}
\vspace{.1in} 

\noindent {\bf Proof.}
The first equality is given in 6.1.3 above. The last 
inclusion follows from the inclusion $\s{Star}(\Delta^l,\sL(\Delta^j,P))\sbs\sL(\Delta^j,P)$.
The middle equality in the statement of the lemma is an equality of hyperbolic 
cone complexes.
We have\vspace{.05in}

\noindent \hspace{.4in}$\rC\Delta^k\times\rC\sL(\Delta^k,P)\,=\,
\rC\Delta^j\times\Big(\rC\Delta^l\times\rC \sL(\Delta^k,P)\Big)\,=
\rC\Delta^j\times\rC\s{Star}(\Delta^l,\sL(\Delta^j,P))
$\vspace{.03in}

\noindent where the first equality follows from 6.1.2
(interchange $j$ and $k$). The second equality 
follows from 5.3.2 and 6.1.3 above,
and the fact that $\Delta^l=\sL(\Delta^j,\Delta^k)$. This proves the lemma.\vspace{.1in}

Here is a metric version of Lemma 6.1.5.
Let $\Delta^j$, $\Delta^k$, and $\Delta^l$ be as in Lemma 6.1.5. Let $h:\bS^{m-k-1}\ra\sL(\Delta^k,P)$
be a homeomorphism
and consider the cone of $h$, $\rC h:\R^{m-k}\ra\rC\sL(\Delta^k,P)$.
Let $f'$ be a metric on $\R^{m-k}$ of the form $f'=f'_r+dr^2$. Thus $f'$ and $\sigma\0{\R^{m-k}}$
have the same ray structure. 
The metric $f=h_*f'$ is a metric on
$\rC\sL(\Delta^k,P)$, and it has the same ray structure as 
$\sigma\0{\rC\sL(\Delta^k,P)}$.
We can consider the (restriction of the) metric $\cE_k(f)$ defined on
$\cE_k(\rC\sL(\Delta^k,P)) $ to
$\cE_{\rC\Delta^k}(\rC\sL(\Delta^k,P))=\rC\Delta^k\times\rC\sL(\Delta^k,P)$.
And, since we have $\sL(\Delta^k,P)=\sL(\Delta^l,\rC\sL(\Delta^j,P))$ (see 5.3.2) 
the metric $f$ is also a metric on $\rC\sL(\Delta^l,\sL(\Delta^j,P))$, and we
can consider the metric $\cE_j(\cE_l(f))$ on $\cE_{\rC\Delta^j}\big(\cE_{\rC\Delta^l}(\rC\sL(\Delta^l,\sL(\Delta^j,P))\big)=\rC\Delta^j\times\rC\Delta^l\times
\rC\sL(\Delta^l,\sL(\Delta^j,P))$.\vspace{.1in}

\noindent {\bf Corollary 6.1.6.} {\it Using the identification in 6.1.5 we get}  \,\,$\cE_k(f)\,=\,\cE_j\big(\cE_l(f))
$\vspace{.1in}

\noindent {\bf Proof.} The proof follows from Proposition 2.5 and
the proof of Lemma 6.1.5. This proves the corollary.\vspace{.1in}

We can allow $f$ above to be non-smooth, e.g. we can take $f=\sigma\0{\rC \sL(\Delta^k,P)}$, and we obtain the following corollary.
It follows from 6.1.4 and 6.1.6.\vspace{.1in}

\noindent {\bf Corollary 6.1.7.} {\it Let $\Delta=\Delta^k,\,
\Delta^j,\,\Delta^l$ as in 6.1.5. We have
\vspace{.03in}

\hspace{.31in}$\rC\s{Star}(\Delta,P)\,\,=\,\,\cE\0{\rC \Delta}\Big(\rC\big( \sL(\Delta,P)\big)\Big)
\,\,=\,\,
\cE_{\rC\Delta^j}\Big(\cE_{\rC\Delta^l}\big(\rC\sL(\Delta^l,\sL(\Delta^j,P)\big)\Big)$\vspace{.05in}

\noindent where we consider $\rC\s{Star}(\Delta,P)\sbs \rC P$ with metric $\sigma\0{\rC P}$, $\rC\big(\sL(\Delta,P)\big)$ with metric $\sigma\0{\rC\,\sL(\Delta,P)}$, and 
$\rC\sL(\Delta^l,\sL(\Delta^j,P))$
with metric $\sigma\0{\rC\sL(\Delta^l,\sL(\Delta^j,P))}$.}
\vspace{.1in}

The next two results will be needed in 6.2.\vspace{.1in}

\noindent {\bf Lemma 6.1.8.} {\it \, Let  $\Delta^j\sbs\Delta^k\in P$. Then}\vspace{.05in}

\hspace{1.5in}$
\sN_s(\rC\Delta^k,\rC P)\,=\,\rC\Delta^j\times \sN_s\Big( \rC\Delta^l,\rC\sL(\Delta^j,P) \Big)
$\vspace{.02in}

\noindent {\it where $\Delta^l=\sL(\Delta^j,\Delta^k)$.
A similar statement holds if we replace $\sN$ by $\stackrel{\circ}{\sN}$.}\vspace{.1in}

\noindent {\bf Remark.} Note that the left-hand side of the equality, $\sN_s(\rC\Delta^k,\rC P)$, is a subset of $\rC\s{Star}(\Delta^k,P)$. The right-hand side
is a subset of $\rC\s{Star}(\Delta^l,\sL(\Delta^j,P))
\sbs\rC\Delta^j\times\rC\sL(\Delta^j,P)$. By Lemma 6.1.5 we can
write  $\rC\s{Star}(\Delta^k,P)\sbs\rC\Delta^j\times\rC\sL(\Delta^j,P)$. Lemma 6.1.8
says that under this inclusion $\sN_s(\rC\Delta^k,\rC P)$ corresponds to
$\rC\Delta^j\times \sN_s\big( \rC\Delta^l,\rC\sL(\Delta^j,P) \big)$.\vspace{.1in}

\noindent {\bf Proof.} We have
\vspace{.05in}

\hspace{.2in}$
\sN_s(\rC\Delta^k,\rC P)\,=\,\rC\Delta^j\times\bigg(\rC\Delta^l\times{\bar{\B}}_{_s}(\rC \sL(\Delta^k,P))\bigg)\,=
\,\rC\Delta^j\times \sN_s\bigg( \rC\Delta^l,\rC\sL(\Delta^j,P) \bigg)
$\vspace{.05in}

\noindent where the first equality follows from 6.1.2 and 6.1.4 and the
last from 5.3.2 and 6.1.4. This proves the lemma.\vspace{.1in}

\noindent {\bf Lemma 6.1.9.} {\it \, Let $s>0$, $\beta\in(0,\pi/2)$ and $\Delta\in P$. 
Then}
\vspace{.05in}

\hspace{1.5in}$
\sN_{\,s_{_\beta}}\big(   \rC\Delta, \rC P \big)\,\,\,\,{\mbox{\LARGE $\cap$}}\,\,\,\,\bS_s(\rC P)\,\,\,{\mbox{{\LARGE =}}}
\,\,\,\sN_{_\beta}\big( \Delta, P \big)\times\{ s\}
$\vspace{.05in}

\noindent {\it where $s_{_\beta}=sinh^{-1}\big(sinh (s)\,sin(\beta)  \big)$ and we are identifying $\bS_s(\rC P)$ with $P\times \{ s\}$
(thus $\sN_{_\beta}\big( \Delta, P \big)\times\{ s\}\,\,\sbs\,\, P\times\{ s\}\,\,=\,\,\bS_s(\rC P)$).}\vspace{.1in}

\noindent {\bf Proof.} Denote the vertex of $\rC\big( \sL(\Delta, P) \big)$ by $o'$. Note that both sides of the equality above
are contained in  $\bS_s(\rC P)$. From 6.1.4 and $\beta<\pi/2$ we also get that both sides are contained in $\rC\s{Star}(\Delta,P)$.
 Let $p\in \bS_s(\rC P)$, then $d_{_{\rC P}}(o,p)=s$. From 6.1.4 we can write
$p=(x,y)\in \rC\Delta\times{\bar{\B}}_{_{s'}}(\rC \sL(\Delta,P))$, for some $s'$. Consider the geodesic segments $a=[o,p]$, 
$b=[(x,o'),p]\sbs\{ x\}\times \rC \sL(\Delta,P)$ 
and $c=[o, (x,o')]\sbs\rC\Delta\times\{o'\}$. The length of $a$ is $s$. Since
each $\{x\}\times\rC\sL(\Delta,P)$ is totally geodesic (see 6.1.3) we get that
the length of $b$ is $s'$. Also since $p\in\rC\s{Star}(\Delta,P)$ we have that
$p\in\rC\Delta^j$ for some $\Delta^j\in P$ containing $\Delta$.
But $\rC\Delta^j=\bar{\HH}_+^{j+1}\sbs \HH^{j+1}$ is totally geodesic in
$\rC P$ hence all three segments $a$, $b$, $c$ are contained in $\rC\Delta^j$.
Therefore we get a hyperbolic geodesic triangle with sides
$a$, $b$, $c$, whose angle at $(x,o')$ is $\pi/2$ (because $\rC\Delta\times\{ o'\}$ and $\{ x\}\times \rC\sL(\Delta,P)$ are perpendicular, see 6.1.3).
Let $\beta'$ be the angle at $o$. Then $p\in \sN_{\,s_{_\beta}}\big(   \rC\Delta, \rC P \big)$ if and only if
$s'\leq s_{_\beta}$. Also $p\in \sN_{_\beta}\big( \Delta, P \big)$ if and only if $\beta'\leq \beta$. But
$s_{_\beta}=sinh^{-1}\big(sinh (s)\,sin(\beta)  \big)$ and by the hyperbolic law of sines we also get that 
$s'=sinh^{-1}\big(sinh (s)\,sin(\beta')  \big)$. Consequently
$s'\leq s_{_\beta}$ is
equivalent to $\beta'\leq \beta$. This proves the lemma.
\vspace{.2in}

\noindent {\bf 6.2. Construction of the Fundamental Neighborhoods in Hyperbolic Cones.}\vspace{.1in}

In this section we construct the fundamental sets
$\cY$ and $\cX$ on the cone of a given all-right spherical
complex $P$. These sets depend on a number of pre-fixed 
variables. This subsection is a bit involved and technical,
but the sets $\cY$, $\cX$ are key objects which will be used
in Section 8 to smooth the metric $\sigma\0{\rC P}$ on $\rC P$.
The results that will be used in Section 8 are propositions
6.2.1, 6.2.3, and 6.2.4.
\vspace{.1in}

Let $\xi>0$, $\varsigma\in (0,1)$ and $c>1$  with $c\,\varsigma<e^{-4}$.
Let  $\sB=\sB(\varsigma; c)=\{\beta_i\}$
and $\sA=\sB(\varsigma)=\{\alpha_i\}$
be set of widths as in 5.2. 
We have $sin\,\beta_i=c\,\varsigma^{i+1}$, $sin\,\alpha_i=\varsigma^{i+1}$. 
Since $e^{-4}<\frac{\sqrt{2}}{2}$, the condition $c\,\varsigma<e^{-4}$ together with corollary 5.2.3 imply that
$\big(\s{B},\s{A}\big)$ and $\big(\s{B},\s{B}\big)$
satisfy condition {\bf DNP} in Section 5.1.\vspace{.1in}

Given a number $r>0$ and an integer $k\geq 0$ we define  $r\0{k}=r\0{k}(r)=
sinh^{-1}\big(\frac{sinh(r)}{sin(\alpha\0{k})}\big)$. By convention we also
set $r\0{-1}=r$. (Alternatively we could declare that every set of widths $\{\alpha_k\}$
has a (-1) term $\alpha_{-1}$ always equal to $\pi/2$.)
Let $k$ and $m$ be integers with $m\geq 2$ and $0\leq k\leq m-2$. Define
$s\0{m,k}
=sinh^{-1}\big(\frac{sinh(r)\,sin(\beta\0{k})}{sin(\alpha\0{m-2})}\big)=
sinh^{-1}\big(sinh(r\0{m-2})\,sin(\beta\0{k})\big)$.
We write $r\0{m,k}=r\0{m-k-3}$. Note that $r\0{m,k}<s\0{m,k}$.\vspace{.1in}

Let $P=P^m$ be an all-right  spherical complex with $m\leq \xi$, and let $r>(6+2\xi)$.
%this is used in 6.2.2, hence in 6.2.1.
 For every $\Delta^k\in P$, $0\leq k\leq m-2$, define the following subsets of $\rC P$:\vspace{.1in}

{\small $\begin{array}{lll}
\cY (P, \Delta^k, r,\xi,(c,\varsigma))&=&
\stackrel{\circ}{\sN}\0{s\0{m,k}}(\rC\Delta^k, \rC P)\,-\, \Bigg(\bigcup_{j<k}
\s{N}\0{r\0{m,j}}(\rC\Delta^j,\rC P)\Bigg)
\,-\, \B_{r\0{m-2}-(4+2\xi)}(\rC P)\\\\
\cY (P, r,\xi,(c,\varsigma))&=&
\rC P\,-\, \Bigg(\bigcup_{j<m-1}
\s{N}\0{r\0{m,j}}(\rC\Delta^j,\rC P)\Bigg)
\,-\, \B_{r\0{m-2}-(4+2\xi)}(\rC P)\end{array}$}\\\\

%For $k=m-1, m$\, we define  $\cY(P,\xi,\Delta^m,r)$
%with the same formula as above, with one change: replace
%$\stackrel{\circ}{\sN}\0{s\0{m,k}}(\rC\Delta^k, \rC P)$
%by $\rC\s{Star}(\Delta^k,P)$\vspace{.1in}

Since $\xi$, $c$ and $\varsigma$ will remain constant, in the rest of this section we will write
$\cY(P,\Delta^k,r)$ instead of $\cY(P,\Delta^m,r,\xi,(c,\varsigma))$.\vspace{.1in}

\noindent {\bf Proposition 6.2.1.} {\it For $r>(6+2\xi)$ we have the following properties}
\begin{enumerate}
\item[{\bf (i)}] $\cY(P,\Delta^k,r)
\,\,\,{\mbox{{\Large $\sbs$}}}\,\,\, \stackrel{\circ}{\sN}\0{s\0{m,k}}(\rC\dDelta^k, \rC P)\,\,\,{\mbox{{\Large $\sbs$}}}\,\,\,int\,\rC\s{Star}(\Delta^k,P)$
\item[{\bf (ii)}] {\it $\cY(P,\Delta^k,r)
\,\,\,{\mbox{{\Large $\cap$}}}\,\,\,
\sN\0{r\0{m,j}}(\rC\Delta^j, \rC P)=\emptyset$, for $j< k$.}

\item[{\bf (iii)}] $\cY(P,\Delta^j,r)
\,\,\,{\mbox{{\Large $\cap$}}}\,\,\,
\,\B_{r\0{m-2}-(4+2\xi)}(\rC P)=\emptyset$, 
\item[{\bf (iv)}] $\rC P\,-\,\B_{r\0{m-2}-(4+2\xi)}(\rC P)\,\,=\,\,\cY(P,r)\,\,\cup\,\,
{\mbox{\large $\bigcup$}}\0{\Delta^k\in P,\,k\leq m-2}\cY(P,\Delta^k,r)$
\item[{\bf (v)}] {\it $\Delta^j\cap\Delta^k=\emptyset$ \,\,\,implies\,\,\,
$\sN\0{s\0{m,j}}(\rC \Delta^j,\rC P)\cap\sN\0{s\0{m,k}}(\rC \Delta^k,\rC P)=\emptyset$}
\item[{\bf (vi)}] {\it $\Delta^j\cap\Delta^k=\emptyset$ \,\,\,implies\,\,\,
$\cY(P,\Delta^j,r)\cap\cY(P,\Delta^k,r)=\emptyset$}
\item[{\bf (vii)}] {\it  $\Delta^k=\Delta^i\cap\Delta^j$, \,with $k<i,\,j$,\,\,
implies \,\,, $\cY(P,\Delta^i,r)\cap\cY(P,\Delta^j,r)=\emptyset$}
%\item[{\bf (viii)}]  {\it for any two different $k$-simplices\,\, $\Delta_1^k,\,\Delta_2^k$ \,\,\,we have\,\,\, $\cY(P,\Delta_1^k,r)\cap\cY(P,\Delta_2^k,r)=\emptyset$}
\item[{\bf (viii)}]  {\it $\cY(P,r)
\,\,\,{\mbox{{\Large $\cap$}}}\,\,\,
\sN\0{r\0{m,j}}(\rC\Delta^j, \rC P)=\emptyset$, for $j< m-1$.}
%{\it $\Delta^j\sbs\Delta^k$ \,\,\,implies\,\,\,
%$\cY(P,\Delta^j,r)\cap\cY(P,\Delta^k,r)\,\,\,{\mbox{{\Large $\sbs$}}}\,\,\,
%\rC\s{Star}(\Delta^j,P)\,\,\,{\mbox{{\Large $\sbs$}}}\,\,\,\rC\s{Star}(\Delta^k,P)$}
\end{enumerate}%\vspace{.2in}
\noindent {\bf Proof.}
The statements (ii), (iii), and (ix) follow from the definition of $\cY$. We prove (i). The second inclusion holds because
$\stackrel{\circ}{\sN}\0{s\0{m,k}}(\rC\dDelta^k, \rC P)$ is open.
We prove the first inclusion. By definition
we have
$\cY(P,\Delta^k,r)\sbs\, \stackrel{\circ}{\sN}\0{s\0{m,k}}(\rC\Delta^k, \rC P)$. If  a point $p\in\,\stackrel{\circ}{\sN}\0{s\0{m,k}}(\rC\Delta^k, \rC P)-
\stackrel{\circ}{\sN}\0{s\0{m,k}}(\rC\dDelta^k, \rC P)$ then its distance to $\rC\p \Delta^k$ is
$<s\0{m,k}$. Hence $p\in \,\stackrel{\circ}{\sN}\0{s\0{m,k}}(\rC\Delta^j, \rC P)$ for some $\Delta^j\sbs\p\Delta^k$; thus $j<k$.
But it can be checked that $r\0{m,j}>s\0{m,k}$, $j<k$
(this follows from $c\varsigma<e^{-4}<1$). Therefore
$p\in  \,\stackrel{\circ}{\sN}\0{r\0{m,j}}(\rC\Delta^j, \rC P)$, which
implies
$p\notin \cY(P,\Delta^k,r)$. This proves (i).
Item (ix) follows from (i). Next we prove (iv). Using
$r\0{m,j}<s\0{m,j}$ and the definition
of $\cY(P,r)$ we have\vspace{.1in}

%\hspace{.1in}
 {\small $\rC P\,-\,\B_{r\0{m-2}-(4+2\xi)}(\rC P)\,\,=\,\,\cY(P,r)\,\,\cup\,\,{\mbox{\large $\bigcup$}}\0{j\leq m-2}
\sN\0{r\0{m,j}}(\rC\Delta^j, \rC P)\,\,\sbs\,\,\cY(P,r)\,\,\cup\,\,{\mbox{\large $\bigcup$}}\0{j\leq m-2}
\sN\0{s\0{m,j}}(\rC\Delta^j, \rC P)$}\vspace{.1in}

\noindent This together with (iii)
imply that we can prove (iv) by showing, by induction on $k$, that $U=
{\mbox{\large $\bigcup$}}\0{ l\leq m-2}\cY(P,\Delta^l,r)$ contains
$\stackrel{\circ}{\sN}\0{s\0{m,k}}(\rC\Delta^k, \rC P)\,\,
-\,\,\B_{r\0{m-2}-(4+2\xi)}(\rC P)$
for every $k$-simplex of $P$, $k\leq m-2$. For $k=0$ this statement holds because $\cY(\Delta^0,P)=
\,\,\stackrel{\circ}{\sN}\0{s\0{m,0}}(\rC\Delta^0, \rC P)
-\B_{r\0{m-2}-(4+2\xi)}(\rC P)$. Assume $U$ contains
every   $\stackrel{\circ}{\sN}\0{s\0{m,j}}(\rC\Delta^j, \rC P)\,\,
-\B_{r\0{m-2}-(4+2\xi)}(\rC P)$, for all $j<k$. By the definition of
$\cY(\Delta^k,P)$ we have that
$\stackrel{\circ}{\sN}\0{s\0{m,k}}(\rC\Delta^k, \rC P)
\,-\,\B_{r\0{m-2}-(4+2\xi)}(\rC P)$ is contained in
%\vspace{.05in}

\hspace{1in}
 $\bigg[\cY(\Delta^k,P)\,\,\,\,\,\,\cup\,\,\,\,\,\,
{\mbox{\large $\bigcup$}}\0{ j<k}
\sN\0{r\0{m,j}}(\rC\Delta^j, \rC P)
\bigg]\,\,\,-\,\,\, \B_{r\0{m-2}-(4+2\xi)}(\rC P)$%\vspace{.05in}

\noindent This together with the fact that
$s\0{m,k}>r\0{m,k}$ and the inductive hypothesis imply that
$\stackrel{\circ}{\sN}\0{s\0{m,k}}(\rC\Delta^k, \rC P)
\,-\,\B_{r\0{m-2}-(4+2\xi)}(\rC P)
\,\,\,{\mbox{{\large $\sbs$}}}\,\,\,U$. This proves (iv).
To prove the other two statements we need a lemma.\vspace{.06in}

\noindent {\bf Lemma 6.2.2.}
{\it For $t\geq r\0{m-2}-(4+2\xi)$ and $r>(6+2\xi)$ we have (see Lemma 6.1.9)}\vspace{.05in}

\hspace{1.3in}
$\sN\0{r\0{m,k}}\big(   \rC\Delta, \rC P \big)\,\,\,\,{\mbox{\LARGE $\cap$}}\,\,\,\,\bS_t(\rC P)\,\,\,{\mbox{{\LARGE =}}}
\,\,\,\sN\0{\theta\0{m,k}(t)}\big( \Delta, P \big)\times\{ t\}
$%\vspace{.05in}

\hspace{1.3in}
$\sN\0{s\0{m,k}}\big(   \rC\Delta, \rC P \big)\,\,\,\,{\mbox{\LARGE $\cap$}}\,\,\,\,\bS_t(\rC P)\,\,\,{\mbox{{\LARGE =}}}
\,\,\,\sN\0{\phi\0{m,k}(t)}\big( \Delta, P \big)\times\{ t\}
$\vspace{.05in}

\noindent {\it where $\theta\0{m,k}(t)$ and $\phi\0{m,k}(t)$ are defined by 
the equations $sin(\theta\0{m,k}(t))=c''sin(\alpha\0{k})$, $sin(\phi\0{m,k}(t))=c''sin(\beta\0{k})$, with $c''=\frac{sinh(r\0{m-2})}{sinh(t)}<2e^{2}$.
Moreover $\theta\0{m,k}(t)$ and $\phi\0{m,k}(t)$ are well defined
and less that $\pi/4$.}\vspace{.1in}

\noindent {\bf Proof.} From Lemma 6.1.9 and the definitions of 
$\alpha\0{k}$ and $\beta\0{k}$
we have\vspace{.05in}

\hspace{1in}
$sin (\theta\0{m,k}(t))\,\,=\,\, \frac{sinh(r\0{m,k})}{sinh(t)}\,\,=\,\,
\frac{sinh(r\0{m-2})}{sinh(t)}\,\,\frac{sinh(r\0{m,k})}{sinh(r\0{m-2})}\,\,=\,\,
c''sin(\alpha\0{k})$

\noindent  and %\vspace{.05in}

\hspace{1in}
$sin (\phi\0{m,k}(t))\,\,=\,\, \frac{sinh(s\0{m,k})}{sinh(t)}\,\,=\,\,
\frac{sinh(r\0{m-2})}{sinh(t)}\,\,\frac{sinh(s\0{m,k})}{sinh(r\0{m-2})}\,\,=\,\,
c''sin(\beta\0{k})$

\noindent Since $\xi>0$, simple calculation shows
that $c''<2e^{2}$, 
%$c''<\frac{sinh 4}{sinh 2}$
provided $t\geq r\0{m-2}-(4+2\xi)$, $r>6+2\xi$ (thus $r\0{m-2}>6+2\xi$). Hence
the definitions of $\alpha_k$ and $\beta_k$ at the beginning
of this section imply
$c''sin(\alpha\0{k})\,=\,c''\varsigma^{k+1}<\frac{\sqrt{2}}{2}$
and $c''sin(\beta\0{k})\,=\,c''c\,\varsigma^{k+1}<\frac{\sqrt{2}}{2}$.
This proves the lemma.\vspace{.1in}

We now finish the proof of Proposition 6.2.1. Statement (v) follows from Lemma
6.2.2 and the fact that  $\beta$-neighborhoods, $\beta<\pi/4$, of disjoint
simplices in an all-right spherical complex are disjoint. Statement (vi) follows
from (v). We prove (vii). Note that $c''=c''(m,t)$.
Using items (i), (ii), and  lemmas 6.2.2 and 5.1.5 it is enough to prove that, for fixed $t$ and $m$, 
the pair of sets of widths $\big(\{ \phi\0{m,k}(t) \},\{\theta\0{m,k}(t)  \}\big)$
satisfies {\bf DNP}. But from the definitions we have 
$\{ \phi\0{m,k}(t) \}=\s{B}(\varsigma,c c'')$ and $\{\theta\0{m,k}(t)  \}=\s{B}(\varsigma, c'')$. Therefore Lemma 5.2.3  and the condition
$c\,\varsigma<e^{-4}$ imply
$\big(\{ \phi\0{m,k}(t) \},\{\theta\0{m,k}(t)  \}\big)$
satisfies {\bf DNP}. This proves (vii).
%Statement (viii) follows from (vii) by taking $i=j$.
This proves Proposition 6.2.1.\vspace{.1in}

Define the sets%\vspace{.1in} 

\hspace{1.5in}$\begin{array}{lll}
\cX(P^m,\Delta^k,r)&=&\cY(P^m,\Delta^k,r)-\B_{r\0{m-2}}(\rC P^m)
\\\\\cX(P^m,r)&=&\cY(P^m,r)-\B_{r\0{m-2}}(\rC P^m)
\end{array}$\vspace{.1in}

Alternatively, we can define  $\cX(P^m,\Delta^k,r)$ by the same formula 
that defines  $\cY(P^m,\Delta^k,r)$ with just one small change: in the last term
replace the radius $r\0{m-2}-(4+2\xi)$ by $r\0{m-2}$. Similarly for
 $\cX(P^m,r)$.\vspace{.1in}

\noindent {\bf Proposition 6.2.3.} {\it For $\Delta^j\sbs\Delta^k\in P$ we have}\vspace{.02in}

\hspace{1.5in}$\cY(P,\Delta^k,r)\,\,\,\,{\mbox{\Large $\sbs$}}\,\,\,\,\,\rC\Delta^j\,\,\,\,\times\,\,\,\,
\cX\Big(\sL(\,\Delta^j,P\,), \Delta^l , r \Big)
$\vspace{.03in}

\noindent {\it where $\Delta^l=\Delta^k\cap\sL(\Delta^j,P)$
is opposite to $\Delta^j$ in $\Delta^k$.}
\vspace{.1in}

\noindent {\bf Remark.} The left term in the proposition is a subset of
$\rC\s{Star}(\Delta^k,P)$, thus also a subset of $\rC\s{Star}(\Delta^j,P)$. The right term is a subset $\rC\Delta^j\times
\rC\sL(\Delta^j, P)$ and,  by  6.1.3,  we can write
$\rC\Delta^j\times
\rC\sL(\Delta^j, P)=\rC\s{Star}(\Delta^j,P)$. Proposition 6.2.3 says that 
$\cY(P,\Delta^k,r)$ is a subset of
$\rC\Delta^j\times
\cX\big(\sL\big(\Delta^j,P\big), \Delta^l , r \big)$ under this identification.\vspace{.1in}

\noindent {\bf Proof.} By the (alternative) definition of $\cX$, 
it is enough to prove the following three statements\vspace{.1in}

{\it (1)}\,\,\,\,\,\,\,\,$\cY(P,\Delta^k,r)\,\,\,{\mbox{\Large $\sbs$}}\,\,\,
\rC\Delta^j\,\,\times\,\,
\stackrel{\circ}{\s{N}}\0{s\0{m-j-1,l}}\Big( \rC\Delta^l,\rC\sL(\, \Delta^j,P\,)\Big)$\vspace{.05in}

{\it (2)} \,\,\,\, For $\Delta^i\in\sL(\Delta^j,P)$,
$i<l=k-j-1$, we have\vspace{.05in}

\hspace{1in}$\cY(P,\Delta^k,r)\,\,\,{\mbox{\Large $\cap$}}\,\,\,
\Big[\rC\Delta^j\,\,\times\,\,\s{N}\0{r\0{m-j-1,i}}\big( \rC\Delta^i,\rC\sL(\, \Delta^j,P\,)\big)
\Big]\,\,=
\,\,\emptyset$\vspace{.05in}

{\it (3)}\,\,\,\,
$\cY(P,\Delta^k,r)\,\,\,{\mbox{\Large $\cap$}}\,\,\,
\big[\rC\Delta^j\,\,\times\,\,\B\0{r\0{m-j-3}}\big(\rC\sL(\, \Delta^j,P\,)\big)
\big]\,\,=\,\,\emptyset$\vspace{.1in}

Statement (1) follows from (i) of Proposition 6.2.1, Lemma 6.1.8 and the
equalities $s\0{m,k}=s\0{m-j-1,k-j-1}$ and $l=k-j-1$. 
Statement (2) follows from (ii) of Proposition 6.2.1, Lemma 6.1.8 and the
statements $r\0{m,i+j+1}=r\0{m-j-1,i}$,  $i+j+1<k$. For (3) note that
(6.1.4) and the definition of $r\0{m,j}$ imply\vspace{.05in}

\hspace{1.4in}$\rC\Delta^j\times\B\0{r\0{m-j-3}}\big(\rC\sL(\, \Delta^j,P\,)\big)
\,\,\,\,\,\,\,{\mbox{\Large $=$}}\,\,\,\,\,\,\,
\sN\0{r\0{m,j}}\big(\rC\Delta^j,\rC P\big)
$\vspace{.05in}

This together with (ii) of Proposition 6.2.1 imply (3). This proves the proposition.\vspace{.1in}

\noindent {\bf Proposition 6.2.4.} {\it For $\Delta^k\in P$, $k\leq m-2$, we have}\vspace{.1in}

\hspace{1.3in}$\cY(P,r)\,\,\,\,{\mbox{\Large $\cap$}}\,\,\,\,
\cY(P,\Delta^k,r)\,\,\,\,{\mbox{\Large $\sbs$}}\,\,\,\,\,\rC\Delta^k\,\,\,\,\times\,\,\,\,
\cX\big(\sL(\,\Delta^k,P\,), r \big)
$\vspace{.1in}

%\noindent {\bf Remark.} The left term in the lemma is a subset of
%$\rC\s{Star}(\Delta^k,P)$, thus also a subset of $\rC\s{Star}(\Delta^j,P)$. The right term is a %subset $\rC\Delta^j\times
%\rC\sL(\Delta^j, P)$. Item (5) in remark 6.1.2 says that we can write
%$\rC\s{Star}(\Delta^j,P)=\rC\Delta^j\times
%\rC\sL(\Delta^j, P)$. Lemma 6.1.2 says that 
%\cY(P,\Delta^j,r)\,\cap\,
%$\cY(P,\Delta^k,r)$ corresponds to a subset of
%$\rC\Delta^j\times
%\cX\big(\sL\big(\Delta^j,P\big), \Delta^l , r \big)$ under this identification\vspace{.1in}

\noindent {\bf Proof.} Using the definition of 
%$\cY(P,\Delta^k,r)$, $\cY(P,r)$
$\cX\big(\sL\big(\Delta^k,P\big), r \big)$, it is enough to prove the following three statements\vspace{.1in}

{\it (1)}\,\,\,\,$\cY(P,\Delta^k,r)\,\,\,{\mbox{\Large $\sbs$}}\,\,\,
\rC\Delta^k\,\,\times\,\,
\rC\sL\big( \Delta^k,P\big)$\vspace{.1in}

{\it (2)} For $\Delta^j\in P$, $\Delta^k\sbs\Delta^j$,
$l\leq m-k-3$, 
and $\Delta^l$ opposite to $\Delta^k$ in $\Delta^j$, we have
\vspace{.1in}

\hspace{1.3in}$\cY(P,r)\,\,\,{\mbox{\Large $\cap$}}\,\,\,
\Big[\rC\Delta^k\,\,\times\,\,\s{N}\0{r\0{m-k-1,l}}\big( \rC\Delta^l,\rC\sL(\, \Delta^k,P\,)\big)
\Big]\,\,=
\,\,\emptyset$\vspace{.1in}

{\it (3)}
$\cY(P,r)\,\,\,{\mbox{\Large $\cap$}}\,\,\,
\big[\rC\Delta^k\,\,\times\,\,\B\0{r\0{m-k-3}}\big(\rC\sL(\, \Delta^k,P\,)\big)
\big]\,\,=\,\,\emptyset$\vspace{.1in}

Statement (1) follows from (i) of Proposition 6.2.1,  and 6.1.3.
Statement (2) follows 6.2.1 (viii), Lemma 6.1.8, the identities
$r\0{m-k-1,j-k-1}=r\0{m,j}$, $k+l+1=j$, the fact that $l\leq m-k-3$ if and only if
$j\leq m-2$,
and the definition of $\cY(P,r)$. Finally (3) 
follows from 6.2.1 (viii), the definition of $r\0{m,k}$
and the definition of  $\cY(P,r)$. This proves the lemma.
\vspace{.2in}

\noindent {\bf \large 6.3. Radial Stability of the Sets $\cY(P,\Delta^k,r)$.}
\vspace{.1in}

In Section 8 we will need a sort of a stable property for the sets
$\cY$. We use the objects and notation used
in Section 6.2. Recall that $\s{Star}(\Delta,P)$ is the simplicial star of
$\Delta$ in $P$, and that an element in $\rC P$ can be written as
$sx$, $s\in [0,\infty)$, $x\in P$.
Let $\theta\in(0,\pi/2)$, and write $a(s)=a\0{\theta}(s)=sinh^{-1}(sinh(s)\,sin\,\theta)$.\vspace{.1in}

\noindent {\bf Lemma 6.3.1.} {\it Let $b\in \R$, $\Delta^k\in P$, and
$x\in\s{Star}(\Delta,P)$. Then {\small $(s+b)x\in\s{N}_{a(s)}(\rC \Delta^k,\rC P)$}
if and only if \,\,{\small$
sin(\gamma)\,\frac{sinh(s+b)}{sinh (s)}\,\,\leq\,\,sin\,\theta 
$}, where $s>0$, $\gamma=\gamma(x)=d\0{P}(x,\Delta^k)$.}\vspace{.1in}

\noindent {\bf Proof.} Note that $\gamma$ is the angle opposite to
the cathetus of length $d(s)=d\0{\rC P}((s+b)x,\rC \Delta^k)$ of the right
hyperbolic triangle with hypotenuse $(s+b)$. 
We want $d(s)\leq a(s)$; equivalently $sinh(d(s))\leq sinh(a(s))$. By the hyperbolic law of sines $sinh(d(s))=sin(\gamma)\,sinh (s+b)$, hence $sinh(d(s))\leq sinh(a(s))$ is equivalent to $sin(\gamma)\,sinh (s+b)\leq sinh (s)\, sin\,\theta$. This proves the lemma.\vspace{.1in}

Note that the lemma also holds if we replace
$\s{N}$ by $\stackrel{\circ}{\s{N}}$ and
$\leq $ by $<$.
Write $R(s)=R_{x,b}(s)=(s+b)x$.\vspace{.1in}

\noindent {\bf Lemma 6.3.2.} {\it Let $\Delta^k$, $P$ and $x$ as in Lemma 6.3.1.
We have three mutually exclusive cases:\vspace{.05in}

\noindent \,\,\,\,{\bf C1:}  $e^b sin(\gamma)<sin\,\theta$, which implies that
$R(s)\in\,\,\stackrel{\circ}{\s{N}}_{a(s)}(\rC \Delta,\rC P)$, $s\geq s_0$, for some $s_0$.

\noindent \,\,\,\,{\bf C2:}  $e^b sin(\gamma)>sin\,\theta$, which implies that
$R(s)\notin\s{N}_{a(s)}(\rC \Delta,\rC P)$, for all $s>0$.

\noindent \,\,\,\,{\bf C3:} $e^b sin(\gamma)=sin\,\theta$, which implies that
$R(s)\notin\s{N}_{a(s)}(\rC \Delta,\rC P)$, for all $s>0$.}
\vspace{.05in}

\noindent {\bf Proof.} The lemma follows from 6.3.1 and the following
two facts: (1) the function $s\mapsto\frac{sinh(s+b)}{sinh(s)}$ is strictly
decreasing for $s>0$, and (2) $\lim_{s\ra \infty}\frac{sinh(s+b)}{sinh(s)}=e^b$.
This proves the lemma.\vspace{.05in}

From the definition of $r\0{k}$ given at the beginning of 6.2
we have $r\0{m-2}=r\0{m-2}(r)=sinh^{-1}(\frac{sinh (r)}{sin\,\alpha\0{m-2}})$, hence we can write $r=r(r\0{m-2})=
sinh^{-1}(sinh(r\0{m-2})\,sin\,\alpha\0{m-2})$. Therefore we can write
$r\0{m,k}=r\0{m,k}(r)$ and $s\0{m,k}=s\0{m,k}(r)$ in terms of the new variable $r\0{m-2}$,
and a calculation shows that $r\0{m,k}=a\0{\alpha\0{k}}(r\0{m-2})$
and $s\0{m,k}=a\0{\beta\0{k}}(r\0{m-2})$.
We will use these identities in the proof of the next result. \vspace{.1in}

\noindent {\bf Proposition 6.3.3.}  {\it Fix $b\in \R$ and let $x\in P$.
Then at least one of the following conditions hold.
\vspace{.05in}

\noindent\,\,\,\, (1) There is $\Delta^k$, $k\leq m-2$, such that
{\small $R_{x,b}(r_{m-2})\in\cY(P,\Delta^k,r(r_{m-2}))$}, for all $r_{m-2}>r'$, for some 
$r'$.

\noindent\,\,\,\, (2) We have that
{\small $R_{x,b}(r_{m-2})\in\cY(P,r(r_{m-2}))$}, for all $r_{m-2}>r'$, for some 
$r'$.}\vspace{.05in}

\noindent {\it Moreover, these two conditions are stable in the following sense.
If $x'$ and $b'$ are sufficiently close to $x$ and $b$, respectively,
and $R_{x,b}$ satisfies (i) then $R_{x',b'}$ also satisfies (i) (with the same
$r'$). Similarly for condition (ii).}\vspace{.1in}

\noindent {\bf Proof.} By induction. Suppose {\bf C1} of
6.3.2 holds for $R=R_{x,b}$ with $\theta=\alpha\0{0}$, for some 
$\Delta^0$.
Then, since $\s{N}_{a\0{\alpha\0{0}}(r\0{m-2})}(\Delta^0,P)=\s{N}_{r\0{m,0}(r\0{m-2})}(\Delta^0,P)\sbs\cY(P,\Delta^0,r)$ we see that $R$ satisfies (1) for $\cY(P,\Delta^0,r)$ and we are done.
Suppose {\bf C3} holds with $\theta=\alpha\0{0}$, for some $\Delta^0$. Then $x\in\s{Star}(\Delta^0,P)$ and $e^b sin(\gamma)=sin(\alpha_0)$, where $\gamma=\gamma(x)$. Since $\alpha\0{k}<\beta\0{k}$, we
have $e^b sin(\gamma)<sin(\beta_0)$, hence by 6.3.2 (with
$\theta=\beta\0{0}$) we have that
$R(r\0{m-2})\in\,\stackrel{\circ}{\s{N}}\0{s\0{m,0}(r\0{m-2})}(\Delta^0,P)$, 
for large $r\0{m-2}$, and  follows 
that $R$ satisfies (1) for $\cY(P,\Delta^0,r)$ and we are done.
Now suppose that {\bf C2} happens for all $\Delta^0$, with $\theta=\alpha\0{0}$.
As before we have three possibilities. First {\bf C1} 
holds for $R=R_{x,b}$ with $\theta=\alpha\0{1}$, for some 
$\Delta^1$.
This, together with the assumption that {\bf C2} holds for all
$\Delta^0$ (with $\theta=\alpha\0{0}$), and the definition of
$\cY(P,\Delta^1,r)$ imply that $R$ satisfies (1) for $\cY(P,\Delta^1,r)$ and we are done.
Suppose {\bf C3} holds for $R$ and $\Delta^1$
(with $\theta=\alpha\0{1}$), for some $\Delta^1$. Using the same argument as in the $\Delta^0$ case 
(when we assumed {\bf C3} some $\Delta^0$)
we get that $R$ satisfies (1) for $\cY(P,\Delta^1,r)$ and we are done.
The third case is that {\bf C2} happens for $R$ and all  $\Delta^1$. Proceeding in this way we obtain that
either $R$ satisfies (1), for some $\Delta^k$, $k\leq m-2$ or 
{\bf C2} holds for $R$ and all  $\Delta^k$, $k\leq m-2$
(with $\theta=\alpha\0{k}$). Hence (2)
holds for $R$. Moreover it does so stably. This proves the proposition.
\vspace{.2in}

\begin{center} {\bf \large 7. Smooth Structures on Cube
and All-Right
Spherical Complexes.}
\end{center}

For the basic definitions and results about cube complexes see
for instance \cite{BH}. Given a (cube or all-right spherical) complex $K$ we use the same notation $K$ for the complex itself (the collection of all
closed cubes or simplices) and its realization (the union of all cubes or simplices).
For $\sigma\in K$ we denote its interior by $\dsigma$.\vspace{.1in}

Let $M^n$ be a smooth manifold of dimension $n$. A {\it smooth cubulation} of  $M$ is a pair $(K,f)$, where 
$K$ is a cube complex and $f:K\ra M$ is a non-degenerate
$PD$ homeomorphism \cite{MunkresLectures}, that is, for all $\sigma\in K$  we have $f|_{\sigma}$ is a smooth embedding. Sometimes we will write $K$ instead
of $(K,f)$. The smooth manifold $M$ together with a smooth cubulation is
a {\it smooth cube manifold} or a {\it smooth cube complex}. A {\it smooth all-right-spherical triangulation}  and a {\it smooth all-right-spherical manifold} (or complex) is defined analogously.\vspace{.1in}

%Note that if $K$ is a smooth cubulation (or all-right
%spherical triangulation) of $M$, then $K\cong_{PL}M$, that is,
%$K$ is $PL$-homeomorphic to the smooth manifold $M$.\vspace{.1in}

In this section $\sL(\sigma^{j},K)$ means the {\it geometric link} 
of an open $j$-cube or $j$-all-right simplex $\sigma^j$, defined as
the union of the end points of straight (geodesic) segments of small length $\epsilon>0$
emanating perpendicularly (to $\dsigma^j$) from some point  $x\in \dsigma^j$. 
%We say that the link is {\it based at $x$.} 
The star $\s{Star}(\sigma,K)$ as the union of such segments. We can identify the star with the cone of the link
$\rC\sL (\sigma,K)$  (or $\epsilon$-cone) defined as
$\rC \, \sL(\sigma, K)=\sL(\sigma,K)\times [0,\epsilon)\,/\, \sL(\sigma,K)\times\{0\}$. Thus a point $x$ in $\rC \, \sL(\sigma, K)$, different from the cone point $o=o\0{\rC \, \sL(\sigma, K)}$, can be written as
$x=t\,u$, $t\in (0,\epsilon)$, $u\in \sL(\sigma, K)$. For $s>0$ we get the
{\it cone homothety} $x\mapsto sx=(st)u$ (partially defined if $s>1$). 
If we want to make explicit the dependence of the link or the cone on $\epsilon$ we shall write $\sL_\epsilon(\sigma,K)$ or $\rC _\epsilon\,\sL(\sigma,K)$ respectively. 
\vspace{.1in}

\noindent {\bf Remark 7.0.1.} As usual we shall identify the $\epsilon$-neighborhood
of $\dsigma$ in $K$ with $\rC _\epsilon\,\sL(\sigma,K)\times
\dsigma$ (or just $\rC \,\sL(\sigma,K)\times
\dsigma$). %Hence a cone homothety induces a {\it neighborhood homothety}
%obtained by crossing it with the identity $1\0{\dsigma}$.
%Note  that $\sL(\sigma,K)$ and $\rC \,\sL(\sigma,K)$ are subsets of $K$. 
\vspace{.1in}

In what follows we assume that $f:K\ra M$ is a smooth cubulation
(or all-right spherical triangulation) of the smooth manifold $M$.
Since the $PL$ structure on $M$ induced by $K$ is Whitehead compatible with $M$ we have that
the link $\sL(\sigma^i,K)$ is $PL$ homeomorphic to $\bS^{n-i-1}$.
A {\it link smoothing for}  $\dsigma^i$ (or $\sigma^i$)
is just a homeomorphism $h_{\sigma^i}:\bS^{n-i-1}\ra\sL(\sigma^i,K)$.
The {\it cone} of $h_{\sigma^i}$ is the map 
$\rC \,h_{\sigma^i}:\D^{n-i}\longrightarrow \rC\sL (\sigma^i, K)$
given by $t\,x= [x,t]\mapsto t\,h_{q^i}(x)=[h_{q^i}(x),  \,t]$,
 where we are canonically identifying the $\epsilon$-cone  of $\bS^{n-i-1}$ with the disc $\D^{n-i}$. We remark that we are not assuming $h_{\sigma^i}$ to be smooth. A link smoothing $h_{\sigma^i}$  induces the following  smoothing of the normal neighborhood of
$\dsigma^i$:
$$h^\bullet_{\sigma^i}=f\,\, \circ\,\,\Big(\rC \,h_{\sigma^i}\times 1_{\dsigma^i}\Big):\D^{n-i}\times \dsigma^i\longrightarrow M $$

The pair $(\,h^\bullet_{\sigma^i}\, ,\,\D^{n-i}\times \dsigma^i\,)$, or simply 
$h^\bullet_{\sigma^i}$, is a {\it normal
chart} on $M$. Note that the collection  $\cA=\big\{\,(\,h^\bullet_{\sigma^i}\, ,\,\D^{n-i}\times 
\dsigma^i\,)\,\big\}_{\sigma^i\in K}$ is a topological atlas for $M$.
Sometimes will just write $\cA=\big\{\,h^\bullet_{\sigma^i}\,\big\}_{\sigma^i\in K}$.
The topological atlas $\cA$ is called a {\it normal atlas}. It depends uniquely on the
the complex $K$, the map $f$ and the collection of link smoothings 
$\{h_{\sigma}\}_{\sigma\in K}$. 
To express the dependence of the atlas on the set of links smoothings
we shall write $\cA=\cA\big(\{h_{\sigma}\}_{\sigma\in K}\,\big)$
(this is different from $\cA=\big\{\,h^\bullet_{\sigma^i}\,\big\}_{\sigma^i\in K}$,
as written above). The most important feature about these normal
atlases is that they preserve the radial and sphere (link) structure given by $K$.\vspace{.1in}

Note that not every collection of link smoothings induce a smooth atlas. But when 
the induced atlas is smooth we call $\cA$ a {\it normal smooth atlas on $M$ with respect to} $K$ and the corresponding smooth structure $\cS'$
a {\it normal smooth structure on $M$ with respect to $K$}.  
In this case we say that the set of link smoothings
$\{h_{\sigma}\}_{\sigma\in K}$ is {\it smooth}.
The following theorem is proved in \cite{O1}.\vspace{.1in}

\noindent {\bf Theorem 7.1.} {\it Let $M$ be a smooth cube or 
all-right spherical manifold,
with smooth structure $\cS$.
Then $M$ admits a normal smooth structure $\cS'$ diffeomorphic to $\cS$.}\vspace{.1in}

Hence if $M^n$ is a smooth manifold with smooth structure $\cS$ and $K$ is a cubulation (or all-right spherical triangulation) of $M$, then there are link smoothings $ h_{\sigma}$, for all $\sigma\in K$, such that the atlas $\cA=\cA\big(\{h_{\sigma}\}_{\sigma\in K}\,\big)$
is smooth (equivalently, $\{h_{\sigma}\}_{\sigma\in K}$ is smooth). Moreover the normal smooth structure $\cS'$,
induced by $\cA$, is diffeomorphic to $\cS$.\vspace{.2in}

%The following is a corollary of the proof of Theorem 1.2 given
%in \cite{O1} (see Lemma 1.3 in \cite{O1}).\vspace{.1in}

%\noindent {\bf Corollary 1.2.3.} {\it Let $f:K\ra (M,\cS)$
%be a smooth cubulation (or all-right spherical triangulation)
%of the smooth manifold $(M,\cS)$. Let $\cS'$ as in Theorem 1.2
%Then, for every $\sigma\in K$ we have that $f(\sL(\sigma,K))$
%is a smooth submanifold of $(M,\cS')$.}\\

%\noindent {\bf Corollary 1.2.4.} {\it Let $M$, $\cS$ and $\cS'$ as in Theorem 1.2.
%(or its addendum).
%Then $K$ is $PL$-homeomorphic to $(M,\cS')$.}\\
%\vspace{.2in}

\noindent {\bf \large 7.2. Induced Link Smoothings.}\vspace{.1in}

Let $K$ be a cubical or all-right spherical complex. Then the links of $\sigma\in K$ are
all-right-spherical complexes. 
We explain here how to obtain from a given a collection of link smoothings for $K$ (and its corresponding normal atlas and structure) a collection of links smoothings for a link in $K$
(and its corresponding normal atlas and structure).\vspace{.1in}

The all-right-spherical structure on 
$\sL(\sigma,K)$ induced by $K$ has all-right-spherical simplices
$\big\{ \,\tau\,\,\cap\,\, \sL(\sigma,K)\,\,  ,\,\, \tau\in K\big\} .$ 
Note that $\tau\,\,\cap\,\, \sL(\sigma,K)$ is non-empty only when $\sigma\subsetneq \tau$, hence we can write\vspace{.05in}

\noindent\hspace{1.8in}$\sL(\sigma,K)=\big\{ \,\tau\,\,\cap\,\, \sL(\sigma,K)\,\,  ,\,\, \sigma\subsetneq \tau\in K\big\} $\vspace{.05in}

Since $\tau\,\,\cap\,\, \sL(\sigma,K)$ is a simplex in the 
all-right spherical complex $\sL(\sigma,K)$ we can consider its
link $\sL\,\Big(\, \tau\,\cap\,\sL\,(\sigma,K)  , \sL\,(\sigma,K)  \,\Big)$.
By definition we have:\vspace{.05in}

\noindent {\bf (7.2.1)} \hspace{1in}$\sL\,\Big(\, \tau\,\cap\,\sL\,(\sigma,K)  , \sL\,(\sigma,K)  \,\Big)\,=\,\sL\,\big( \, \tau   , K  \,\big)$\vspace{.05in}

\noindent provided we choose the radii and bases of the links
properly. In the formula above  radii and bases
are not specified but the radii are certainly not equal. 
The simple relationship between these radii is given by equation (1)
in the proof of Lemma 1.3 \cite{O1} (or the corresponding one in the spherical case; see Remark 1 after the proof of Lemma 1.3 \cite{O1}). By 7.2.1 we can say that the set of link smoothings $\{ h_\sigma\}_{\sigma\in K}$ for
$K$ induces, just by restriction, a set of link smoothings for $\sL(\sigma, K)$,
$\sigma\in K$. That is, we set
 $h\0{\tau\cap\sL\,(\sigma,K)}=h\0{\tau}$,
 $\sigma\subsetneq \tau\in K$.
The next result is proved in \cite{O2}.
\vspace{.1in}

\noindent {\bf Proposition 7.2.2.} {\it Let $\{h_\sigma\}_{\sigma\in K}$ be a set of
link smoothings on $K$, and let $\sigma^k\in K$.
Assume $\{h_\sigma\}_{\sigma\in K}$ is smooth; that is, the atlas $\cA=\cA\big(\{h_\sigma\}_{\sigma\in K}\big)$ is smooth. 
Let $\cS'$ be the normal smooth structure
on $K$ induced by $\cA$. Then:
\vspace{.05in}

\noindent \,\,\,\,(1) The set of link smoothings $\{h_
{\sigma^i\cap\sL\,(\sigma,K)}\}_{\sigma^k\subsetneq \sigma^i}$ 
for the links of $\sL(\sigma^k,K)$ is smooth; that is,
the 

\noindent \,\,\,\,\,\,\,\,\,\,\,\,\,\,atlas
$\cA_{\sigma^k}=\cA_{\sL(\sigma^k,K)}=\big\{h^\bullet_
{\sigma^i\cap\sL\,(\sigma^k,K)}\big\}_{\sigma^k\subsetneq \sigma^i}$ 
is a smooth normal atlas on $\sL(\sigma^k,K)$. 
 \vspace{.05in}

\noindent \,\,\,\,(2) The link smoothing%\vspace{.02in}

\hspace{1.5in}$h_{\sigma^k}:\bS^{n-k-1}\ra\Big(\sL(\sigma^k,K)\,\,,\,\,
\cS_{\sigma^k}\Big)$%\vspace{.02in}

\noindent \,\,\,\,\,\,\,\,\,\,\,\, is a diffeomorphism. Here $\cS_{\sigma^k}$ is the smooth structure induced by the atlas $\cA_{\sigma^k}$.
\vspace{.05in}

\noindent \,\,\,\,(3)  We have that $\sL(\sigma^k,K)$ is a smooth submanifold of
$(K,\cS')$. Moreover
\vspace{.05in}

\hspace{1in}$\cS'\big|_{\sL(\sigma^k,K)}\,\,=\,\,\cS_{\sigma^k}
$\vspace{.05in}

\noindent \,\,\,\,\,\,\,\,\,\,\,\, where $\cS'\big|_{\sL(\sigma^k,K)}$
denotes the restriction of $\cS'$ to $\sL(\sigma^k,K)$.}

%And, since a link smoothing $h_\sigma$ is the restriction of the embedding $h^\bullet_\sigma$ we get the following corollary.\vspace{.1in} 

%\noindent {\bf Corollary 2.4.} {\it For every $\sigma^k\in K^n$, the link smoothing $h_{\sigma}:\bS^{n-k-1}\ra \big(\sL(\sigma,K),\cS'_\sigma\big)$ is a diffeomorphism.}
\vspace{.2in}

\noindent {\bf \large 7.3. The Case of Manifolds with Codimension Zero Singularities.}
\vspace{.1in}

Here we treat the case of manifolds with a one point singularity.
The case of manifolds with many (isolated) point singularities
is similar. \vspace{.1in}

Let $Q$ be a smooth manifold with a one point singularity $q$, that is
$Q-\{q\}$ is a smooth manifold and there is a topological embedding
$\rC_1 N\ra Q$, with $o\0{\rC N}\mapsto q$, that is a smooth embedding outside 
the vertex $o\0{\rC N}$. Here $N=(N,\cS_N)$ is a closed smooth manifold (with smooth
structure $\cS_N$).
Also $\rC_1 N $ is the (closed ) cone of width 1 and we identify
$\rC_1 N-\{o\0{\rC N}\}$ with $N\times (0,1]$. We write $\rC_1N\sbs Q$.
We say that the {\it singularity $q$ of $Q$ is  modeled on  $\rC N$.}\vspace{.1in}

Assume $(K,f )$ is a smooth cubulation of $Q$, that is

{\bf (i)}
$K$ is a cubical complex.

{\bf (ii)} $f:K\ra Q$ is a homeomorphism. Write $f(p)=q$ and $L=\sL(p,K)$.

{\bf (iii)} $f|_{\sigma}$ is a smooth embedding for every
cube $\sigma$ not containing $p$.

{\bf (iv)} $f|_{\sigma-\{p\}}$ is a smooth embedding for every
cube $\sigma$ containing $p$.

{\bf (v)}  $L$ is $PL$  homeomorphic to $(N,\cS_N)$.
\vspace{.1in}

Many of the definitions and results given before
for smooth cube manifolds still hold (with minor changes)
in the case of manifolds with a one point singularity:%\vspace{.1in}

\begin{enumerate}
\item[{\bf (1)}] A {\it link smoothing} for $L=\sL(p,K)$ (or $p$) is just a
homeomorphism $h_p:N\ra L$.
Since all but one of the links of $K$ are spheres, sets of  link smoothings for $K$
are defined, that is, they are sets of link smoothings for the sphere links plus
a link smoothing for $L$.
\item[{\bf (2)}] Given a set of link smoothings for $K$ we get a set of normal charts
as before. For the vertex $p$ we mean the cone map 
$h_p^\bullet=f\circ\rC h_p:\rC N\ra Q$. We will also denote the restriction
of $h_p^\bullet$ to $\rC N-\{o\0{\rC N}\}$ by the same notation $h_p^\bullet$.
As before $\{h^\bullet_\sigma\}_{\sigma\in K}$ is a {\it (topological) normal atlas on $Q$
with respect to $K$}. The atlas on $Q$ is {\it smooth} if all transition functions are smooth, where for the case $h_p^\bullet:\rC N-\{o\0{\rC N}\}\ra Q-\{q\}$ we are identifying $\rC N-\{
o\0{\rC N}\}$ with $N\times (0,1]$ with the product smooth structure obtained from
\s{some} smooth structure ${\tilde{\cS}}_N$ on $N$.
A smooth normal atlas on $Q$ with respect to $K$ induces, by restriction,  a smooth normal structure on $Q-\{q\}$ with respect to $K-\{p\}$ (this makes sense even though
$K-\{p\}$ is not, strictly speaking, a cube complex).
\item[{\bf (3)}] We say that the set $\{h_\sigma\}$ is {\it smooth} 
if the atlas $\cA=\{h^\bullet_\sigma\}_{\sigma\in K}$ is smooth.  In this case we say that
the smooth atlas  $\cA$ (or the induced smooth structure, or the set $\{h_\sigma\}$) 
is {\it correct with respect to $N$} if $\cS_N$ and ${\tilde{\cS}}_N$ are diffeomorphic.
\item[{\bf (4)}] Also it is straightforward to verify that
Proposition 7.2.2 holds in our present case. 
\item[{\bf (5)}] In \cite{O2} the following version of Theorem 7.1 is
proved:\end{enumerate}

\noindent {\bf Theorem 7.3.1.} {\it  Let $Q$ be a smooth manifold
with one point singularity $q$ modeled on $\rC N$, where $N$ is a closed
smooth manifold. Let $(K,f)$ be a smooth cubulation of $Q$. Then
$Q$ admits a normal smooth structure with respect to $K$, which restricted to $Q-\{q\}$ is diffeomorphic to $Q-\{q\}$.
Moreover this normal smooth structure is correct with respect to $N$ if
\vspace{.05in}

\noindent \,\,\,\,(a) $dim\, N\leq 4$.

\noindent \,\,\,\,(b) $dim\, N\geq 5$ and the Whitehead group $Wh(N)$ of $N$ vanishes. }
\vspace{.1in}

\begin{center} {\bf \large  Section 8. Smoothing Hyperbolic Cones}\end{center} 
%\vspace{.1in}

Given an all-right spherical complex $P^m$ of dimension $m$ and a compatible smooth structure $\cS_P$ on $P$,
by Theorem 7.1 (see also remark 1 after the statement of 7.1) we can assume that $\cS_P$ is a normal smooth structure,
and $\cS_P$ has a normal atlas $\cA_P$.  The atlas $\cA_P$ and its induced differentiable structure $\cS_P$ are constructed (uniquely) from a set of link smoothings
$\cL=\{ h\0{\Delta}\}\0{\Delta\in P}$. To express this dependence we will sometimes write $\cA_P=\cA_P\big(\cL\big)$ and $\cS_P=\cS_P\big(\cL\big)$.\vspace{.1in}

Recall that the cone $\rC P$ has a piecewise hyperbolic metric 
induced by the piecewise spherical metric on $P$.
We denote these metrics by $\sigma\0{\rC P}$ and $\sigma\0{P}$ respectively.
As mentioned in Section 6, the piecewise hyperbolic metric $\sigma\0{\rC P}$ has a well defined ray structure.\vspace{.1in}

\noindent {\bf (8.0.1.)}\,\,\, Consider the following data.
\begin{enumerate}
\item[{\bf 1.}] A positive number $\xi$.
\item[{\bf 2.}] A collection $\s{d}=\{ d_2, d_3,....\}$ of real numbers, with $d_i>(6+2\xi)$.
We write $\s{d}(k)=\{ d_2, d_3,...,d_k\}$
\item[{\bf 3.}] 
A positive number  $r$, with $r>2d_i$, $i=2,...,m+1$, and $m$ as in item 5.
\item[{\bf 4.}] 
Real numbers $\varsigma\in (0,1)$, $c>1$, 
with $c\,\varsigma<e^{-4}$.
Hence we get sets of widths  (see 5.2 and 5.3)
$\sA=\sB(\varsigma)=\{\alpha_i\}$
and $\sB=\sB(\varsigma;c)=\{\beta_i\}$,
where $sin\,\alpha_i=\varsigma^{i+1}$, $sin\,\beta_i=c\,\varsigma^{i+1}$. 
\item[{\bf 5.}] An all-right spherical complex $P^m$, $dim\, P=m$, with compatible smooth normal atlas $\cA\0{P}\big(\cL_P\big)$, where $\cL_P$ is a smooth set
of link smoothings on $P$.
\item[{\bf 6.}] A diffeomorphism $\phi\0{P}=\phi\0{P,\cL\0{P}}: (P,\cS\0{P}(\cL_P))\ra \bS^{m}$ to the standard $m$-sphere. The map $\phi\0{P}$ is called a {\it global
smoothing for $P$, with respect to $\cS_P$ (or $\cA_p$, or $\cL_P$)}.
For $m=1$ we shall take $\phi\0{P}$ in a canonical way
(that is, depending only on $P$); see 8.1.
\end{enumerate}%\vspace{.1in}

The smooth atlas $\cA\0{P}(\cL_P)$ on $P$ induces, by coning, a smooth atlas
on $\rC P-\{o\0{\rC P}\}$, and, by item 6, this atlas together with  the coning
$\rC\phi\0{P}:\rC P\ra\R^{m+1}$ of the map
$\phi\0{P} $
induce a smooth atlas $\cA\0{\rC P}=\cA\0{\rC P}\big(\cL\0{P},\phi\0{P}\big)$ on $\rC P$. We denote the corresponding smooth structure by $\cS\0{\rC P}=\cS\0{\rC P}\big(\cL\0{P},\phi\0{P}\big)$. Note that
we get a diffeomorphism
$\rC\phi\0{P}:(\rC P,\cS\0{\rC P})\ra\R^{m+1}$.\vspace{.1in}

With the data given in items 1-6 above we will construct  by induction on the dimension $m$, the {\it smoothed} Riemannian metric
$\cG\big(P,\cL\0{P},\phi\0{P},r,\xi,\s{d},(c,\varsigma)\big)$ on the cone $\rC P$ of $P$, where we consider $\rC P$ with smooth structure $\cS\0{\rC P}$. \vspace{.1in}

In sections 8.1 and 8.2 we will assume $\xi$, $\s{d}$,
$c$, $\varsigma$ fixed. In particular we shall assume
$\s{A}$, $\s{B}$ fixed. 
So, to simplify our notation, we shall denote the smoothed metric by $\cG(P,\cL\0{P},\phi\0{P},r)$ or just $\cG(P,r)$ or
$\cG(P)$. In sections 8.3 and 8.4 we need to make explicit the dependence of the smoothed metric on  the other variables, and we will show that, given $\epsilon>0$, we can choose 
$r$ and $d_i$, $i=2,...m$, large so that $\cG\big(P,\cL\0{P},\phi\0{P},r,\xi,\s{d}, (c,\varsigma))$ has curvatures $\epsilon$-pinched to -1, provided the variables
satisfy certain conditions.
Before we begin with dimension 1 we need to discuss induced structures.\vspace{.1in}

Let $\Delta=\Delta^k\in P$. The {\it restriction of $\cL_P$ to $\sL(\Delta,P)$}
is the set $\cL_P|_{\sL(\Delta,P)}=\{h_{\Delta'}\}_{\Delta\subsetneqq\Delta'}$, see 7.2. Sometimes we will
just write $\cL_{\sL(\Delta,P)}$ or, more specifically,
$\cL_{\sL(\Delta,P)}(\cL_P)$. The corresponding induced atlas on $\sL(\Delta,P)$ is $\cA_{\sL(\Delta,P)}(\cL_P)=\{h^\bullet_{\Delta'}\}_{\Delta\subsetneqq\Delta'}$,
and sometimes we will
simply write $\cA_{\sL(\Delta,P)}$. The
smooth structure on $\sL(\Delta,P)$ induced by $\cA_{\sL(\Delta,P)}$ will be denoted
by $\cS_{\sL(\Delta,P)}(\cL_P)$, or simply by $\cS_{\sL(\Delta,P)}$. By Proposition 7.2.2
we have that, for $\Delta\in P$, the link smoothing $h_\Delta$ is a global
smoothing for $\sL(\Delta,P)$ with respect to $\cS_{\sL(\Delta,P)}$. Write $\phi_{\sL(\Delta,P)}=\phi_{\sL(\Delta,P)}(\cL_P)=h_\Delta$.
Therefore we obtain
the following {\it restriction rule}:\vspace{.1in}

\noindent {\bf (8.0.2.)}\hspace{1.3in} $\cL_P\longrightarrow \Big(\,\cL_{\sL(\Delta,P)}
(\cL_P)\,,
\,\phi_{\sL(\Delta,P)}(\cL_P)\, \Big)$ \vspace{.1in}

\noindent where $\cL_P$ satisfies 5 in (8.0.1) for $P$, and the objects	$\cL_{\sL(\Delta,P)}$, $\phi_{\sL(\Delta,P)}$ satisfy 5, 6 of (8.0.1) for $\sL(\Delta,P)$. 
The smooth structure on $\rC\sL(\Delta,P)$ constructed from
the data $\big(\cL_{\sL(\Delta,P)},\phi_{\sL(\Delta,P)}\big)$ will be denoted by
$\cS_{\rC\sL(\Delta,P)}(\cL_P)$, or $\cS_{\rC\sL(\Delta,P)}(\cL_{\sL(\Delta,P)},\phi_{\sL(\Delta,P)})$,
or simply by $\cS_{\rC\sL(\Delta,P)}$. 
The next lemma says that the restriction rule
(8.0.2) is transitive, that is, it respects the identity $\sL(\Delta^l,\sL(\Delta^j,P))=\sL(\Delta^k,P)$,
 where $\Delta^l=\sL(\Delta^j,\Delta^k)$ (see 5.3.2). \vspace{.1in}

\noindent {\bf Lemma 8.0.3.} {\it Let $\Delta^j\sbs\Delta^k\in P$ and let $\Delta^l=\sL(\Delta^j,\Delta^k)$. Then we have}
\vspace{.05in}

\hspace{1.5in}$
\cL_{\sL\big(\Delta^l,\sL(\Delta^j,P)\big)}\bigg(    \cL_{\sL(\Delta^j,P)} \big(  \cL_P \big)            \bigg)\,=\,\cL_{\sL(\Delta^k,P)}\big(  \cL_P \big)$\vspace{.05in}

\hspace{1.5in}$
\phi_{\sL\big(\Delta^l,\sL(\Delta^j,P)\big)}\bigg(    \cL_{\sL(\Delta^j,P)} \big(  \cL_P \big)            \bigg)\,=\,\phi_{\sL(\Delta^k,P)}\big(  \cL_P \big)
$\vspace{.05in}

\noindent {\bf Proof.} If we use the simplicial definition of link 
the identity $\sL(\Delta^l,\sL(\Delta^j,P))=\sL(\Delta^k,P)$ is an equality
of sets; hence the lemma follows from the definition
of $\cL$ and $\phi$. This proves the lemma.\vspace{.1in}

Recall that we have an identification $\rC\s{Star}(\Delta,P)=\rC\Delta\times\rC\sL(\Delta,P)$
(see 6.1.3). The ``open" version of this identification is
$\stackrel{\circ}{\rC\s{Star}}(\Delta,P)=\rC\dDelta\times\rC\sL(\Delta,P)$, where
$\stackrel{\circ}{\rC\s{Star}}(\Delta,P)=\rC(\stackrel{\circ}{\s{Star}}(\dDelta,P))$.
Here $\stackrel{\circ}{\s{Star}}(\dDelta,P)=\,\stackrel{\circ}{\s{N}}_{\pi/2}(\dDelta,P)$. Note that  $\stackrel{\circ}{\rC\s{Star}}(\Delta,P)$ as an open subset of
$\rC P$ has the induced smooth structure $\cS_{\rC P}|_{\stackrel{\circ}{\rC\s{Star}}(\Delta,P)}$, and, for simplicity, we will just write $\cS_{\rC P}$. On the other hand note
that $\rC\dDelta=\HH^{k+1}_+\sbs\HH^{k+1}$ has the natural smooth structure
$\cS_{\HH^{k+1}}$, and $\rC\sL(\dDelta,P)$ has the smooth structure
$\cS_{\rC\sL(\Delta,P)}$. Therefore we can give
$\rC\dDelta\times\rC\sL(\Delta,P)$ the ``product" smooth structure
$\cS_\times=\cS_{\rC\dDelta\times\rC\sL(\Delta,P)}$.\vspace{.1in}

\noindent {\bf Lemma 8.0.4.} {\it The following identification is a diffeomorphism}\vspace{.05in}

\hspace{1.5in}$
\big(\,\stackrel{\circ}{\rC\s{Star}}(\Delta,P)\,,\, \cS_{\rC P}\,   \big)\,\,=\,\, \big(  \,\rC\dDelta\times\rC\sL(\Delta,P)\,,\,   \cS_\times \, \big)
$\vspace{.1in}

\noindent {\bf Proof.} We use the variables $s$, $t$, $r$, $y$, $v$, $x$, $w$, $u$, $\beta$ defined
in Section 2. We also use the notation from 6.1.2, 6.1.3
and 7.0.1.
We assume that the image of the chart $h^\bullet_\Delta$
is $\stackrel{\circ}{\s{N}}_{\pi/2}(\Delta,P)$. By rescaling,
and using the notation in 6.1.2 and 6.1.3
we can write 
 \begin{equation*}\begin{array}{lccc}
h^\bullet_\Delta\,\,:&\D^{m-k}(\pi/2)\times\dDelta&\longrightarrow&P\\\\
&(\,\beta\,u'\,,\, w\,)&\mapsto&\big[\,w,  h_\Delta(u')  \,\big](\beta)
\end{array}
\tag{1}
\end{equation*}
\noindent where $\D^{m-k}(\pi/2)$ is the disc of radius $\pi/2$, and 
we are expressing and element $\D^{m-k}(\pi/2)$ as $\beta u'$, with
$\beta\in [0,\pi/2)$, $u'\in\bS^{m-k-1}$. A chart for
$(\,\stackrel{\circ}{\rC\s{Star}}(\Delta,P)\,,\, \cS_{\rC P}\,  )$
is the cone of $h^\bullet_\Delta$, which we shall denote by $h^\ast_\Delta$. Explicitly,
from (1) we have (see 7.0.1)
 \begin{equation*}\begin{array}{lccc}
h^\ast_\Delta\,\,:&\R_+\times\D^{m-k}(\pi/2)\times\dDelta&\longrightarrow&\rC P\\\\
&(\,s\,,\,\beta\,u'\,,\, w\,)&\mapsto&s\,\big[\,w,  h_\Delta(u')  \,\big](\beta)
\end{array}
\tag{2}
\end{equation*}
\noindent And for $(  \,\rC\dDelta\times\rC\sL(\Delta,P)\,,\,   \cS_\times \, )$  we can take the following chart
\begin{equation*}\begin{array}{lccc}
h^\dag_\Delta\,\,:&\R_+\times \R^{m-k}\times\dDelta&\longrightarrow&\rC\dDelta\times\rC\sL(\Delta,P)\\\\
&(\,t\,,\,r\,u'\,,\, w\,)&\mapsto&\big(\,t\,w\,,\,  r\,h_\Delta(u')  \,\big)
\end{array}
\tag{3}
\end{equation*}
\noindent where we  write an element in $\R^{m-k}$ as
$r u'$, $r\in[0,\infty)$, $u'\in\bS^{m-k-1}$. From (2) and (3)  and 6.1.3 we get
\begin{equation*}
\Big(h^\dag_\Delta\Big)^{-1}\circ h^\ast_\Delta\big(\,s\,,\, \beta\,u'   \,,\, w  \,\big)
\,\,=\,\, \Big( t\,,\,   r\,u'\,,\, w  \,\Big)
\tag{4}
\end{equation*}
\noindent where the relationship  between the variables
$s$, $\beta$, $t$, $r$ is the following (see Section 2). There is a right hyperbolic triangle
with catheti of length $t$, $r$, hypotenuse of length $s$ and angle $\beta$
opposite to the cathetus of length $r$. Using hyperbolic trigonometry we can
find an invertible transformation $(s,\beta)\ra (t,r)$. In particular 
$r=sinh^{-1}(sin\beta\, sinh(s))$.  The variables $s$ and $t$ are never zero,
but $\beta$ and $r$ could vanish. Note that $\beta=0$ if and only if $r=0$. 
To get differentiability at $\beta=0$ note that the map $(s,\beta u')\ra ru'$
can be rewritten as $(s,z)\ra(\frac{r(s,\beta)}{\beta}z)$, $\beta=|z|$, which is smooth
because $\frac{r(s,\beta)}{\beta}$ is a smooth even function on $\beta$. Similarly,
the smoothness of the inverse of the map in (4) follows from the fact that
the map $(t,r)\ra\frac{\beta(r,t)}{r}$ is a smooth even function on $r$.
This proves the lemma.
\vspace{.2in}

\noindent {\bf \large 8.1 Dimension One.} 
\vspace{.1in}

An all-right spherical complex $P^1$ of dimension one satisfying (8.0.1) (6) 
is formed by a finite number
$k'$ of segments of length $\pi/2$ glued successively forming a circle. Hence we shall canonically take (up to rotation)
$\phi=\phi\0{P}:P\ra\bS^1$ so that $\phi$ maps each 1-simplex to an arc of length $2\pi/k'$, and it does so with constant speed. Using $\phi$  we shall identify $P$ with a circle of length $k'\pi/2$.   Therefore $P$ with metric $\sigma\0{P}$, is isometric to $\bS^1$
with metric $k\sigma\0{\bS^1}$, $k=k'/4$. Consequently we identify $\rC P$ with $\R^2$,
and $\rC P-\{ o\0{\rC P}\}$ to $\R^2-\{0 \}$ with hyperbolic metric
$\sigma\0{\rC P}=sinh^2(t)\,k\,\sigma\0{\bS^1}+dt^2$. Notice that this metric is smooth on $\R^2$ away from the cone point $o\0{\rC P}=0\in\R^2$,
and it does have a singularity at 0 unless $k=1$.\vspace{.1in}

Recall we are assuming  $r> d_2$, where both $r$ and $d_2$ are given. Next we give two constructions of $\cG(P)$. The first one
is very explicit (see Gromov-Thurston \cite{GT}) and does not use 
(directly) any of the methods introduced
previously.  The second one looks more like
the inductive construction in 8.2, and uses the construction
given in Section 4. These two constructions are slightly different but both
satisfy the two properties {\bf P'1}, {\bf P'2} given below. Here is the first construction.\vspace{.1in}

 Let $\rho$ be as in Section 1.
Define \vspace{.05in}

\hspace{1.5in} $\mu(t)=\mu\0{d_2,r,k}(t)=k\,\rho({\mbox{\scriptsize $\frac{t}{d_2}-\frac{r-d_2}{d_2}$}})\, + \,\big(1-\rho({\mbox{\scriptsize $\frac{t}{d_2}-\frac{r-d_2}{d_2}$}})\big)$\vspace{.1in}
 
\noindent hence
$\mu(t)=1$, for $t\leq r-d_2$ and $\mu(t)=k$ for $t\geq r$. Define
\vspace{.1in}

\hspace{2.15in}$\cG(P,r)\,\,=\,\,sinh^2(t)\,\mu(t)\,\sigma\0{\bS^1}\,+\,dt^2
$\vspace{.1in}

Since the metric $\cG(P,r)$ is equal to the canonical hyperbolic warp metric
$sihn^2(t)\sigma\0{\bS^1}+dt^2$ on the ball of radius $r-d_2$, we can
extend $\cG(P,r)$ to the cone point $o\0{\rC P}=0\in\R^2$.
It is straightforward to verify that $\cG(P,r)$ satisfies the following three properties:
\begin{enumerate}
\item[{\bf P'1.}] The metrics $\cG(P,r)$ and $\sigma\0{\rC P}$ have the same ray structure.
\item[{\bf P'2.}] The metric $\cG(P,r)$ coincides with $\sigma\0{\rC P}$ outside
the ball of radius $r$.
\item[{\bf P'3.}] The metric $\cG(P,r)$ coincides with $sinh^2(t)\sigma\0{\bS^1}
+dt^2$ on
the ball of radius $r-d_2$.
\item[{\bf P'4.}] The family of metrics $\{\cG(P,r)\}\0{r>d_2}$ has cut limits (that is, it has cut limits on $I=\R$, see Section 4). 
Notice that $d_2$ is fixed and $r$ is the index of the family.\end{enumerate}

\noindent Actually from the definition of cut limit we have that the cut
limit of $\cG(p,r)$ at $b$ is \vspace{.1in}

\noindent {\bf (8.1.2.)}\hspace{.5in}$\Big(\,\lim_{r\ra\infty}\mu\0{d_2,r,k}(r+b)\,\Big)\,\sigma\0{\bS^1}\,=\,\Big(\, 1\,+\,\big(k-1\big)\rho\big( 1+\frac{b}{d_2}  \big)\,\Big)\,\sigma\0{\bS^1}$.\vspace{.1in}

Here is the second construction.  Recall that we are identifying $\rC P-\{o\0{\rC P}\}$ with
$\R^2-\{0\}$. Consider the constant $\odot$-family of metrics $\{g\0{r}\}_{r-\frac{1}{2}>d_2}$
given by $g\0{r}=k\sigma\0{\bS^1}+dt^2=\sigma\0{\rC P}$. Now just define 
\vspace{.1in}

\hspace{1.6in}$\cG(P,r)\,=\,   \cG(P,r,\s{d})  \,=\, \cH\0{r\,,\,d\0{2}}(g\0{r})$
\vspace{.1in} 

\noindent In this case we also have that $\cG(P)$ satisfies 
{\bf P'1}, {\bf P'2},  {\bf P'3} and {\bf P'4} (see Proposition 4.4).
\vspace{.2in}

\noindent {\bf \large 8.2 The Inductive Step .} 
\vspace{.1in}

Recall that in this section we are assuming $\xi$, $c$, $\varsigma$ (hence $\s{A}$, $\s{B}$),
and $\s{d}$ constant.
With the data $\xi$, $\s{A}$, $\s{B}$, $r>0$ and an all-right spherical complex $P$ we constructed in Section 6.2 the numbers $r\0{k}=r\0{k}(r)$ and the sets 
$\cY(P,\Delta^k,r)$, $\cY(P,r)$, $\cX(P,\Delta^k,r)$, $\cX(P,r)$, 
where $\Delta^k\in P$. The inverse of the function
$r\0{k}=r\0{k}(r)$ shall be denoted by $r=r(r\0{k})$.
Recall also that in  6.1.4 we identified $\rC\s{Star}(\Delta^k, P)$, with metric $\sigma\0{\rC P}|\0{\rC\s{Star}(\Delta^k, P)}$, with
$\cE_{\rC\Delta^k}(\rC\sL(\Delta^k))$, with metric $\cE_k(\sigma\0{\rC\sL(\Delta^k,P)})$.
We will use these objects in this section.\vspace{.1in}

Let $m\geq 2$ and
suppose that for every triple $(P,\cL\0{P},\phi\0{P})$,  $j=dim\,P\leq m-1$, as in items 
5 and 6 of (8.0.1), and $r>d_i$, $i=2,...,m+1$\,
there are couple of Riemannian metrics: the {\it smoothed} metric 
$\cG(P,\cL\0{P},\phi\0{P},r,\xi,\s{d},(c,\varsigma))$ (or simply $\cG(P,\cL\0{P},\phi\0{P},r)$, or even $\cG(P,r)$), 
and the {\it patched} metric $\ccP(P,\cL\0{P},r)$ (or just $\ccP(P,r)$), satisfying the 
following properties
\begin{enumerate}
\item[{\bf  P1.}] The smoothed metric $\cG(P,r)$ is a Riemannian metric on $(\rC P,\cS_{\rC P})$, and it has
the same ray structure as $\sigma\0{\rC P}$.
\item[{\bf  P2.}] The patch metric $\ccP(P,r)$ is a Riemannian metric defined outside the ball in $\rC P$ of radius
$r\0{j-2}-(4+2\xi)$ (with smooth structure $\cS_{\rC P}$), and it has
the same ray structure as $\sigma\0{\rC P}$.
%\item[{\bf  P3.}] The family of metrics $\Big\{\,\cG\big(P,r(r\0{j-2})\,%\big)\Big\}\0{r\0{j-2}}$ has cut limits.
\item[{\bf  P3.}] On $\cY(P,\Delta^k,r)$, $k\leq j-2= dim\,P-2$,
 the patched metric $\ccP(P,r)$ coincides with the
metric    
$$\cE_{\rC \Delta^k}\bigg(\cG\Big(\sL(\Delta^k,P),r\Big)\bigg)$$
where $\cG\Big(\sL(\Delta^k,P),r\Big)=\cG\Big(\sL(\Delta^k,P),\cL_{\sL(\Delta^k,P)}(\cL_P),
\phi_{\sL(\Delta^k,P)}(\cL_P), r\Big)$
is defined on  \newline $(\rC\s{Star}(\Delta,P), \cS_{\rC P})$.
(Recall $ \cY(P,\Delta^k,r)\sbs\rC\s{Star}(\Delta^k,P)=\rC\Delta^k\times\rC\sL(\Delta^k,P)$, see 6.1.3, 6.1.4, and 8.0.4.)

\item[{\bf  P4.}] On $\cY(P,r)$
 the patched metric $\ccP(P,r)$ coincides with 
$\sigma\0{\rC P}$ (which is hyperbolic on $\cY(P,r)$\,).
\item[{\bf  P5.}] The metrics $\cG(P,r)$ and $\ccP(P,r)$ coincide outside the ball  in $\rC P$ of radius $r\0{j-2}+\frac{1}{2}$.
\end{enumerate}%\vspace{.2in}

\noindent Note that the patched metric $\ccP(P,\cL_P,r)$ does not depend on
$\phi\0{P}$.
Properties {\bf P3}, {\bf P4}, {\bf P5} and the definition of the sets
$\cX(P,\Delta^k,r)$, $\cX(P,r)$ imply
\begin{enumerate}
\item[{\bf  P6.}] On $\cX(P,\Delta^k,r)$, $k\leq j-2= dim\,P-2$,
 the smoothed metric $\cG(P,r)$ coincides with the
metric 
$$\cE_{\rC \Delta^k}\bigg(\cG\Big(\sL(\Delta^k,P),r\Big)\bigg)$$
where $\cG\Big(\sL(\Delta^k,P),r\Big)=\cG\Big(\sL(\Delta^k,P),\cL_{\sL(\Delta^k,P)}(\cL_P),
\phi_{\sL(\Delta^k,P)}(\cL_P), r\Big)$
is defined on \newline
$(\rC\s{Star}(\Delta,P), \cS_{\rC P})$.
\item[{\bf  P7.}] On $\cX(P,r)$
the smoothed metric $\cG(P,r)$ coincides with the
metric  $\sigma\0{\rC P}$.
\end{enumerate}

Note that the metrics $\cG(P^1,r)$ constructed for spherical all-right 
1-complexes in 8.1, together with the choice $\ccP(P^1,r)=\sigma\0{\rC P^1}$ 
satisfy properties {\bf P1-P5}. Indeed {\bf P1'} implies {\bf P1},
{\bf P2'} implies {\bf P5} (recall $r_{-1}=r$, see 6.2) and 
{\bf P2}, {\bf P3}, {\bf P4} are trivially satisfied.\vspace{.1in}

Now, assume we are given the data: $P$,
$dim\, P=m$, $\cL\0{P}$, $\phi\0{P}$, $r$ as  items 5 and 6 in (8.0.1). Define the patched metric $\ccP(P,r)=\ccP(P,\cL\0{P},r)$
on $\rC P-\B\0{r\0{m-2}-(4+2\xi)}(\rC P)$ as in {\bf P3} and {\bf P4} above.
That is, we define $\ccP(P,r)$ by demanding that:
\begin{enumerate}
\item[{\bf  P''3.}] On $\cY(P,\Delta^k,r)$, $k\leq dim\,P-2$,
  $\ccP(P,r)$ coincides with the
metric   $\cE_{\rC \Delta^k}\bigg(\cG\Big(\sL(\Delta^k,P),r\Big)\bigg)$.
\item[{\bf  P''4.}] On $\cY(P,r)$, the patched metric $\ccP(P,r)$ coincides with the
metric  $\sigma\0{\rC P}$.
\end{enumerate}

\noindent {\bf Lemma 8.2.1.} {\it The patched metric $\ccP(P,r)$ defined by
properties {\bf P''3} and {\bf P''4} is well defined.}\vspace{.1in}

\noindent {\bf Proof.} The metric $\ccP(P,r)$ is defined on the ``patches''
$\cY(P,\Delta,r)$, $\Delta\in P$, and $\cY(P,r)$. We have to prove that these definitions coincide on the
intersections $\cY(P,\Delta^k,r)\cap\cY(P,\Delta^j,r)$,
$\cY(P,r)\cap\cY(P,\Delta^j,r)$. If $\Delta^j\cap\Delta^k=
\emptyset$ then (vi) of Proposition 6.2.1 implies $\cY(P,\Delta^j,r)\cap\cY(P,\Delta^k,r)=
\emptyset$. Also if $\Delta^j\not\sbs\Delta^k$ and $\Delta^k\not\sbs\Delta^j$ by (vii) of Proposition 6.2.1, we also get
$\cY(P,\Delta^j,r)\cap\cY(P,\Delta^k,r)=\emptyset$. Therefore we assume
$\Delta^j\sbs\Delta^k$, $j<k$.\vspace{.1in}

Recall that $\cY(P,\Delta^j,r)\sbs\rC\s{Star}(\Delta^j,r)$
and  $\cY(P,\Delta^k,r)\sbs\rC\s{Star}(\Delta^k,r)$ (see 6.2.1 (i)).
The metrics  
\begin{equation*}
\begin{array}{ccc}
h
%&=&\cE_{\rC \Delta^j}\bigg(\cG\Big(\sL(\Delta^j,P),r\Big)\bigg)
&=&\cE_{\rC \Delta^j}\bigg(\cG\Big(\sL(\Delta^j,P),\cL_{\sL(\Delta^j,P)}(\cL_P), \phi_{\sL(\Delta^j,P)}(\cL_P),r\Big)\bigg)
\end{array}
\tag{1}
\end{equation*}
\begin{equation*}
\begin{array}{ccc}g
&=&\cE_{\rC \Delta^k}\bigg(\cG\Big(\sL(\Delta^k,P),\cL_{\sL(\Delta^k,P)}(\cL_P), \phi_{\sL(\Delta^k,P)}(\cL_P),r\Big)\bigg)
\end{array}
\tag{2}
\end{equation*}
\noindent are defined on the whole of $\rC\s{Star}(\Delta^j,P)$ and $\rC\s{Star}(\Delta^k,P)$,
respectively. 
From  6.1.3 we have that $\rC\s{Star}(\Delta^j,P)=\rC\Delta^j\times
\rC\sL(\Delta^j,P)$. And from Lemma 6.2.3 we have that
$\cY(P,\Delta^k,r)\sbs\rC\Delta^j\times
\cX\Big(\sL\big(\Delta^j,P\Big), \Delta^l , r \Big)$, where $\Delta^l=\Delta^k\cap
\sL(\Delta^j,P)$ (alternatively $\Delta^l$ is opposite to $\Delta^j$ in $\Delta^k$, or $\Delta^l=\sL(\Delta^j,\Delta^k)$).
Hence it is enough to prove that the metrics $h$ and $g$ coincide on 
$\rC\Delta^j\times \cX\Big(\sL\big(\Delta^j,P\Big), \Delta^l , r \Big)$.
But (2) and 6.1.6 (see also 6.1.5, or 6.1.2, and Proposition 2.5) imply 
\begin{equation*}
g=\cE_{\rC \Delta^j}\bigg[
\cE_{\rC \Delta^l}\bigg(\cG\Big(\sL(\Delta^k,P),  \cL_{\sL(\Delta^k,P)}(\cL_P), \phi_{\sL(\Delta^k,P)}(\cL_P),       r\Big)\bigg)\bigg]
\tag{3}
\end{equation*}
\noindent Note that
the inductive hypothesis (that is, properties {\bf P3}, {\bf P5}, which imply {\bf P6}) applied to the data $\sL(\Delta^j,P)$ and $\Delta^l$
%$\cL_{\sL(\Delta^k,P)}(\cL_P)$ and  $\phi_{\sL(\Delta^k,P)}(\cL_P)$,
gives us that on the set $\cX\Big(\sL\big(\Delta^j,P\big), \Delta^l , r \Big)$
we have 
\begin{equation*}\cG\Big(\sL(\Delta^j,P),
\cL_{\sL(\Delta^j,P)}(\cL_P),\phi_{\sL(\Delta^j,P)}(\cL_P),
r\Big)\,=\,\cE_{\rC \Delta^l}(f)
\tag{4}
\end{equation*}
\noindent where
{\small\begin{equation*}f=
\cG\bigg(\sL(\Delta^l,\sL(\Delta^j,P)),
\cL_{\sL(\Delta^l,\sL(\Delta^j,P))}(\cL_{\sL(\Delta^j,P)}),
\phi_{\sL(\Delta^l,\sL(\Delta^j,P)}(\cL_{\sL(\Delta^j,P)}),
r\bigg)
\tag{5}
\end{equation*}}
\noindent From (5), (3), (2) and transitivity of the restriction rule (see 8.0.3) we get
\begin{equation*}f=
\cG\Big(\sL(\Delta^k,P),
\cL_{\sL(\Delta^k,P)}(\cL_{P}),
\phi_{\sL(\Delta^k,P)}(\cL_{P}),
r\Big)
\tag{6}
\end{equation*}
\noindent Putting together (1), (4) and (6) we obtain an equation with the same right-hand side as in (3) but with $h$ instead of $g$ on the left-hand side. This proves
that $g=h$ on $\cY(P,\Delta^j,r)\cap\cY(P,\Delta^k,r)$.\vspace{.1in}

The proof that the patched metric is well defined on $\cY(P,\Delta^k,r)\cap\cY(P,r)$
uses a similar argument and it follows from 6.2.4, the inductive hypothesis
applied to $\sL(\Delta^k,P)$ (that is properties {\bf P4}, {\bf P5} which imply {\bf P7})
and 6.1.4. This proves the lemma.\vspace{.1in}

Recall that $r\0{m-2}=r\0{m-2}(r)$. Let $r=r(r\0{m-2})$ be the inverse,
where we consider $r\0{m-2}$ as a large real variable.
For $P=P^m$ using $\rC\phi\0{P}$ we get an identification between
$\rC P$ and $\R^{m+1}$. Therefore we can consider the
family of metrics
$\Big\{\ccP\big(P,r(r\0{m-2})\big)\Big\}_{r\0{m-2}}$ as a $\odot$-family of
metrics on $\R^{m+1}$. 
%Because of the unpleasant constant 1/2 in
%the warp forcing process we rescale this family and
%now consider the $\odot$-family of metrics.%\vspace{.1in}
%\noindent {\bf (8.2.2.)}\hspace{2in} $\Big\{\ccP\big(P,r(r\0{m-2})\big)\Big\}_{r\0{m-2}-
%{\mbox{\tiny$\frac{1}{2}$}}}$ 
%\noindent (recall that the indexation
%of the family tells us where we take the spherical cuts).
We define\vspace{.1in}

\noindent {\bf (8.2.2.)}\hspace{1.3in}
$\cG(P,r)\,\,=\,\,\cH\0{r\0{m-2}
%-{\mbox{\tiny$\frac{1}{2}$}}
\,,\,d\0{m+1}}
%-\frac{1}{2}}
\Big(\ccP\big(P,r(r\0{m-2})\big)   \Big)$\vspace{.1in}

\noindent 
By construction the patch metric satisfies {\bf P3} and {\bf P4}. 
Property  {\bf P2} for the patch metric can be proved by induction on the dimension of
$P$ (using properties {\bf P1-P5}), together with Remark 2.3.
Property {\bf P5} for $\cG(P,r)$ holds by construction and by (ii) of 4.4. Property {\bf P1} follows from {\bf P2} and {\bf P5}.
This concludes the construction of the smoothed
Riemannian metric $\cG(P,r)=\cG(P,\cL_P,\phi_P,r,\xi,\s{d}, (c,\varsigma))$.\vspace{.1in}

Hence, by construction, we have the following properties.\vspace{.1in}

\noindent {\bf P8.} The smoothed metric $\cG(P^m,r)$ is hyperbolic on 
$\B_{r\0{m-2}-d\0{m+1}}(\rC P)$.\vspace{.1in}

Note that the patched metric $\ccP(P^m,r)$ does not depend 
$d_i$, $i>m$. Also the smoothed metric $\cG(P^m,r)$ does not depend 
$d_i$, $i>m+1$.\vspace{.2in}

\noindent {\bf \large 8.3. On the Dependence of $\cG(P,r)$ on the Variable $c$.}\vspace{.1in}

In this section we show that the smoothed metric $\cG(P,r)=\cG(P,\cL_P,\phi_P,\xi,
r,(c,\varsigma))$ does not depend on the variable  $c$, provided
$c\,\varsigma$ is small enough. In the next section we will show that, assuming $\s{d}$ and $r$
large, the metric $\cG(P,r)$ is $\epsilon$-hyperbolic. However the excess of the 
$\epsilon$-hyperbolic charts does depend on $c$.
In the next result assume $\varsigma$, $\xi$ and $\s{d}$ fixed. We shall write
$\cG(P,r,c)=\cG(P,\cL_P,\phi_P,r,\xi,(c,\varsigma))$
and similarly for the patch metric. \vspace{.1in}

\noindent {\bf Proposition 8.3.1.} {\it Let $c'>c>1$ 
%and $\xi'>\xi>0$ 
be such that  $c'\varsigma <e^{-4%+\xi'
}$.
%,  $c\varsigma <e^{-(4+\xi')}$. 
Then on $\rC P-\B_{r\0{m-2}-(4+2\xi)}(\rC P)$  we have
$\ccP(P,r,c')\,\,=\,\,\ccP(P,r,c)$.
Also, on $\rC P$
we have $\cG(P,r,c')\,\,=\,\,\cG(P,r,c)$}.\vspace{.1in}

\noindent {\bf Proof.} Write $\s{A}'=\s{B}(c',\varsigma)$. Denote by
$\cY'(P,\Delta,r)=\cY(P,\Delta,r,\xi,(c',\varsigma))$ the sets obtained by replacing $c$ in the definition
of $\cY(P,\Delta,r)=\cY(P,r,\xi,(c,\varsigma))$ (see 6.2) by $c'$. 
Similarly we obtain $\cY'(P,r)$. 
%Assume, say, that $c'\geq c$.
%We claim that we can also assume $\xi'\geq \xi$. To see this suppose the %proposition holds
%when $c'\geq c$ and $\xi'\geq \xi$. Hence if
%$\xi'<\xi$, then $\cG(P,r,\xi',c')=\cG(P,r,\xi',c)=\cG(P,r,\xi,c)$ (here we use the %hypothesis
%$c\varsigma <e^{-(4+\xi')}$). Therefore in what follows
%we assume $c'\geq c$ and $\xi'\geq \xi$. In particular 
Also, let $s'\0{m,k}$ be obtained from $s\0{m,k}$ by replacing $c$ by $c'$
(see 6.2). Then $s'\0{m,k}>s\0{m,k}$. Since $c'>c$
we have 
\begin{equation*}
\cY(P,\Delta,r)\sbs\cY'(P,\Delta,r) {\mbox{\,\,\,and\,\,\,}}\cY(P,r)\sbs \cY'(P,r)
\tag{1}
\end{equation*}
We will prove
the proposition by induction on the dimension $m$ of $P^m$.
It can be checked from Section 8.1 that the case $m=1$ does not depend
on the variable $c$.
Assume $\cG(P^k,r,c')\,\,=\,\,\cG(P^k,r,c)$, for every $P^k$, $k<m$.
Consider $P^m$. First we prove that the corresponding patched metrics
$\ccP(P,r,c')$ and $\ccP(P,r,c)$ coincide. But it follows from properties {\bf P3} 
and {\bf P4} applied to both metrics, the inductive hypothesis and (1) that 
$\ccP(P^m,r,c'))\,\,=\,\,\ccP(P^m,r,c)$ on $\cY(P,\Delta^k,r)$, for all $\Delta^k\in P$,
$k\leq m-2$, and on $\cY(P,r)$. Therefore, by  6.2.1 (iv), 
the metrics $\ccP(P^k,r,c'))$, $\ccP(P^k,r,c)$ coincide on 
$\rC P-\B_{r\0{m-2}-(4+2\xi)}(\rC P)$.
Finally note that the smoothed metrics $\cG(P,r,c)$, $\cG(P,r,c')$ are obtained
from the corresponding patched metrics by using the hyperbolic forcing process of
Section 4. But this process depends only on $\s{d}$ and $r\0{m-2}=sinh^{-1}(\frac{sinh(r)}{sin(\alpha\0{m-2})})$. The former
is fixed and the latter, since $sin(\alpha\0{m-2})=\varsigma^{m-1}$ (see 6.2),
is independent of $c$ and $c'$. This proves the proposition.\vspace{.1in}

In the next section we will need the following result.
We use the notation in the proof of the previous proposition.
Recall $s'\0{m,k}$ is obtained from $s\0{m,k}$ by replacing $c$ by $c'$.\vspace{.1in}

\noindent {\bf Lemma 8.3.2.} {\it If $c'>c$ we have
\,\,{\small $(\,s'\0{m,k}\,-\,s\0{m,k}\,)\,\,>\,\,ln\Big(\frac{c'}{c}\Big)\,-\,1$},
 provided $r>1$}.%\vspace{.1in}

\noindent {\bf Proof.} From the definition at the beginning of 6.2 we have $s\0{m,k}=sinh^{-1}(c\,\frac{sinh(r)}{\varsigma^{m-k-2}})$ and
$s'\0{m,k}=sinh^{-1}(c'\,\frac{sinh(r)}{\varsigma^{m-k-2}})$.
A simple calculation shows that the function $t\mapsto sinh^{-1}(c't)-sinh^{-1}(c t)$
is increasing. And another calculation shows that the value of this function
at $t=1$ has value at least $ln(c')-ln(c)-1$. This proves the lemma.
\vspace{.2in}

\noindent {\bf 8.4. On the $\epsilon$-Close to Hyperbolicity of $\cG(P,r)$ .} \vspace{.1in}

In this section we prove that the smoothed
metrics on $\rC P^m$ are $\epsilon$-close to hyperbolic, provided $d_{2},...,  d_{m+1}$ and $r$
are large enough. 
Recall that an element of $\rC P$ can be written as $sx$, $s\geq 0$, $x\in P$.\vspace{.1in}

\noindent {\bf Lemma 8.4.1.}
{\it The family of metrics \, 
{\small $\Big\{\ccP\big(P^m,r(r\0{m-2})\big)\Big\}\0{r\0{m-2}}$}
\, has cut limits on $[0,\infty)$. Also, the family of metrics \, 
{\small $\Big\{\cG\big(P^m,r(r\0{m-2})\big)\Big\}\0{r\0{m-2}}$}
\, has cut limits on $\R$.}%\vspace{.1in}

\noindent {\bf Proof.}
First note 4.4 (v), 4.4 (vi) and 8.2.2 imply that the first
statement in the lemma implies the second. 
We prove the first statement  by induction on the dimension $m$ of $P^m$. For $m=1$ the lemma follows from (v), (vi) and (vii) of Proposition 4.4 (alternatively we can use {\bf P4'} and 8.1.2 in Section 8.1)
\vspace{.1in}

\noindent {\bf Claim.} {\it Suppose the $\odot$-family of metrics
$\Big\{\cG\big(\sL(\Delta^k,P^m),r(r\0{m-k-3})\big)\Big\}\0{r\0{m-k-3}}$
has cut limits on $\R$. Then the $\odot$-family of metrics $
\Big\{\,\,\cE_{\rC\Delta^k}\Big(\,\cG\big(\sL(\Delta^k,P),r(r\0{m-2})\,\Big)\,\,\Big\}\0{r\0{m-2}}
$ also has cut limits on $\R$.}%\vspace{.1in}

\noindent {\bf Proof of claim.} 
By construction  (see 4.4 (i) or {\bf P8}), the family
$\big\{\cG\big(\sL(\Delta^k,P),r(r\0{m-k-3})\big)\big\}\0{r\0{m-k-3}}$
satisfies the hypothesis of Proposition 4.7.2:
the family is hyperbolic around the origin.
Since $r\0{m-k-3}=sinh^{-1}(sinh(r\0{m-2}) sin(\alpha\0{k}))$
the claim follows from Proposition 4.7.2. This proves the claim.\vspace{.1in}

We continue with the proof of Lemma 8.4.1.
Assume the lemma holds for $P^k$, $k<m$. Let $P=P^m$
%We now prove that property {\bf P3} holds for $\{\cG(P^m,r(r\0{m-2}))\}$.
and $\Delta^k\in P$. Suppose that the lemma does not hold for 
the family ${\cal{F}}=\big\{\ccP\big(P^m,r(r\0{m-2})\big)\big\}\0{r\0{m-2}}$.
We will show a contradiction. To simplify our notation write $s=r_{m-2}$
and $g_s=\cG\big(P^m,r(s)\big)$. We have $g_s=sinh^2(r)(g_s)_t+dt^2$,
where $t$ is the distance to $o\0{\rC P^m}$.
Since $\cF$ does not have cut limits on $[0,\infty)$ the is a bounded closed interval
$I\sbs [0,\infty)$ such that $\cF$ does not have cut limits on $I$. For $(x,b)\in P^m\times I$ write $g^\ast_s(x,b)=\widehat{(g_s)}_{s+b}(x)$. Note that
$sinh^2(s+b)g^\ast_s(x,b)+dt^2\,=\, g_s\big((s+b)x\big)$. Since we are assuming that $\cF$ does not have cut limits on $I$ we have that
the sequence $\{g^\ast_s\}$ defined on $P\times I$
does not converge in the $C^2$ topology.  Hence there is a derivative $\p^J$,
for some multi-index of order $\leq 2$, and sequences $s_n\ra\infty$, $x_n\ra x$,
$b_n\ra b$ such that $|\p^Jg^\ast_{s_n}(x_n,b_n)-\p^Jg^\ast_{s_{n+1}}(x_n,b_n)|\geq a$
for some fixed $a>0$, and $n$ even. By Proposition
6.3.3 we have that $R_{x,b}(s)=(s+b)x\in\cY(P,\Delta^k,r(s))$, for some $\Delta^k$,
$k\leq m-2$, and $s>s'$, for some $s'$; or 
$R_{x,b}(s)=(s+b)x\in\cY(P,r(s))$, $s>s'$, for some $s'$.
We assume the first case: $R_{x,b}(s)=(s+b)x\in\cY(P,\Delta^k,r(s))$, for some $\Delta^k$,
$k\leq m-2$. The other case is similar.  Since $s$ is large, we also get that 
$R_{x,b}(s)=(s+b)x\in\cX(P,\Delta^k,r(s))$, $s>s'$.
Moreover, also by 6.3.3, we can assume
$R_{s_n,b_n}(s)=(s+b_n)x_n\in\cX(P,\Delta^k,r(s))$, for $s>s'$. 
But by property {\bf P6}, on $\cX(P,\Delta^k,r(s))$ the metric $g_s$
is equal to $\cE_{\rC\Delta^k}\Big(\,\cG\big(\sL(\Delta^k,P),r(s)\,\Big)$.
Consequently the family of metrics $\big\{ \cE_{\rC\Delta^k}\Big(\,\cG\big(\sL(\Delta^k,P),r(s)\,\Big) \big\}_s$ does not have cut limits on $I$ either.
But the claim, together with the inductive hypothesis, imply that this family
does have cut limits, which leads to a contradiction. This proves the lemma.\vspace{.1in}

For a positive real number $\xi$ and a positive integer write $\xi\0{k}=\xi-k+\frac{1}{k}$. Note that $\xi\0{1}=\xi$.\vspace{.1in}

%We now assume $\s{d}=\{d_i\}$ variable\vspace{.1in}

\noindent {\bf Proposition 8.4.2.} {\it Let $\varsigma\in (0,1),\,,\xi>0$, $c>1$, and consider
$\big(P^m,\,\cL\0{P},\,\phi\0{P}\big)$. Assume
\vspace{.05in}

\noindent\,\,\,\,(i)  $c\,\varsigma<e^{-(8+\xi)}$

\noindent\,\,\,\,(ii)  $c\,\geq e^{4+\xi}$

\noindent\,\,\,\,(iii)  $m+1\leq \xi$.\vspace{.1in}

Let
$\epsilon>0$. Then we have that $\cG\big(P,\cL\0{P},\phi\0{P},\xi,r,\s{d}, (c,\varsigma))$ is  $(B_a,\epsilon)$-close to hyperbolic ($a=r\0{m-2}-d\0{m+1}$), with charts of excess $\xi\0{m}$, provided  $d_i$ and $r-d_i$, $i=2,...m+1$, are sufficiently large.}\vspace{.1in}

\noindent {\bf Remarks.} \\
\noindent {\bf 1.} By ``sufficiently large"  we mean that there are
$r_i(P,\epsilon)$ and $d_i(P,\epsilon)$, $i=2,...,m+1$, such that the proposition holds
whenever we choose $r-d_i\geq r_i(P,\epsilon)$ and $d_i\geq d_i(P,\epsilon)$. We will write $r_i(P)=r_i(P,\epsilon)$,
and $d_i(P)=d_i(P,\epsilon)$, if the context is clear.

\noindent {\bf 2.} The choices of $c$, $\xi$ and $\varsigma$ do not depend on
$\epsilon$.

\noindent {\bf 4.} 
If we want the smoothed metric on a cone $\rC P^m$
to be $(B_a,\epsilon)$-close to hyperbolic we can choose
$\xi=m+1$, $c=e^{4+\xi}$ and $\varsigma<e^{-(12+2\xi)}$. With these choices the method would not work for $P$ of dimension $>m$.

\noindent {\bf 5.} The condition $c\,\varsigma=e^{-(8+2\xi)}$ is stronger than the
condition $c\, \varsigma<e^{-4}$. The latter is used to construct the
smoothed metric but it is not strong enough to give us $\epsilon$-close to hyperbolicity.\vspace{.1in}

\noindent {\bf Proof.} 
We assume $c$, $\xi$, $\varsigma$ fixed and satisfying (i)
and (ii), that is, $c\,\varsigma<e^{-(8+\xi)}$ and $c\geq e^{4+\xi}$. We will only mention the relevant objects to our argument
in the notation for the smoothed metrics. That is, we will write
$\cG(P,\s{d},r,\xi, (c,\varsigma))$ or just $\cG(P,\s{d},r)$.
Our proof is by induction on the dimension $m$ of $P^m$, with $m+1<\xi$. Without loss of generality we can assume every $\epsilon$ we take satisfies:
{\footnotesize \begin{equation*}
\epsilon\,<\,\frac{1}{(1+\xi)^2}
\tag{1}
\end{equation*}}
For $m=1$ we have that the proposition follows from 8.1 and
Theorem 4.5 by writing $\lambda=r$, choosing $g_r=\sigma\0{\rC P}$,
replacing $\xi$ by $\xi+1$, and taking $\epsilon'=\epsilon$. 
Also, since $g_r=\sigma\0{\rC P}$ is $\epsilon$-close to hyperbolic, for every $\epsilon$,
we can take the $\epsilon$ in 4.5 to be zero. With all these choices Theorem 4.5 implies that
$\cG(P,d_2,r)$ is $\epsilon$-close to hyperbolic, with charts of excess $\xi=\xi\0{1}$, provided $r-d_2$ and $d_2$ are large enough.\vspace{.1in}

Let $m$ such that $m+1\leq\xi$. We write $a\0{k}=r\0{k-2}-d\0{k+1}$,
and note that $\cG(P^k,r,\s{d})$ is, by construction (see {\bf P8}),
radially hyperbolic on the ball of radius $a\0{k}$.
We now assume that the proposition holds for all $k<m$. That is, given $\epsilon>0$
and $P^k$,  the smoothed metric $\cG(P^k,r,\s{d})$ is $(B\0{a\0{k}},\epsilon)$-close to hyperbolic,
with charts of excess $\xi\0{k}$,
provided $r-d_i$ and $d_i$, $i=2,...,d_{k+1}$ are large enough. Note that, since we
are assuming $k<m$, we get that $k+1<\xi$. 
We use the following notation
{\footnotesize \begin{equation*}A\0{k}=C(m-k,k+1,\xi)\hspace{.5in}
B=C_2(\xi)\hspace{.5in}
\epsilon\0{k}=\frac{\epsilon}{3A\0{k}B}
\tag{2}
\end{equation*}}
\noindent where $C$ is as in Theorem 2.6 and $C_2$ as in Theorem 4.5. Let $P=P^m$. For $k<m$ write $L_k=\big\{ \sL(\Delta^k,P) \big\}\0{\Delta^k\in P}$.  A generic element in $L_k$ will be denoted by $Q=Q^j$, $j+k=m-1$. By inductive
hypothesis, for each $Q^j$ there are $r_i(Q^j)=r_i(Q^j,\epsilon\0{k})$ and $d_i(Q^j)=d_i(Q^j,\epsilon\0{k})$, $i=2,...,j+1$ such that
$\cG(Q,r,\s{d})$ is $(B\0{a\0{j}},\epsilon\0{k})$-close to hyperbolic, with charts of excess $\xi\0{j}$,
provided $r-d_i\geq r_i(Q^j)$ and $d_i\geq d_i(Q^j)$. For $i\leq m$, let
$d_i(P)$ be defined by\vspace{.05in}

\hspace{2.2in}{\footnotesize$d_i(P)\,=\, {\mbox{max}}\0{Q^j,\,i\leq j+1}\Big\{ d_i(Q^j)   \Big\}$}\vspace{.1in}

\noindent We write $\s{d}(P)=\{ d_2(P),...,d_{m}(P),...\}$ where
$d_i(P)$, $i>m+1$, is any positive number. This is just for notational purposes
and the arguments given below will not depend the $d_i(P)$, $i>m+1$.
We do reserve the right to later choose $d_{m+1}(P)$ larger.
Also write\vspace{.05in}

\hspace{1.6in}
{\footnotesize$ r_i(P)=d_i(P)\,+\, {\mbox{max}}\0{{Q^j,\, i\leq j+1}}\big\{\,4\,ln(m)\,,\, r_i(Q^j)\,  \big\}$}\vspace{.1in}
%\,\,\,\,\,\, and \,\,\,\,\,\, $r\0{P}={\mbox{max}}\0{0\leq i\leq m}r_i(P)$\\

\noindent Therefore %, from the inductive hypothesis and property {\bf P8} 
we get that\vspace{.1in}

\noindent {\bf (8.4.3.)} {\it 
For every $Q^j\in L_k$, the metric 
$\cG(Q^j,r,\s{d})$ is $(B\0{a\0{j}},\epsilon\0{k})$-close to hyperbolic, with charts of excess, 
\newline \hspace*{.61in}$\xi\0{j}$ provided $r-d_i\geq r_i(P)$ and $d_i\geq d_i(P)$, $i=2,...,k+1$.}\vspace{.1in}

\noindent  By definition we have $r\0{i}(P)\geq 4\,ln(m)$. Hence, if $r>r\0{i}(P)$ and $0\leq j\leq m-1$
we get that $\xi\0{j}-e^{-(r/2)}>\xi\0{j}-\frac{1}{m^2}>\xi-j+\frac{1}{j+1}\geq\xi-(m-1)+\frac{1}{m}$. This together with 
(8.4.3), Theorem 2.6 and the definitions given in
(2) imply  that\vspace{.1in}

\noindent {\bf (8.4.4.)} {\it 
For every \,$\sL(\Delta^k,P)\in L_k$, the metric 
$\cE_{k+1}\big(\cG(\sL(\Delta^k,P),r,\s{d})\big)$, defined  on
the space \newline \hspace*{.58in} $\cE_{k+1}\big(\rC\sL(\Delta^k,P)  \big)$, is 
$(B\0{a\0{j}},\frac{\epsilon}{3\,B})$-close to hyperbolic, with
charts of excess 
$\xi-(m-1)+\frac{1}{m}$, \newline \hspace*{.58in}  
provided $r-d_i\geq r_i(P)$ and $d_i\geq d_i(P)$, 
$i=2,...,k+1$.}
\vspace{.06in}

\noindent {\bf Lemma 8.4.5.} {\it The patched metric $\ccP(P,r,\s{d})$ is
radially $(\frac{\epsilon}{3\, B})$-close to hyperbolic 
on $\rC P-\B\0{r\0{m-2}-(1+\xi)}$, with
charts of excess 
$\xi-(m-1)+\frac{1}{m}$,
provided $r-d_i\geq r_i(P)$,
$d_i\geq d_i(P)$, $i=2,...,m$.}\vspace{.1in}

\noindent {\bf Proof.} 
Before we prove the lemma we need some preliminaries. Recall that we have functions $r\0{m,k}=r\0{m,k}(r)$
and $s\0{m,k}(r)$. For $\Delta=\Delta^k\in P$ write
$Y_\Delta=\cY(P,\Delta,r,\xi,(c,\varsigma))$ and 
$Y=\cY(P,r,\xi,(c,\varsigma))$ (see 6.2). For $\Delta=\Delta^k$, $k\leq m$ define

{\footnotesize \begin{equation*}
N_\Delta\,=\,\s{N}_{s\0{m,k}}(\rC\Delta,\rC P)\,-\, \bigcup_{\Delta^l\in P,\, l<k}\s{N}_{s\0{m,k}}
(\rC\Delta^l,\rC P)\,-\,\B_{r\0{m-2}-(1+\xi)}(\rC P)
\tag{3}
\end{equation*}}

{\footnotesize\begin{equation*}
N\,=\,\rC P\,-\, \bigcup_{\Delta^l\in P,\, l\leq m-2}\s{N}_{s\0{m,k}}
(\rC\Delta^l,\rC P)\,-\,\B_{r\0{m-2}-(1+\xi)}(\rC P)
\tag{4}
\end{equation*}}

\noindent 
Write $N_k=\cup_{\Delta^k\in P}N_{\Delta^k}$. 
It is straightforward to show that $\rC P-\B\0{r\0{m-2}-(1+\xi)}=
N\cup\bigcup_{k\leq m-2} N_k$.
Let
$c'=e^{4+\xi}c$. From  hypothesis (i), that is from $c\varsigma<e^{-(8+\xi)}$, we get that $c'\varsigma<e^{-4}$,
hence we can define the sets $Y'_\Delta=\cY(P,\Delta,r,\xi,(c',\varsigma))$
and $Y'=\cY(P,r,\xi,(c',\varsigma))$ (see 6.2). That is

{\footnotesize \begin{equation*}
Y'_\Delta\,=\,\stackrel{\circ}{\s{N}}_{s'\0{m,k}}(\rC\Delta^k,\rC P)\,-\, \bigcup_{\Delta^l\in P,\, l<k}\s{N}_{r\0{m,k}}
(\rC\Delta^l,\rC P)\,-\,\B_{r\0{m-2}-(4+2\xi)}(\rC P)
	\end{equation*}}

\noindent and analogously for $Y'$. Here $s'\0{m,k}$ is defined by replacing
$c$ by $c'$ in the definition of $s\0{m,k}$. Note that if we replace $c$ by 1
in the definition of $s\0{m,k}$ we obtain $r\0{m,k}$. This together with
hypothesis (ii),  the definition of $c'$, and Lemma 8.3.2 imply

{\small \begin{equation*}
\begin{array}{lll}
(\,s'\0{m,k}\,-\,s\0{m,k}\,)\,>\,\,3+\xi\\
(\,s\0{m,k}\,-\,r\0{m,k}\,)\,>\,\,3+\xi
\end{array}
\tag{5}
\end{equation*}}

\noindent It follows from $c'\varsigma<e^{-4}$, Lemma
8.3.1 and {\bf P3}, {\bf P5} that for each $\Delta^k\in P$, $k\leq m-2$, we have that\vspace{.1in}

\noindent {\bf (8.4.6.)} {\it  the metrics $\ccP(P,r,c')$,
$\cE_{k+1}\big(\cG(\sl(\Delta^k,P), r, c'\big)$ and 
$\cE_{k+1}\big(\cG(\sl(\Delta^k,P), r,c\big)$ coincide on $Y'_{\Delta^k}$}.\vspace{.1in}

For $p\in\rC P$ denote the ball of radius $s$ centered at $p$ by
$\B_{s,p}(\rC P)$, with respect to the metric $\sigma\0{\rC P}$.\vspace{.1in}

\noindent {\bf Claim 1.}
{\it For $\Delta=\Delta^k$, $k\leq m-2$, we have that}\hspace{.2in}
$d\0{\ccP}\Big(\,  N_\Delta    \,,\, \rC P\,-\, Y'_\Delta      \,\Big)\,\, \geq\, 
\, 3+\xi$\vspace{.1in}

\noindent {\bf Remark.} Here $d\0{\ccP}(.,.)$ denotes path distance with respect to the metric  $\ccP(P,r)$.\vspace{.1in}

\noindent {\bf Proof of claim.} Write $\Delta=\Delta^k$. From the definitions we have
$N_\Delta\sbs Y_\Delta\sbs Y'_\Delta\sbs\rC\s{Star}(\Delta,P)$. Note that
$\rC\s{Star}(\Delta,P)\sbs \rC P$ but we can also  consider
$\rC\s{Star}(\Delta,P)\sbs E=\cE_{k+1}(\rC\sL(\Delta,P))$.
We claim that we can work on $E$, that is, it is enough to prove that
$d\0{E}\Big(\,  N_\Delta    \,,\, \rC P\,-\, Y'_\Delta      \,\Big)>3+\xi$ ($E$ with metric $\cE_{k+1}(\cG(\sl(\Delta,P),r))$).
To see this note first that $\overline{N_\Delta}\sbs int\,Y'_\Delta$ (this follows from
$s\0{m,k}<s\0{m,j}$, $j<k$). Hence
if there is $p\in N_\Delta$ and $q\notin Y'_\Delta$
and $\alpha$ a path joining $p$ to $q$ of $\ccP$-length $< 3+\xi$ then
there a restriction $\beta$ of $\alpha$ such that: (1) it begins at $p$, \,(2) it
ends at some point $q'\in\p Y'_\Delta$,\, (3) it is totally contained in $\overline{Y'_\Delta}$,\,
(4) its $\ccP$-length is $< 3+\xi$. By property {\bf P3} the path
$\beta$ (now considered in $Y'\sbs E$, with metric $\cE_{k+1}(\cG(\sL(\Delta,P),r))$) has
the same properties (1)-(4). This shows that we can work on $E$ with metric
$\cE_{k+1}(\cG(\sL(\Delta,P),r))$ instead of $\rC P$ with metric $\ccP(P,r)$.\vspace{.1in}

Let $o\in \rC\s{Star}(\Delta,P)\sbs E$ correspond to $o\0{\rC P}$.
Let $p\in N_\Delta$, and $q\notin Y'_\Delta$.
From (3) we have 3 cases.\vspace{.05in}

\noindent {\bf Case 1.} $q\notin \s{N}_{s'\0{m,k}}(\rC\Delta,\rC P)$. Since $\s{N}_{s'\0{m,k}}(\rC\Delta,\rC P)=\bar{\HH}^{k+1}_+\times
\B_{s'\0{m,k}}(\rC\sL(\Delta,P))$ and $N_\Delta\sbs\s{N}_{s\0{m,k}}(\rC\Delta,\rC P)=
\bar{\HH}^{k+1}_+\times
\B_{s\0{m,k}}(\rC\sL(\Delta,P))$, this case follows from the first inequality in (5).
Note that the $s$ neighborhoods of $\HH^{k+1}$ in $\cE_{k+1}(\rC\sL(\Delta,P))$
with respect to the metrics $\cE_k(\cG(\sL(\Delta,P),r))$ and $\cE_{k+1}(\sigma\0{\rC\sL(\Delta,P)})$ coincide (see Remark 2.3).\vspace{.1in}

\noindent {\bf Case 2.} $q\in \s{N}_{r\0{m,j}}(\Delta^j,P)$, $j<k$.
Since $p\in N_\Delta$ w have that $p\notin \s{N}_{s\0{m,j}}(\Delta^j,P)$,
this case is similar to case 1, but uses the second inequality in (5), instead
of the first.\vspace{.1in}

\noindent {\bf Case 3.} {\it $q\in\B_{r\0{m-2}-(4+2\xi)}(E)$, where the ball is centered
at $o$.}
Since $p\notin \B_{r\0{m-2}-(1+\xi)}(E)$ this case follows from the fact
that $\big(r\0{m-2}-(1+\xi)\big)-\big(r\0{m-2}-(4+2\xi)\big)=3+\xi.$
This proves the claim.\vspace{.1in}

\noindent {\bf Claim 2.}  
%{\it Let $p\in N$. Then for $r$ large enough
%$\B_{(4+\xi),p}(\rC P)\sbs Y_\Delta$.}\\
{\it We have that}\hspace{.2in}$
d\0{\rC P}\big(\,  N    \,,\, \rC P\,-\, Y     \,\big)\,\, >\, 
\, 5+\xi$ \,(note that $Y=Y'$).\vspace{.1in}

\noindent {\bf Proof of claim.} First, we can extend the argument given in the
proof of claim 1 one more step, that is for $\Delta=\Delta^{m-1}$. The only problem
here is that we have not defined the sets $\cY(P,\Delta^k,r)$ for
$k=m-1$  (see 6.2). So just define $s'\0{m,m-1}=sinh^{-1}(c'sinh(r)\varsigma)$,
and $Y'_\Delta$ accordingly (as in 6.2). It can be checked, using the results in
Section 6, that the argument above goes through and we get
$d\0{\ccP}\Big(\,  N_\Delta    \,,\, \rC P\,-\, Y'_\Delta      \,\Big)\,\, \geq\, 
\, 3+\xi$ (in this case the patched metric is just $\sigma\0{\rC P}$). Since $Y'_\Delta\sbs Y$ it remains to prove that for $\Delta=\Delta^m$
we have $d\0{\ccP}\Big(\,  N_\Delta    \,,\, \rC P\,-\, Y     \,\Big)\,\, \geq\, 
\, 3+\xi$. To deal with this case define $Y^*=\rC P-\bigcup_{k< m}\s{N}_{r\0{m,k}}(\rC\Delta^k,\rC P)$ and let $Y^*_\Delta$ be the component of $Y^*$ contained in
$\rC\Delta$. Here we are taking 
$r\0{m,m-1}=sinh^{-1}(sinh(r)\varsigma)$.
Then $N_\Delta\sbs Y^*_\Delta\sbs Y$. We want to prove that
$d\0{\ccP}\Big(\,  N_\Delta    \,,\, \rC P\,-\, Y^*_\Delta     \,\Big)\,\, \geq\, 
\, 3+\xi$, but this case now
can be reduced to the case $\Delta=\Delta^m\sbs\bS^m\sbs\HH^{m+1}=\rC \bS^m$,
which can easily be dealt with. This proves the claim.\vspace{.1in}

We are now ready to prove the lemma.
By 8.3.1 it is enough to prove the lemma for
$\ccP'(P,r,\s{d})=\ccP(P,r,\s{d},c')$. Recall that
$\rC P-\B\0{r\0{m-2}-(1+\xi)}=
N\cup\bigcup_{k\leq m-2} N_k$. We prove by induction on $k$ that the patched metric
$\ccP'(P,r,\s{d})$ is radially $(\frac{\epsilon}{3\,B})$-close to hyperbolic, with excess
$\xi''=\xi-(m-1)+\frac{1}{m}$ on $N_k$. 
We begin with $k=0$. Assume $p\in N_0$. Then $p\in N_{\Delta^0}$ for some $\Delta^0$. From (8.4.4) there is a radially $(\frac{\epsilon}{3\,B})$-close to hyperbolic chart $\phi:\T_{\xi''}\ra \cE_k(\rC\sL(\Delta^0,P))$. But it follows from
Lemma 1.3,  claim 1, (1), and assumptions (i), (iii), that the image $\phi(\T_{\xi''})\sbs Y'_{\Delta^0}$.
This together with (8.4.6) imply $\phi$ is also  chart for $\ccP'(P,r,\s{d})$. This proves
the case $k=0$. The inductive step  $k\leq m-2$ is similar. It remains to prove
the case $p\in N$. But this case follows from a similar argument as 
above (in this case fitting a chart in $Y$) and using Claim 2 and Property {\bf P5}. 
This proves Lemma 8.4.5.\vspace{.1in}

We now finish the proof of Proposition 8.4.2. Set $\epsilon'=\frac{\epsilon}{3}$. By Lemma 8.4.5 we can
apply Theorem 4.5 to the family
$\big\{ \ccP(P,r(r\0{m-2}),\s{d}) \big\}_{r\0{m-2}}$.
Note that we have to use Lemma 8.4.1
to satisfy one of the hypothesis of Theorem 4.5.  Since $\epsilon'+B\frac{\epsilon}{3\,B}<\epsilon$
(recall $B=C_2$, see (2)) from Theorem 4.5 we
 obtain a number $r\0{m+1}(P)$ and a (possibly larger) number  $d_{m+1}(P)$ such that
$\cG(P,r,\s{d})$ is radially $\epsilon$-close to hyperbolic, provided $r-d_i\geq r_i(P)$ and
$d_i\geq d_i(P)$, $i=2,...,m+1$. Finally note that the excess of the charts given
by Theorem 4.5 is $(\xi-(m-1)+\frac{1}{m})-1=\xi\0{m}$.
This proves Proposition 8.4.2.
\vspace{.2in}

\noindent {\bf \large 8.5.  Smoothing Cones Over Manifolds.}
\vspace{.1in}  

As in the beginning of Section 8, let $P^m$ be an all-right spherical complex
and $\cS_P=\cS(\cL_P)$ a compatible normal smooth structure on $P$.
In the previous sections we have canonically constructed a 
Riemannian metric $\cG(P,\cL_P, \phi_P,r,\xi,\s{d},(c,\varsigma))$ on the cone
$\rC P$. An important assumption was that $(P,\cS_P)$ was diffeomorphic 
(by $\phi\0{P}$) to the sphere $\bS^m$. We cannot expect to do the same
construction on a general manifold $P$ because $\rC P$ is not in 
general a manifold. But we will canonically construct a complete  Riemannian 
metric on $\rC P-o\0{\rC P}$ that has some of the previous properties.\vspace{.1in}

We consider the same data as before: $P^m$, $\cL_P$, $r$, $\xi$, $\s{d}$, $(c,\varsigma)$ satisfying (8.0.1)
but with one change, replace the map $\phi\0{P}$ by a Riemannian metric
$h\0{P}$ on the closed smooth manifold $(P,\cS_P)$. Hence we begin with the 
following data:  $P^m$,  $\cL_P$, $h\0{P}$,
$r$, $\xi$,  $\s{d}$, $(c,\varsigma)$.\vspace{.1in}

First note that, by Theorem
7.1, a compatible normal smooth structure on $P$ exists. 
We will assume that $P$ has either dimension $\leq 4$ or $Wh(\pi_1P)=0$,
so that we can apply 7.3.1.
Note also that the sets $\cY(P,\Delta, r)$, $\cY(P,r)$ are defined for general $P$ (no just
for $P=\bS^m$) and satisfy all the properties given in Section 6.2. 
Now, since all the links of $P$ are spheres we can define, as in 8.1 and 8.2
the patch metric $\ccP(P,r)=\ccP(P,\cL_P,r, \xi,\s{d},(c,\varsigma))$ on $\rC P-\B_{r\0{m-2}-(4+2\xi)}(\rC P)$, and this metric satisfies properties {\bf P2}, {\bf P3}, {\bf P4} given in section
8.2.\vspace{.1in}

Recall that in 8.2 this construction is completed by applying hyperbolic forcing to the $\odot$-family of metrics $\Big\{\ccP\big(P,r(r\0{m-2})\big)\Big\}_{r\0{m-2}}$.
This method consists of two parts: warp forcing and then the
two variable deformation. In our more general setting here we can still apply warp forcing, but we cannot
directly apply hyperbolic forcing (at least not in the way given in Section 3)
because we do not have $P=\bS^m$. In our case, 
to finish our construction we apply first warp forcing and then a version of the two variable deformation for general $P$; this new version will use the metric 
$h\0{P}$ instead of the canonical metric $\sigma\0{\bS^m}$ on the sphere $\bS^m$.\vspace{.1in}

Consider now the  $\odot$-family of metrics $\Big\{\ccP\big(P,r(r\0{m-2})\big)\Big\}
_{r\0{m-2}}$ and apply warp forcing to obtain
\vspace{.1in}

\hspace{2in}{\footnotesize $
g\0{r\0{m-2}}\,\,=\,\,\cW_{r\0{m-2}}\Big(\,\cG\big(P,\,r(r\0{m-2})  \big)\,\Big)
$}\vspace{.1in}

\noindent and we have that $g\0{r\0{m-2}}$ is warped on 
$\B_{r\0{m-2}-\frac{1}{2}}(\rC P)-o\0{\rC P}$, specifically we have $g\0{r\0{m-2}}=sinh^2(t)\, g+dt^2$, where $g$
is a Riemannian metric on $P$ (it is the spherical cut of
$\ccP\big(P,r(r\0{m-2}))$ at $r\0{m-2}$) and $t$ is the distance-to-the-vertex
function on $\rC P$. Let $\rho\0{a,d}$ be the function in 3.1.
Also define the metric $g_t=h\,+\, \big(\rho\0{r\0{m-2}-d\0{m+1},d\0{m+1}}( t)\big)\,\big(g\,-\,h\big)$.
Now define the metric 
$\cG(P,h,r)=\cG(P,\cL_P,h,r,\xi,\s{d},(c,\varsigma))$ by\vspace{.1in}

\hspace{.4in} {\small$
\cG(P,h,r)\,\,=\,\, \left\{
\begin{array}{ll}
sinh^2(t)\,g_t +dt^2 & {\mbox{on}}\,\,\, 
\B_{r\0{m-2}}(\rC P)\,-\, \B_{r\0{m-2}-d\0{m+1}}(\rC P)\\ \\
\mu ^2(t)\, h\,+\, dt^2 & {\mbox{on}}\,\,\, 
 \B_{r\0{m-2}-d\0{m+1}}(\rC P)
\end{array}\right.
$}\vspace{.1in}

\noindent where $\mu(t)=\frac{e^{t}-e^{\lambda(t)}}{2}$, and
$\lambda=\rho\0{r\0{m-2}-2d\0{m+1},d\0{m+1}}$. 
Also we are assuming $r\0{m-2}-2d\0{m+1}>0$. Note that
$\cG(P,h,r)=\frac{1}{2}e^t h+dt^2$ on $\B_{r\0{m-2}-2d\0{m+1}}(\rC P)-o\0{\rC P}$,
that is for  $0< t\leq r\0{m-2}-2d\0{m+1}$.
We write $\rC P-o\0{\rC P}=P\times (0,\infty)$
and extend the metric $\cG(P,h,r)$ to $P\times \R$ by 
$\frac{1}{2}e^t h+dt^2$ for $-\infty< t\leq r\0{m-2}-2d\0{m+1}$.\vspace{.1in}

\noindent {\bf Corollary 8.5.1.} {\it The metrics $\cG(P,h,r)$ and $\ccP(P,r)$ have the following properties

\noindent\,\,\,\,(i) \,\,\,\,\,\,$\cG(P,r)$ is a Riemannian metric on $P\times\R$ that has
the same ray structure as $\sigma\0{\rC P}$ (on $P\times (0,\infty)$).

\noindent\,\,\,\,(ii)\,\,\, Properties  {\bf P3} and  {\bf P4}.

\noindent\,\,\,\,(iii) \,We have $\cG(P,h,r)=\frac{1}{2}e^t h+dt^2$ for $-\infty< t\leq r\0{m-2}-2d\0{m+1}$.

\noindent\,\,\,\,(iv)\,\, Given $\epsilon>0$ we have that the sectional curvatures of
$\cG(P,h,r)$ are $\epsilon$-pinched to -1 for $ t\geq$

\noindent\hspace{.47in}$ r\0{m-2}-2d\0{m+1}$ provided
$r-d_i$, $d_i$, $i=2,...,m+1$, and $r-2d\0{m+1}$ are large enough.}
\vspace{.1in}

\noindent {\bf Proof.} Item (i) follows %from Lemma 1.2 and 
the same argument used for {\bf P1}  in the spherical case.
Item (ii) is true by construction (see also Lemma 8.2.1). Item (iii)
follows from the discussion above, and (iv) from 8.4.2 and Bishop-O'Neill
warp curvature formula \cite{BisOn}, p.27. 
\vspace{.2in}

\begin{center} {\bf \large 9. On Charney-Davis Strict Hyperbolization Process.}
\end{center}

We use some of the notation in \cite{ChD}. In particular the canonical $n$-cube
$[0,1]^n$ will be denoted by $\square^n$ and $\dsquare^n=(0,1)^n$.
(This differs with the notation used in Section 7, where an $n$-cube was denoted by
$\sigma^n$.)
Also $B_n$ is the isometry group
of $\square^n$.\vspace{.1in}

A {\it Charney-Davis strict hyperbolization piece of dimension $n$} is a compact
connected orientable hyperbolic $n$-manifold with corners satisfying the properties stated in Lemma 5.1 of \cite{ChD}.
The group $B_n$ acts by isometries on $X^n$
and there is a smooth map $f:X^n\ra\square^n$ constructed in Section 5 of
\cite{ChD} with certain properties. 
We collect some facts from \cite{ChD}.

{\bf (1)} For any $k$-face $\square^k$ of $\square^n$ we have that $f^{-1}
(\square^k)$ is totally geodesic in $X^n$. Moreover
$X^n$  is a 

Charney-Davis hyperbolization
piece of dimension $k$. The  submanifold (with 
corners)
$f^{-1}(\square^k)$ is a 

$k$-$face$ of $X^n$. Note that the intersection of
faces is a face and every $k$-face is the intersection of 

exactly $n-k$ distinct
$(n-1)$-faces.

{\bf (2)} The map $f$ is $B_n$-equivariant.

{\bf (3)} The faces of $X^n$ intersect orthogonally.

{\bf (4)} The map $f$ is transversal to the $k$-faces of $\square^n$, $k<n$.

\noindent The $k$-face $f^{-1}(\square^k)$ of $X$ will be 
denoted by $X_{\square^k}$. The interior $f^{-1}(\dot{\square^k})$ will be denoted by $\dX_{\square^k}$. The following is proved
in \cite{O2}.\vspace{.1in}

\noindent {\bf Proposition 9.1.} {\it For every $n$ and $r>0$ the is a 
Charney-Davis hyperbolization piece of dimension $n$ such that the 
widths of the normal neighborhoods of every $k$-face, $k=0,...n-1$,
are larger that $r$.} \vspace{.1in}

For a $k$-face $X_{\square^k}$ and $p\in X_{\square^k}$, the set of inward normal vectors to
$X_{\square^k}$ at $p$ can be identified with the canonical all-right $(n-k-1)$-simplex
$\Delta_{\bS^{n-k-1}}$. In this sense we consider $\Delta_{\bS^{n-k-1}}\sbs T_pX$.
Similarly we can consider $\Delta_{\bS^{n-k-1}}\sbs T_q\square^n$, for 
$q\in \square^k$. We make the convention that
the two identifications above are done with respect to an ordering of
the $(n-1)$-faces $X_{\square^{n-1}}$ of $X$ and the corresponding
ordering for $\square^n$. For a proof of the following proposition see \cite{O2}.
\vspace{.1in}

\noindent {\bf Lemma 9.2.} {\it For $p\in \dX_{\square^k}$,
we have that $Df_p$ sends non-zero normal vectors to non-zero normal vectors; thus $Df_p|\0{\Delta_{\bS^{n-k-1}}}:\Delta_{\bS^{n-k-1}}\ra\Delta_{\bS^{n-k-1}}$. Moreover, {\bf n} $\circ \,(Df_p|\0{\Delta_{\bS^{n-k-1}}}):\Delta_{\bS^{n-k-1}}
\ra\Delta_{\bS^{n-k-1}}$
is the identity, where {\bf n}$(x)=\frac{x}{|x|}$ is the normalization map.}
\vspace{.1in}

The strict hyperbolization process of Charney and Davis is done by gluing copies
of $X^n$ using the same pattern as the one used to obtain the cube complex $K$
from its cubes. This space is called $K_X$ in \cite{ChD}. 
%We call this space the {\it piece-by-piece strict hyperbolization of $K$}.
Note that we get a map
$F:K_X\ra K$, which restricted to each copy of $X$ is just the map $f:X^n\ra\square^n$. We will write $X_{\square^k}=F^{-1}(\square^n)$,  and $\dX_{\square^k}=X_{\dsquare^k}=F^{-1}(\dsquare^k)$,
for a $k$-cube $\square^k$
of $K$.\vspace{.1in}

By Lemma 9.2 we can use the derivative of
the map $F:K_X\ra K$ (in a piecewise fashion) to identify $\sL(X_{\square^k},K_X)$ with $\sL(\square^k,K)$, where in both cases we consider the ``direction" definition of
link, that is, the
link $\sL(X_{\square^k},K_X)$ (at $p\in \dX_{\square^k}$) is the set of normal vectors to $X_{\square^k}$ (at $p$) and  the
link $\sL(\square^k,K)$ (at $q\in\dsquare^k$) is the set of normal vectors to $\square^k$ (at $q$). Hence we write
$\sL(X_{\square^k},K_X)=\sL(\square^k,K)$; thus the set of links for $K$ coincides
with the set of links for $K_X$.\vspace{.1in}

In what follows we assume that
the width of the normal neighborhoods of all $X_{\square}$ to be larger than $s\0{0}$, for some $s\0{0}$. Also let $r$ such that
$s\0{0}>2r$. By 9.1 we can take $s\0{0}$ and $r$ arbitrarily large.\vspace{.1in}

Let $X_{\square^k}\sbs X_{\square^{n}}$ be a $k$-face of $K_X$,
contained in the copy $X_{\square^n}$ of $X$ over $\square^n$. For a non-zero vector $u$ normal to $X_{\square^k}$ at $p\in X_{\square^k}$, and pointing inside 
$X_{\square^n}$, we have that $exp_p(tu)$ is defined and contained in
$X_{\square^n}$, for $0\leq t< s\0{0}/|u|$. Let $h_{\square^k}:\bS^{n-k-1}\ra\sL(\square^k,K)=\sL(X_{\square^k},K_X)$ be a link smoothing of the link corresponding to $\square^k\in K$.
We define the map \vspace{.1in}

\hspace{1.9in}{\small $H\0{\square^k}\,\,\,:\,\,\D^{n-k}\times \dX_{\square^k}
\,\,\,\,\longrightarrow\,\,\,\, K_X $}\vspace{.05in}

\hspace{2.8in}{\small $(tv,p)\,\,\,\longmapsto H\0{\square^k}(t\,v,p)=exp\0{p}\,\big(\,2r\,t\,\,h_{\square^k}(v)\,\big)$}\vspace{.1in}

\noindent where $v\in\bS^{n-k-1}$ and $t\in [0,1)$.
For $k=n$ we have that $H\0{\square^n}$ is the inclusion $\dX_{\square^n}\sbs K_X$
(or we can take this as a definition).
Note that $H\0{\square^k}$ is a topological embedding because we are assuming
the width of the normal neighborhood of $X_{\square}$ to be larger than $s\0{0}>2r$.
We call a chart of the form of $H\0{\square^k}$ (for some link smoothing
$h\0{\square^k}$) a {\it normal chart for the $k$-face $X_{\square^k}$}.
A collection $\big\{ H\0{\square} \big\}\0{\square\in K}$ of normal charts
is a {\it normal atlas}, and if this atlas is smooth (or $C^k$) the induced
differentiable structure is called a {\it normal smooth (or $C^k$) structure.} The following theorem is proved in \cite{O2}.\vspace{.1in}

\noindent {\bf Theorem 9.3.} {\it 
Let $\cL=\{h\0{\square}\}\0{\square\in K}$ be a set of link smoothigs for
$K$. If $\cL$ is smooth then the normal  atlas $\big\{ H\0{\square} \big\}\0{\square\in K}$ on $K_X$
is smooth.}\vspace{.1in}

We will write $\cA\0{K_X}=\big\{ H\0{\square} \big\}\0{\square\in K}$.
Note that the normal atlas $\cA\0{K_X}$ depends uniquely on the smooth set of
link smoothings $\cL=\{h\0{\square}\}\0{\square\in K}$ for $K$ (hence for $K_X$). To express this dependence we will sometimes write 
 $\cA\0{K_X}=\cA\0{K_X}(\cL)$.
We will denote by $\cS\0{K_X}=\cS\0{K_X}\big(\cL \big)$ the smooth structure on $K_X$ induced by 
the smooth atlas $\cA\0{K_X}$. The following Theorem is proved in \cite{O2}.\vspace{.1in}

\noindent {\bf Theorem 9.4.} {\it The smooth manifold $\big(K_X,\cS\0{K_X}\big)$ smoothly embeds in $(K,\cS')\times X$,
with trivial normal bundle. Here $\cS'$ is the normal smooth
structure on $K$ induced by $\cL$.}\vspace{.2in}

\noindent {\bf   9.5.  Hyperbolized Manifolds with Codimension Zero Singularities.}\vspace{.1in}

In this section we treat the case of manifolds with a one point singularity.
The case of manifolds with many (isolated) point singularities
is similar. We assume the setting and notation of Section 7.3.
Let $K_X$ be the Charney-Davis strict hyperbolization of $K$.
Denote also by $p$ the singularity of $K_X$.
Many of the definitions and results given before still hold (with minor changes)
in the case of manifolds with a one point singularity (see \cite{O2}
for more details).
\begin{enumerate}
\item[{\bf (1)}] Given a set of link smoothings for $K$ (hence for $K_X$) 
we also get a set of charts $H_\square$.
For the vertex $p$ we mean the cone map 
$H_p=\rC h_p:\rC N\ra \rC L\sbs K_X$. We will also denote the restriction
of $H_p$ to $\rC N-\{o\0{\rC N}\}$ by the same notation $H_p$.
As in item (2) of 7.3 here we are identifying $\rC N-\{o\0{\rC N}\}$ with $N\times (0,1]$ with the product smooth structure obtained from
\s{some} smooth structure ${\tilde{\cS}}_N$ on $N$.
As before $\{H_\square\}_{\square\in K}$ is a {\it normal atlas} for $K_X$
(or $K_X-\{p\}$). A normal atlas for $K-\{p\}$ induces a {\it normal smooth
structure on $K_X-\{p\}$.}
\item[{\bf (2)}]  Again we say that
the smooth atlas  $\{H_\square\}$ (or the induced smooth structure, or the set $\{h_\sigma\}$) 
is {\it correct with respect to $N$} if $\cS_N$ is diffeomorphic to ${\tilde{\cS}}_N$.
\item[{\bf (3)}] Let the set $\cL=\{h_\square\}_{\square\in K}$ 
induce a smooth structure on $K-\{p\}$, that is, $\cL$ is smooth. As in Theorem 9.3 we get that $\{H_\square\}_{\square\in K}$
is a smooth atlas on $K_X-\{p\}$ and it induces a normal smooth
structure $\cS_{K_X}$ on $K_X-\{p\}$. Moreover, from Theorem 7.3.1 we get that $\cS_{K_X}$ is correct with respect to $\cS_N$ when $dim\, N\leq 4$ (always) or 
when $dim\, N>4$, provided $Wh(N)=0$. Note that in this case we can take
the domain $\rC N-\{o\0{\rC N}\}=N\times (0,1]$ of $H_p$ with smooth
product structure $\cS_N\times \cS_{(0,1]}$.
\item[{\bf (4)}] It can be verified that a version of Theorem 9.4 also
holds in this case: $(K_X-\{p\},\cS_{K_X})$ smoothly embeds in  $(K-\{p\},\cS')\times X$
with trivial normal bundle.
\end{enumerate}\vspace{.1in}

\begin{center} {\bf \large Section 10. Proof of the Main Theorem.}
\end{center}

In Section 2  the concept of hyperbolic extension 
over hyperbolic space was introduced.
We  next extend, in the obvious way, this concept to hyperbolic extensions
over hyperbolic manifolds.\vspace{.1in}

As in Section 2, let $(N,h)$ be a complete Riemannian manifold with center
$o=o\0{N}$. Let $(P,\sigma\0{P})$ be a hyperbolic manifold.
The {\it hyperbolic extension of $h$  over $P$} is the Riemannian
metric $g=(cosh^2r)\sigma\0{P}+h$ on $P\times N$, where $r:N\ra[0,\infty)$ is the 
distance to $o$ function on $N$. We write $g=\cE_{P}(h)$
and $(P\times N, g)=\cE_{P}(N,h)$ (or simply $\cE_P(N)$) and we call
$\cE_P(N)$ the {\it hyperbolic extension of $N$ over $P$}.\vspace{.1in}
 
We now begin the proof of the Main Theorem. Let $M^n$ be a closed smooth
manifold. 
Let $K$ be a smooth cubulation of $M$ %(see appendix G) 
and  $K_X$ the Charney-Davis strict hyperbolization of $M$, as in Section 9.
We can assume that the Charney-Davis hyperbolization piece $X$ has normal
bundles with large widths (see 9.1), all larger than a large number $2s\0{0}>0$.
Let $\cA\0{K\0{X}}=\big\{  H_{\square} \big\}\0{\square\in K}$ be a
smooth normal atlas for $K_X$, and $\cS\0{K\0{X}}$ the induced
normal smooth structure on $K_X$. Recall that the
$H_\square$ are constructed from a smooth set
of link smoothings $\cL\0{K}=\{h\0{\square}\}\0{\square\in K}$ for the links of $K$ (or $K_X$).\vspace{.1in}

\noindent {\bf Remark.} The domains of the charts $H_{\square^k}$ are 
the sets $\D^{n-k}\times\dX_{\square^k}$. But in this section, for notational purposes,
we will consider the rescaling of $H_{\square^k}$ given by 
by $H_{\square^k}(tv,p)=exp_p(t\,h_{\square^k}(v))$,
defined on $\D^{n-k}(s\0{0})\times\dX_{\square^k}$. We shall
denote this chart also by $H_{\square^k}$. \vspace{.1in}

In what follows, to simplify our notation, we write $\sL(X_\square)=\sL(X_\square,K_X)$.
Recall that given $\square\in K$,  the set $\cL_K$ of link smoothings for the links $\sL(X_\square)$
of $K_X$ (and of $K$) induce, by restriction (see 7.2), the set of link smoothings 
$\{h\0{\square'}\in\cL_K\, ,\, \square'\subsetneq\square\}$ for the links of
$\sL(X_\square)$. We denote this induced set
of smoothings by $\cL_{\sL(X_\square)}$ or just $\cL\0{\square}$.\vspace{.1in}

The space $K_X$ has a natural piecewise hyperbolic metric which we denote by
$\sigma\0{K_X}$. The piecewise hyperbolic metric on the cones
$\rC\sL(X_{\square})$ of the all-right spherical simplices $\sL(X_{\square})$
will be denoted by $\sigma\0{\rC\sL(X_{\square})}$.
The restriction of $\sigma\0{K_X}$ to the totally geodesic
space $X_\square$ shall be denoted by $\sigma\0{X_\square}$.\vspace{.1in}

For $\square^k\in K$, the {\it (closed) normal neighborhood of $\dX_{\square^k}$ in $K_X$ of width $s<s\0{0}$} is the set $\s{N}\0{s}(\dX_{\square^k},K_X)=H_{\square^k}
\big(\D^{n-k}(s)\times \dX_{\square^k}\big)$. That is, it is the union of the images of
all geodesics of length $\leq s$ in each copy of $X$ containing $X_{\square^k}$,
that begin in (and are normal to) $\dX_{\square^k}$.
Similarly the {\it open
normal neighborhood of $\dX_{\square^k}$ 
of width $s<s\0{0}$}
is the set $\stackrel{\circ}{\s{N}}\0{s}(\dX_{\square^k},K_X)=H_{\square^k}
\big(\,int\,\D^{n-k}(s)\times \dX_{\square^k}\big)$.
Sometimes we will just write  $\s{N}\0{s}(\dX_{\square^k})=\s{N}\0{s}(\dX_{\square^k},K_X)$ and 
$\stackrel{\circ}{\s{N}}\0{s}(\dX_{\square^k})=\,
\stackrel{\circ}{\s{N}}\0{s}(\dX_{\square^k},K_X)$. 
Note that normal neighborhoods respect faces, that is:
\vspace{.1in}

\noindent {\bf (10.1)}\hspace{1.5in}
$\s{N}\0{s}(\dX_{\square^k},K_X)\cap\dX_{\square^l}=
\s{N}\0{s}(\dX_{\square^k\cap\square^l},\dX_{\square^l})$
\vspace{.1in}

Here $\s{N}\0{s}(\dX_{\square^k\cap\square^l},\dX_{\square^l})$
is the union of 
all geodesics of length $\leq s$ in the hyperbolization piece
$\dX_{\square^l}$ that begin in (and are normal to) $\dX_{\square^k
\cap\square^l}$.\vspace{.1in}

Since the
normal bundles of the $X_\square$ are canonically trivial (see construction of $X$  in \cite{ChD}, or Section 2 in \cite{O2}) we can canonically identify the neighborhood $\s{N}\0{s}(\dX_{\square})$
with $\dX_{\square^k}\times\rC_s\sL(X_{\square^k})$,
where  $\rC_s\sL(X_{\square^k})=\B_s\big( \rC\sL(X_{\square^k})  \big)$
is the closed $s$-cone  of length $s$, that is, it is the ball of radius $s$ on
the (infinite) cone $\rC\sL(X_{\square^k})$ centered at the vertex, see 6.1.
Similarly we have the identification $\stackrel{\circ}{\s{N}}\0{s}(\dX_{\square})= \dX_{\square^k}\times\stackrel{\circ}{\rC}_s\sL(X_{\square^k})$,
where  $\stackrel{\circ}{\rC}_s\sL(X_{\square^k})$
is the open $s$-cone  of length $s$.
Moreover these
identifications are also metric identifications, where we consider 
$\s{N}\0{s}(\dX_{\square},K_X)\sbs K_X$ with the (restricted) piecewise hyperbolic metric 
$\sigma\0{K_X}$ and $ \dX_{\square^k}\times\rC_s\sL(X_{\square^k})$
with the hyperbolic extension metric \vspace{.05in}

\hspace{1.7in}
$\cE_{\dX_{\square}}(\sigma\0{\rC\sL(X_\square)})=cosh^2(t)\,\sigma\0{\dX_{\square}}+
\sigma\0{\rC\sL(X_\square)}$\vspace{.05in}

\noindent where $t$ is the distance-to-the-vertex function on the cone
$\rC\sL(X_\square)$.\vspace{.1in}

\noindent {\bf Remarks.} 

\noindent{\bf 1.} The metric $\sigma\0{\rC\sL(X_\square)}$  is not smooth
but the formula above makes sense, giving a well defined piecewise hyperbolic
metric.

\noindent {\bf 2.} Since we are identifying $\s{N}\0{s}(\dX_{\square^k})$
with $\dX_{\square^k}\times\rC_s\sL(X_{\square^k})$ we will consider
 $\s{N}\0{s}(\dX_{\square})$ also as a subset of
$\dX_{\square^k}\times\rC\sL(X_{\square^k})$, where
$\rC\sL(X_{\square^k})$ is  the (infinite) cone over $\sL(X_{\square^k})$.
Note that the metric $\cE_{\dX_{\square}}(\sigma\0{\rC\sL(X_\square)})$
is defined on the whole of $\dX_{\square^k}\times\rC\sL(X_{\square^k})$.\vspace{.1in}

\noindent {\bf Lemma 10.2.} {\it Let 
 $\square^k=\square^i\cap\square^j$,  $k\geq 0$.
Let $s\0{1}, s\0{2}, s < s\0{0}$ be positive real numbers such that
$\frac{sinh\, s\0{1}}{sinh\, s},\frac{sinh\, s\0{2}}{sinh\, s} \leq\frac{\sqrt{2}}{2}$. Then}\,\,\, $
\sN\0{s\0{1}}\big(\dX_{\square^i}  \big)\,\, \cap\,\,\sN\0{s\0{2}}\big(\dX_{\square^j}  \big)
\,\,\, \sbs\,\,\, \sN\0{s}\big(\dX_{\square^k}  \big)
$.\vspace{.1in}

\noindent {\bf Proof.} Using 10.1 we can reduced 
the lemma to the case were $K_X$ is just a
hyperbolization piece $X$. This case is proved in \cite{O2}, 2.3. This proves the lemma.\vspace{.1in}

Suppose $\square^j\sbs\square^k\in K$. Then $\square^k$ determines
the all-right spherical simplex $\Delta\0{\sL(\square^j,K)}(\square^k)
=\square^k\cap\sL(\square^j,K)$ in $\sL(\square^j,K)=\sL(X_{\square^j})$.
We will just write $\Delta(\square^k)$ if there is no ambiguity.
(Other definition previously used: $\Delta\0{\sL(\square^j,K)}(\square^k)=
\sL(\square^j,\square^k)$.)
\vspace{.1in}

\noindent {\bf Lemma 10.3.} 
{\it Let $\square^j\sbs \square^k$ and $s\0{1},\, s\0{2}<s\0{0}$.
Then}\vspace{.05in}

\hspace{1in}
{\small $
\sN\0{s\0{1}}(\dX_{\square^j})\,\,\cap\,\,\sN\0{s\0{2}}(\dX_{\square^k})\,\,=\,\,
\dX_{\square^j}\,\,\times\,\,\sN\0{s\0{2}}\Big(\,  \rC\Delta(\square^k), \,\rC\0{s\0{1}} \sL\big(  X_{\square^j} \big)   \,\Big)$}\vspace{.05in}

\noindent {\bf Proof.} Let $p\in \sN\0{s\0{1}}(\dX_{\square^j})$. 
By the identification $\sN\0{s}(\dX_\square)=\dX_\square\times
\rC_s\big(\sL(X_{\square})\big)$ we can write 
$p=(u,x)\in \dX_\square^j\times\rC_{s\0{1}}\big(\sL(X_{\square^j})\big)$.
Let $\square^l$ such that $p\in \dX_\square^l$, and write
$\square^{j'}=\square^j\cap \square^l$,
$\square^{k'}=\square^k\cap \square^l$. 
Since $\sN\0{s\0{1}}(\dX_{\square^j})\cap \dX_{\square^{l}}
=\sN\0{s\0{1}}(\dX_{\square^{j'}},\dX_{\square^{l}})$ (see 9.),
there is a geodesic
segment $[x',p]$ be from $x'\in X_{\square^{j'}}$ perpendicular to
$X_{\square^{j'}}$. Therefore $x=x'$. 
%these segments are unique%

$X_\square^k=X_\square^k\cap X_\square^l$,

Since $\{x\}\times\rC_s\big(\sL(X_{\square^j})\big)$ is convex in
$\dX_\square^j\times\rC_s\big(\sL(X_{\square^j})\big)$,
we have \vspace{.05in}

\hspace{1in}
{\small $
d\0{\rC\0{s\0{1}}(\dX_\square^j\times\sL(X_{\square^j}))}\Big(\,p\,,\,
X_{\square^k} \, \Big)\,\,=\,\,
d\0{\rC\0{s\0{1}}(\{x\}\times\sL(X_{\square^j}))}\Big(\,u\,,\,
\rC\Delta(\square^j)\,  \Big)
$}\vspace{.05in}

\noindent where, as usual, $d\0{S}$ denotes distance on a space $S$.
Consequently the first term in the equality above is $<s\0{2}$ if and only
if the second term  is $<s\0{2}$. This proves the lemma.\vspace{.1in}

\noindent {\bf Remark.} Clearly the open  version of 
Lemma 10.3 also holds:\vspace{.05in}

\hspace{1in}{\small $
\stackrel{\circ}{\sN}\0{s\0{1}}(\dX_{\square^j})\,\,\cap\,\,\stackrel{\circ}{\sN}\0{s\0{2}}(\dX_{\square^k})\,\,=\,\,
\dX_{\square^j}\,\,\times\,\,
\stackrel{\circ}{\sN}\0{s\0{2}}\Big(\,  \rC\Delta(\square^j), \,
\stackrel{\circ}{\rC}\0{s\0{1}} \sL\big(  X_{\square^k} \big)   \,\Big)$}
\vspace{.1in}

Now, let $\s{d}$, $r$, $\xi$, $c$ and $\varsigma$ be as in items 1, 2, 3, 4 at the beginning of Section 8, and let the numbers $s\0{m,k}=s\0{m,k}(r)$, $r\0{m,k}=r\0{m,k}(r)$
be as in Section 6.2. For each  $\square^k\in K$ define the sets

{\small $$\begin{array}{rcl}
\cZ(X_{\square^k})&=& \stackrel{\circ}{\s{N}}_{s\0{n,k}}\big( X_{\square^k} \big)\,\, -\,\, \bigcup\0{i<k}\sN_{r\0{n,i}}
\big(  X_{\square^i}\big)\\\\
\cZ&=& K_X -\,\, \bigcup\0{i<n-1}\sN_{r\0{n,i}}
\big(  X_{\square^i}\big)
\end{array}$$}

\noindent By 9.1 we can take $s\0{0}$ as large as needed, hence we can assume that $\cZ(X_{\square^k})\sbs\,\,\stackrel{\circ}{\s{N}}_{s\0{0}}(X_{\square^k})$.\vspace{.1in} 

We next use the sets $\cX(P,\Delta,r)$ and $\cX(P,r)$ of Section 6.2.
The sets $\cX(\sL(X_{\square^k}),\Delta(\square^j),
r)$ and  $\cX(\sL(X_{\square^k}),r)$
are a subsets of the (infinite) cone $\rC\sL(X_{\square^k})$.\vspace{.1in}
%It restriction the the $s$-cone ($s<s\0{0}$) will be denoted by 
%$\cY_s(\sL(X_{\square^k}),\Delta(\square^j),
%r\0{j,k})$. Similarly $\cX_s(\sL(X_{\square^k}),\Delta(\square^j),
%r\0{j,k})$ is $\cX(\sL(X_{\square^k}),\Delta(\square^j),
%r\0{j,k})\cap\cY_s(\sL(X_{\square^k}),\Delta(\square^j),r\0{j,k})$

\noindent {\bf Lemma 10.4.} {\it We have the following properties
\vspace{.05in}

\noindent\,\,\,\,\,\,(i) If $\square^i\cap\square^j=\emptyset$ then
$\cZ(X_{\square^i})\cap\cZ(X_{\square^j})=\emptyset$.
\vspace{.05in}

\noindent\,\,\,\,\,(ii) If $\square^k=\square^i\cap\square^j$, $0\leq k<i,j$,  then
$\sN\0{s\0{n,i}}(X_{\square^i})\cap\sN\0{s\0{n,j}}(X_{\square^j})
\sbs \sN\0{r\0{n,k}}(X_{\square^k})$.
%\vspace{.05in}

\noindent\,\,\,\,(iii) If $\square^k=\square^i\cap\square^j$, $0\leq k<i,j$,  then
$\cZ(X_{\square^i})\cap\cZ(X_{\square^j})=\emptyset$.
\vspace{.05in}

\noindent\,\,\,\,(iv) If $\square^j\sbs\square^k$ then we have (see Remark 2 before 10.2)\vspace{.05in}

\hspace{1.4in}{\small $\cZ(X_{\square^j})\cap\cZ(X_{\square^k})\,\,\sbs\,\,\dX_{\square^j}\,\,\times\,\,
\cX\Big(\rC\sL(X_{\square^j}) ,\, \Delta(\square^k),\,r\Big) $}\vspace{.05in}

\noindent\,\,\,\,\,(v) For $k<n-1$ we have\,\, $\cZ\cap\cZ(X_{\square^k})\,\,\sbs\,\,\dX_{\square^k} \,\,\times\,\,
\cX\Big(\rC\sL(X_{\square^k}) ,\,r\Big)$.}\vspace{.05in}

\noindent {\bf Proof.} Let  $\square^i\cap\square^j=\emptyset$.
Then the distance in $K_X$ from $X_{\square^i}$ to $X_{\square^j}$ is 
at least $2s\0{0}$. This proves (1).
%if $[p,q]$ is shortest then $p\in int X_{\square^{i'}}$
%$q\in int X_{\square^{j'}}$ and perpendicular....
Statement (ii) follows from Lemma 10.2, item (4) at the beginning
of Section 8, and the following calculation
for $l=i,j$
(see 6.2 for the definition of $s\0{n,l}$ and $r\0{n,l}$)
\vspace{.05in}

\hspace{1.4in}{\small $
\frac{sinh\, s\0{n,l}}{sinh\, r\0{n,k}}\,\,=\,\, \frac{\Big(\frac{sinh\, r\,\, sin\,\beta\0{l}}
{sin\,\alpha\0{n-2}}\Big)}{\Big(\frac{sinh\, r}{sin\,\alpha\0{n-k-3}} \Big)  }
\,\,=\,\, c\,\varsigma^{l-k}\,\,\leq c\,\varsigma\,\,<\,\,e^{-4}
\,\,<\,\,\frac{\sqrt{2}}{2}
$}\vspace{.05in}

Statement (iii) follows from (ii) and the definition of the sets $\cZ$.
We next prove (iv). Write $Z=\cZ(X_{\square^j})\cap\cZ(X_{\square^k})$.
By the definition of the sets $\cZ$ we have\vspace{.05in}

\hspace{.5in}{\small $
\begin{array}{ccl}Z&=&
\stackrel{\circ}{\sN}\0{s\0{n,j}}(X_{\square^j})\,\,\cap\,\stackrel{\circ}{\sN}\0{s\0{n,k}}(X_{\square^k})\,\,-\,\, \bigcup\0{l<k}\sN\0{r\0{n,l}}(X_{\square^l})\\
&\sbs&\stackrel{\circ}{\sN}\0{s\0{n,j}}(X_{\square^j})\,\,\cap\,\stackrel{\circ}{\sN}\0{s\0{n,k}}(X_{\square^k})\,\,-\,\, \bigcup\0{j\leq l<k}\sN\0{r\0{n,l}}(X_{\square^l})\\
&\sbs&\stackrel{\circ}{\sN}\0{s\0{n,j}}(X_{\square^j})\,\,\cap\,\stackrel{\circ}{\sN}\0{s\0{n,k}}(X_{\square^k})\,\,-\,\, \bigcup\0{j<l<k}\sN\0{r\0{n,l}}(X_{\square^l})
\,\,-\,\,\sN\0{r\0{n,j}}(X_{\square^j})\\
&\sbs&\stackrel{\circ}{\sN}\0{s\0{n,j}}(X_{\square^j})\,\,\cap\,\stackrel{\circ}{\sN}\0{s\0{n,k}}(X_{\square^k})\,\,-\,\, \bigcup\0{j<l<k}\Big(
 \sN\0{s\0{0}}(X_{\square^j}) \,\cap\,\sN\0{r\0{n,l}}(X_{\square^l})\,    \Big)
\,\,-\,\,\sN\0{r\0{n,j}}(X_{\square^j})
\end{array}$}\vspace{.05in}

\noindent This together with Lemma 10.3 imply
$Z\sbs \dX_{\square^j}\times A$ where\vspace{.05in}

\hspace{.3in}{\scriptsize $
A\,\,=\,\,\stackrel{\circ}{\sN}\0{s\0{n,k}}\Big( \rC \Delta(\square^k),\,
\stackrel{\circ}{\rC}\0{s\0{n,j}}\big(\sL(X_{\square^j})\big)\Big)
\,\,-\,\,\bigcup\0{j<l<k}\sN\0{r\0{n,l}}\Big( \rC \Delta(\square^l),\,
\rC\0{s\0{0}}\big(\sL(X_{\square^j})\big)\Big)\,\,-\,\,\B_{r\0{n,j}}\big(\rC\sL(X_{\square^j})\big)
$}%\vspace{.05in}

\noindent hence

\hspace{.02in}{\small $
A\,\,\sbs \,\,\stackrel{\circ}{\sN}\0{s\0{n,k}}\Big( \rC \Delta(\square^k),\,
\rC\sL(X_{\square^j})\Big)
\,\,-\,\,\bigcup\0{j<l<k}\sN\0{r\0{n,l}}\Big( \rC \Delta(\square^l),\,
\rC\sL(X_{\square^j})\Big)\,\,-\,\,\B_{r\0{n,j}}\big(\rC\sL(X_{\square^j})\big)
$}\vspace{.05in}

\noindent But for $i>j$ we have $s\0{n,i}=s\0{n-j,i-j}$, $r\0{n,i}=r\0{n-j-1,i-j-1}$
 and $r\0{n,j}=r\0{n-j-3}$ (see definitions in 6.2).
Therefore $A\sbs \cX\big(\rC\sL(X_{\square^j}) ,\, \Delta(\square^k),\,r\big)$. 
This proves (iv). The proof of (v) is similar to the proof of (iv) with
minor changes.  This proves Lemma 10.4.\vspace{.1in}

We now smooth the metric $\sigma\0{K_X}$. 
For each $\square\in K$ using the construction in Section 8 we get a Riemannian metric
$\cG\Big(\sL(X_\square), \cL_\square, h\0{\square}, r, \xi, \s{d}, (c, \varsigma)\Big)$
on $\sL(X_\square)$,
which we shall simply denote by $\cG\big(\sL(X_\square)\big)$. 
Define  the Riemannian metric $\cG(X_\square)$ on $\stackrel{\circ}{\s{N}}_{s\0{0}}\big(\dX_\square\big)$ by \vspace{.1in}

\hspace{2in}{\small $
\cG\big( X_\square\big)\,\,=\,\,\cE\0{\dX_\square}\Big(\, \cG\big(\sL(X_\square\big)  \, \Big) 
$}%\vspace{.1in}

\noindent {\bf Remark.} Recall that we can consider $\stackrel{\circ}{\s{N}}_{s\0{0}}\big(\dX_\square\big)$ contained in the infinite cone
$ \dX_\square\times\rC\sL(X_\square)$ (see remark 2 before 10.2), and note that
the definition of $\cG(X_\square)$ makes sense in the whole of
$ \dX_\square\times\rC\sL(X_\square)$.\vspace{.1in}

\noindent {\bf Proposition 10.5.} {\it The Riemannian metrics 
$\cG(X_{\square^j})$ and $\cG(X_{\square^k})$ coincide on the intersection
$\cZ(X_{\square^j})\cap\cZ(X_{\square^k})$, $i,j<n-1$.
Also the Riemannian metric $\cG(X_{\square^k})$ coincides with $\sigma\0{K\0{X}}$ on  $\cZ\cap\cZ(X_{\square^k})$.}\vspace{.1in}

\noindent {\bf Proof.} For the first statement items (i) and (iii) of Lemma 10.4 imply that we only need to
consider the case $\square^j\sbs\square^k$, $j<k<n-1$. 
By item (iv) of Lemma 10.4 it is enough to prove that 
$\cG(X_{\square^j})$ and $\cG(X_{\square^k})$ coincide on
$\dX_{\square^j}\,\times\,\cX\Big(\rC\sL(X_{\square^j}) ,\, \Delta(\square^k),\,r\Big)$
(see remark above). Property {\bf P6} in 8.2 implies that the metric $\cG(X_{\square^j})$
coincides with the metric\vspace{.05in}

\hspace{1.55in}{\small $\cE_{\dX_{\square^j}}\bigg[\cE_{\rC\Delta(\square^k)}\bigg( \,\,\cG\Big[\,\,\,\sL\Big(\,\Delta(\square^k)\,,\,\sL(X_{\square^j})\,\Big)\,\,\,\Big]\, \,  \bigg)  \bigg]$}\vspace{.05in}

\noindent on
$\dX_{\square^j}\,\times\,\cX\Big(\rC\sL(X_{\square^j}) ,\, \Delta(\square^k),\,r\Big)$.
But\vspace{.05in} 

\hspace{.75in}$\sL\big(\,\Delta(\square^k),\sL(X_{\square^j})\,\big)=
\sL\big(\,\Delta(\square^k),\sL(\square^j,K\,)\,\big)=\sL(\square^k,K)=\sL(X_{\square^k})$\vspace{.05in}

\noindent Hence we have to prove that \,{\small $\cE_{\dX_{\square^j}}\big(\cE_{\rC\Delta(\square^k)}(  g)\,   \big) =\cE_{\dX_{\square^k}}(  g)$}, where $g= \cG\big(\sL(X_{\square^k}))$. This follows 
from applying Proposition 2.5 locally. 
To prove the second statement in  10.5,
using a similar argument as above (with {\bf P7} instead of {\bf P6})
we reduce the problem to showing that on $\dX_{\square^k}\,\times\,\rC\sL(X_{\square^k})$
we have
$\cE_{\dX\0{\square^k}}(\sigma\0{\rC\sL(X\0{\square^k})})=\sigma\0{K\0{X}}$.
And this follows from applying 6.1.4 locally. This proves 
the proposition.\vspace{.1in}

Finally define the metric $\cG(K_X)=\cG\big(K_X,\cL,r,\xi,\s{d}, (c,\varsigma)\big)$
to be equal to $\cG(X_{\square^k})$ on $\cZ(X_{\square^k})$, for $\square^k\in K$,
$k<n-1$. And equal to $\sigma\0{K\0{X}}$ on $\cZ$.
By lemmas 10.2 and 10.5 the metric $\cG(K_X)$ is a well defined Riemannian metric
on the smooth manifold $(K_X, \cS\0{K_X})$.\vspace{.1in}

\noindent {\bf Corollary 10.6.} {\it Let $\epsilon>0$ and $M^n$  closed.
Choose $\xi$, $c$, $\varsigma$ satisfying (i) and (ii) in 8.4.2, and $\xi\geq n$.
Then the metric
$\cG(K_X)$ has all sectional curvatures $\epsilon$-pinched to -1,
provided  $d_i$, $r-d_i$, $i=2,...,n$, are sufficiently large.}\vspace{.1in}

\noindent {\bf Proof.} 
Choose $\epsilon'$, as in Remark 1.2(2) and 1.4(3), so that
a $(B_a,\epsilon')$-close to hyperbolic metric with charts of excess 1 has
sectional curvatures $\epsilon$-pinched to -1.
Take $A$ so that $A\geq C(n,k,\xi)$ (see 2.6), for all $k\leq n-1$.
Since $M$ is compact we only have finitely many 
cubes in a cubulation $K$ of $M$. Hence the set of links of $K$ (hence of $K_X$)
is finite. This together with Proposition 8.4.2 imply that all 
$\cG(\sL(\square^k,K),r,\s{d})$, are $(B_{a\0{j}},\frac{\epsilon'}{A})$-close to hyperbolic (here $j=n-k-1$ and $a\0{j}=r\0{j-2}-d\0{k+1}$),
with charts of excess $\xi_{j}$. All this provided $d_i$, $r-d_i$, $i=2,...,n$, are sufficiently large.
We can apply Theorem 2.6 (locally, see remark below) to get that
the metrics $\cG(X_{\square^k})$ are $(B_{a\0{j}},\epsilon')$-close to hyperbolic on $\cZ(X_{\square^k})$, with charts of excess $\xi_j-1$,
provided $d_i$, $r-d_i$, $i=2,...,n$, are sufficiently large.
Here $\xi_j-1\geq 1$
and 1.4(4) imply that we can take the excess to be 1.
Therefore all the $\cG(X_{\square^k})$ have curvatures $\epsilon$-close to -1. This proves the corollary.\vspace{.1in}

The corollary proves (i) of the main Theorem. Items (ii), (iii) follow from \cite{ChD}.
Item (iv) follows from  Proposition 9.4.
This proves the Main Theorem.\vspace{.1in}

\noindent {\bf Remark.} Note that it does not make sense to say that
$\cG(X_{\square^k})$ is $\epsilon'$-close to hyperbolic because neither $\dX_{\square^k}$
nor $\dX_{\square^k}\times \rC\sL(X_{\square^k})$ have a center. What we mean by
the ``local application of Theorem 2.6" mentioned in the proof above is the following.
Take $p\in \cZ(X_{\square^k})$ and let $B\sbs\dX_{\square^k}$ be an open ball
centered at $p$. Note that we can also consider $B\times\rC\sL(X_{\square^k})\sbs
\HH^{n-k}\times\rC\sL(X_{\square^k})=\cE_k(\rC\sL(X_{\square^k}))$ and we can now apply 2.6 to $\cE_k(\rC\sL(X_{\square^k}))$, where we are considering $p$
as the center.

\vspace{.2in}

\begin{center} {\bf \large Section 11. Proof of Theorem A.}
\end{center}

Let $N$ be a closed smooth manifold that bounds a compact smooth manifold
$M^m$. Denote the given smooth structure of $N$ by $\cS_N$.
Let $Q$ be the smooth $m$-manifold with one point singularity
formed by gluing the cone $\rC_1 N$ to $M$ along $N\sbs M$.
Let $q$ be the singularity of $Q$ and note that it is modeled on $\rC N$
(see 7.3). A triangulation of $Q$ is obtained by coning a smooth triangulation
of the manifold with boundary $M$, and let $f:K\ra Q$ be the induced
cubulation (see appendix G). Write $f^{-1}(q)=p$. Note that $(K,f)$ is a smooth cubulation of $Q$ in the sense of Section 7.3. By item (2) of 7.3 we have that
$Q-\{q\}$ has a a normal smooth structure $\cS'$ for $K$, induced
by a set of links smoothings $\cL$.\vspace{.1in} 

Let $K_X$ be the Charney-Davis strict hyperbolization of $K$. Also denote by
$p$ the singularity of $K_X$. By item (1) of 9.4, the space $K_X-\{p\}$
has a normal smooth atlas $\{H_\square \}_{\square\in K}$ and normal 
smooth structure $\cS_{K_X}$. Moreover, since
we are assuming $Wh(N)=0$ (if $dim\, N>4$) we have that
we can take the domain $\rC N-\{o\0{\rC N}\}=N\times (0,1]$ of $H_p$
with product smooth structure $\cS_N\times\cS_{(0,1]}$ (see 9.4).\vspace{.1in}

We can now proceed exactly as in Section 10 and define the sets
$\cZ(X_\square)$, $\cZ$, and the metrics $\cG(X_\square)$ depending on
$\cL,r,\xi,\s{d}, (c,\varsigma)$. For the special case $\square^0=p$
we use the results in Section 8.5. We obtain in this way a
Riemannian metric $\cG(K_X)=\cG(K_X,\cL,r,\xi,\s{d}, (c,\varsigma))$ on $K_X-\{p\}$.
Theorem A and its addendum now follow from 8.5.1 (iii), (iv) and the result of
Belegradek and Kapovitch \cite{BK} mentioned in the introduction (before the addendum
to Theorem A). To be able to apply 8.5.1 we need to satisfy the hypothesis
made at the beginning of 8.5: that the Whitehead group $Wh(\pi_1N)$ vanishes.
But this follows from \cite{FH}. This proves Theorem A.%\vspace{.1in}

Pedro Ontaneda

SUNY, Binghamton, N.Y., 13902, U.S.A.

\end{document}